\documentclass[11pt]{amsart}
\usepackage{amsfonts,amsmath,amsthm,amssymb,wasysym,rotating,mathrsfs}
\usepackage{tikz,tikz-cd}
\usepackage{graphicx,enumitem}

\usepackage[all]{xy}
\usepackage{hyperref}

\usepackage[margin=1in]{geometry}

\usepackage[capitalize,nameinlink,noabbrev,nosort]{cleveref}

\graphicspath{{images/}{../images/}}

\definecolor{darkred}{rgb}{1,0,0}

\newtheorem{maintheorem}{Theorem}	
 
\crefname{maintheorem}{Theorem}{Theorems}
\newtheorem{theorem}{Theorem}[section]
\newtheorem{thm}[theorem]{Theorem}
\newtheorem{lemma}[theorem]{Lemma}
\newtheorem{lem}[theorem]{Lemma}
\crefname{lem}{Lemma}{Lemmas}
\newtheorem{proposition}[theorem]{Proposition}
\newtheorem{prop}[theorem]{Proposition}

\newtheorem{cor}[theorem]{Corollary}
\newtheorem{conj}[theorem]{Conjecture}
\crefname{conj}{Conjecture}{Conjectures}

\theoremstyle{remark}

\theoremstyle{definition}
\newtheorem{rmk}[theorem]{Remark}
\newtheorem{example}[theorem]{Example}
\newtheorem{defn}[theorem]{Definition}
\crefname{defn}{Definition}{Definitions}

\crefname{question}{Question}{Questions}
\newtheorem*{claim*}{Claim}
\newtheorem*{defn*}{Definition}
\newtheorem*{thm*}{Theorem}

\newcommand{\bbc}{\mathbb{C}}

\newcommand{\bbf}{\mathbb{F}}

\newcommand{\bbz}{\mathbb{Z}}

\newcommand{\mca}{\mathcal{A}}

\newcommand{\mcc}{\mathcal{C}}

\newcommand{\mcf}{\mathcal{F}}
\newcommand{\mcg}{\mathcal{G}}

\newcommand{\mcr}{\mathcal{R}}
\newcommand{\mcs}{\mathcal{S}}

\newcommand{\mcu}{\mathcal{U}}

\newcommand{\lam}{\lambda}

\DeclareMathOperator{\Piv}{minVal}
\DeclareMathOperator{\Ing}{Ing}
\DeclareMathOperator{\Lec}{Lec}
\DeclareMathOperator{\Lus}{Lus}
\DeclareMathOperator{\obj}{obj}
\DeclareMathOperator{\add}{add}
\DeclareMathOperator{\Spr}{Spr}

\newcommand{\bw}{\mathbf{w}}

\newcommand{\bv}{\mathbf{v}}
\newcommand{\bx}{\mathbf{x}}
\newcommand{\bA}{\mathbf{A}}
\newcommand{\bB}{\mathbf{B}}
\newcommand{\pth}[1]{L_{#1}}
\newcommand{\young}[1]{\lambda_{#1}}
\def\skew_#1^#2{#1/#2}
 
\newcommand{\rmin}[1]{\Delta^{\rho}_{#1}} 
\newcommand{\lmin}[1]{\Delta^{\lambda}_{#1}} 
\newcommand{\ch}{\chi} 
\newcommand{\wl}[1]{w_{(#1)}}  
\newcommand{\wu}[1]{w^{(#1)}} 
\newcommand{\vl}[1]{v_{(#1)}} 
\newcommand{\vu}[1]{v^{(#1)}}

\newcommand{\jhol}[1]{J^{+}_{\mathbf{#1}}}
\newcommand{\jsol}[1]{J^\bullet_{\mathbf{#1}}} 
\newcommand{\wiring}[1]{W_{#1}} 
\newcommand{\jc}[1]{\Spr(#1)}
\newcommand{\iquiv}[1]{Q_{#1}^{\Ing}}
\newcommand{\lquiv}[1]{Q_{#1}^{\Lec}}
\newcommand{\wquiv}[1]{Q^W_{#1}}
\newcommand{\iclus}[1]{\mathbf{A}_{#1}}
\newcommand{\lclus}[1]{\mathbf{B}_{#1}}
\newcommand{\ivar}[1]{A_{#1}}
\newcommand{\lvar}[1]{B_{#1}}
\newcommand{\iseed}[1]{\Sigma_{#1}^{\Ing}}
\newcommand{\lseed}[1]{\Sigma_{#1}^{\Lec}}
\newcommand{\vw}{v, \bw}
\newcommand{\rich}[1]{\mathring{\mcr}_{#1}}
\def\Fl{\mcf\ell}

\newcommand{\minor}[1]{\Delta_{#1}}

\def\gridn{\mcg_n}

\usepackage{subfiles}

\title{Leclerc's conjecture on a cluster structure for type A Richardson varieties}
\author{Khrystyna Serhiyenko}
\author{Melissa Sherman-Bennett}
\thanks{KS and MSB were supported by the National Science Foundation under Award No.~DMS-2054255 and Award No.~DMS-2103282 respectively. Any opinions, findings, and conclusions or recommendations expressed in this material are
	those of the authors and do not necessarily reflect the views of the National Science
	Foundation.}

\begin{document}
	
\maketitle

\begin{abstract}Leclerc \cite{Leclerc} constructed a conjectural cluster structure on Richardson varieties in simply laced types using cluster categories.
We show that in type A, his conjectural cluster structure is in fact a cluster structure. We do this by comparing Leclerc's construction with another cluster structure on type A Richardson varieties due to Ingermanson \cite{gracie}. Ingermanson's construction uses the combinatorics of wiring diagrams and the Deodhar stratification. Though the two cluster structures are defined very differently, we show that the quivers coincide and clusters are related by the twist map for Richardson varieties, recently defined by Galashin--Lam \cite{GLRichardson}.
\end{abstract}

	\section{Introduction}
	
	In this article, we consider cluster structures on open Richardson varieties $\rich{v, w}$ in the complete flag variety $\Fl_n$. For $v \leq w \in S_n$, the open Richardson variety $\rich{v,w}$ is the intersection of the Schubert cell $C_w$ with the opposite Schubert cell $C^v$, and is smooth, affine, and irreducible. Open Richardson varieties are related to the geometric interpretation of Kazhdan-Lusztig polynomials \cite{Deodhar,KL79,KL80}; open Richardson varieties also arise in total positivity for $\Fl_n$ \cite{LusTPGr,Rietsch}. Special cases include \emph{(open) positroid varieties} \cite{KLS,Postnikov}, which are Richardson varieties $\rich{v,w}$ where $w$ has a single descent.
	
	A cluster structure on $\rich{v,w}$ is an identification of the coordinate ring $\bbc[\rich{v,w}]$ with a cluster algebra $\mca(\Sigma)$. Cluster algebras were introduced by Fomin--Zelevinsky \cite{CA1} to provide an algebraic and combinatorial framework for Lusztig's dual canonical bases and total positivity \cite{LusCanonicalBases,LusTPGr}. Cluster algebras have appeared in a wide range of fields, including Teichm\"uller theory \cite{FGDualityConj}, mirror symmetry \cite{GHKK}, Poisson geometry \cite{GSV}, symplectic geometry \cite{STWZ}, knot theory \cite{FPST}, and scattering amplitudes in high energy physics \cite{GGSVV}. One major direction of research is to understand when varieties naturally ocurring in representation theory have a cluster structure; examples of varieties with cluster structures include Grassmannians \cite{Scott}, double Bruhat cells in semisimple Lie groups \cite{CA3}, and unipotent cells in Kac-Moody groups \cite{GLS}.
	
	Cluster algebras are commutative rings with a distinguished set of generators called \emph{cluster variables}, which are grouped together into \emph{clusters}. A cluster can be \emph{mutated} into another cluster, and any two clusters are related by a sequence of mutations. The information of all mutations of a cluster is encoded in a \emph{quiver}, i.e. a directed graph. A cluster and its quiver together form a \emph{seed}.

	Leclerc \cite{Leclerc} used categorification to construct a conjectural cluster structure\footnote{Leclerc's results and conjecture are in types ADE. We will deal only with the type A case in this paper.} for $\rich{v, w}$. In particular, he defined a cluster category inside the module category of a preprojective algebra and identified certain cluster tilting objects in this category. Each cluster tiliting object naturally gives rise to a seed $\lseed{v,\bw}=(\bB_{\vw}, \lquiv{\vw})$, where the cluster $\bB_{\vw}$ is obtained via a cluster character map and the quiver $\lquiv{v, \bw}$ records irreducible morphisms between indecomposable summands of the cluster tilting module. In this contruction, it is relatively easy to obtain the cluster, but quite difficult to compute the quiver. 
Leclerc showed that the cluster algebra $\mca(\lseed{v, \bw})$ is a subring of $\bbc[\rich{v,w}]$. 

	Leclerc conjectured that $\mca(\lseed{v, \bw})$ is equal to $\bbc[\rich{v, w}]$, and showed that equality holds in some special cases. One of the obstacles in proving Leclerc's conjecture is the difficulty in computing the quiver $\lquiv{v,\bw}$ in general. In \cite{SSBW}, the authors together with L. Williams showed that for certain positroid varieties, the quiver $\lquiv{v, \bw}$ coincides with a \emph{plabic graph} quiver. Later, Galashin--Lam \cite{GLPositroid} extended this result to all positroid varieties, and used this to show Leclerc's conjecture in the positroid variety case.
	
Our main result is a proof of Leclerc's conjecture in type A.
	
	\begin{maintheorem}\label{thm:LecConjIsTrue}
		Let $v \leq w$ and let $\bw$ be a reduced expression for $w$. Then 
		\[\bbc[\rich{v, w}]=\mca(\lseed{\vw}). \]
		Moreover, the cluster algebra $\mca(\lseed{\vw})$ does not depend on the choice of $\bw$.
	\end{maintheorem}

We prove \cref{thm:LecConjIsTrue} by comparing $\mca(\lseed{\vw})$ with another cluster structure on $\rich{v,w}$, defined by Ingermanson \cite{gracie}. Ingermanson constructed a seed $\iseed{v, \bw}=(\bA_{v, \bw}, \iquiv{\vw})$ using the wiring diagram of a \emph{unipeak} expression for $w$. The cluster variables $\bA_{v, \bw}$ are particular factors of the \emph{chamber minors} of Marsh--Rietsch \cite{MR}, and the quiver $\iquiv{\vw}$ can be read off from the wiring diagram. In her construction, determining the factorization of chamber minors into cluster variables is quite involved, but once this has been done, it is easy to write down the quiver. Ingermanson showed that the upper cluster algebra $\mcu(\iseed{v, \bw})$ is equal to $\bbc[\rich{v,w}]$; recent results of \cite{GLSBS1} imply that Ingermanson's quiver is locally acyclic, and so $\mca(\iseed{v, \bw})= \bbc[\rich{v,w}]$. 

We show the following relationship between Ingermanson's seed and Leclerc's.

\begin{maintheorem}\label{thm:compareStructures}
	Let $v \leq w$ and let $\bw$ be a unipeak expression for $w$. Let $\tau_{v,w}$ be the right twist map for $\rich{v,w}$ from \cite{GLRichardson}. Then 
	\[(\bA_{v, \bw}, \iquiv{\vw})= (\bB_{\vw}\circ \tau_{v,w}, \lquiv{\vw}).\]
\end{maintheorem}

We separately show that all of the seeds Leclerc defines are related by mutations. As a corollary of these result, we show that the positive part of $\rich{v,w}$ defined by Leclerc's cluster structure agrees with the totally positive part $\rich{v, w}^{>0, \Lus}$ defined by Lusztig \cite{LusTPGr}.

\begin{maintheorem}\label{thm:posPartLeclerc}
	Let $v \leq w$. The subset 
	\[\rich{v, w}^{>0, \Lec}:=\{F \in \rich{v, w}: \text{ all cluster variables in }\mca{(\lseed{\vw})}\text{ are positive on }F\}\]
	coincides with $\rich{v,w}^{>0, \Lus}$.
\end{maintheorem}

\cref{thm:compareStructures} has the effect of simplifying the definitions of both Leclerc's and Ingermanson's seeds. We obtain a much more straightforward method to factor chamber minors into Ingermanson's cluster variables, and an elementary method to draw Leclerc's quiver from a wiring diagram. This alternate description of Leclerc's quiver is in the same vein as the descriptions in the positroid variety case given by \cite{SSBW,GLPositroid}. As such, we hope our results make cluster structures on Richardson varieties more accesible.

Our results show that Leclerc's and Ingermanson's cluster structure on $\rich{v,w}$ are related by the twist map for Richardsons, which generalizes the twist map for positroid varieties \cite{MSTwist}. In the positroid variety case, the twist map is conjectured to be a \emph{quasi-cluster transformation} \cite{Fraser}; i.e. the twist map is conjectured to be a sequence of mutations followed by rescaling by Laurent monomials in frozens. We make the same conjecture in the Richardson variety case.

\begin{conj}
	The twist map $\tau_{v,w}: \rich{v, w} \to \rich{v, w}$ is a quasi-cluster transformation. As a result, any cluster of $\mca(\iseed{v, \bw})$ is related to any cluster of $\mca(\lseed{v, \bw})$ by a sequence of mutations and rescaling by Laurent monomials in frozens.
\end{conj}

We briefly discuss related work on Richardson variety cluster structures, which has been a very active topic of late.  In types ADE, M\'enard \cite{Menard} constructed another cluster tilting object for each reduced word $\bw$ of $w$; his construction has the advantage that the quiver is constructed algorithmically. Cao--Keller \cite{CK} recently showed that, if $\Sigma^M_{v, \bw}$ is the seed obtained from M\'enard's cluster tilting object via the cluster character map, then $\mcu(\Sigma^M_{\vw})=\bbc[\rich{v, w}]$ (again, in types ADE). The relation between $\Sigma_{v, \bw}^M$ and $\lseed{\vw}$ is as yet unknown; the quivers are conjectured to be mutation-equivalent. In a separate direction,  Casals--Gorsky--Gorsky--Le--Shen--Simental \cite{CGGLSS} and  Galashin--Lam--SB--Speyer \cite{GLSBS1,GLSBS2} have independently given cluster structures on \emph{braid varieties} in arbitrary type, which generalize Richardson varieties; it is not known how these two cluster structures are related. Ingermanson's construction is a special case of the construction in \cite{GLSBS1}; M\'enard's seeds are special cases of those in \cite{CGGLSS}.

We begin with background on Richardson varieties and various combinatorial constructions in \cref{sec:background}. We review the constructions of $\lseed{v, \bw}$ and $\iseed{v, \bw}$ in \cref{sec:background2}. In \cref{sec:clusterVar}, we give the relationship between Leclerc's cluster variables and Ingermanson's. In \cref{sec:LQuivFromWiring} we describe Leclerc's quiver using wiring diagrams, and in \cref{sec:quiversEqual} we use this description to prove that Leclerc's quiver coincides with Ingermanson's quiver. Finally, in \cref{sec:LSeedsMutationEquiv}, we complete the proofs of \cref{thm:LecConjIsTrue,thm:compareStructures,thm:posPartLeclerc}.

\vspace{12pt}

\noindent\textbf{Acknowledgements:} MSB thanks Pavel Galashin, Thomas Lam, and David Speyer for helpful conversations related to Ingermanson's construction.

	\section{Background} \label{sec:background}
	
	We use the following standard combinatorial notation: $[n]:=\{1, \dots, n\}$, $\binom{[n]}{h}$ is the set of cardinality $h$ subsets of $[n]$, $w_0$ is the longest permtuation of $S_n$, $s_i \in S_n$ is the transposition exchanging $i$ and $i+1$, for $v, w\in S_n$, $v \leq w$ means $v$ is less than $w$ in the Bruhat order and $\ell(w)$ denotes the length of $w$.
	
	\subsection{Background on Richardson varieties}
	
	Let $G=SL_n(\bbc)$ and let $B, B_- \subset G$ denote the Borel subgroups of upper and lower triangular matrices, respectively. Let $N, N_-$ denote the corresponding unipotent subgroups of upper and lower unitriangular matrices, respectively. For $g \in G$, let $g_i$ denote the $i$th column of $g$. We denote the minor of $g$ on rows $R$ and columns $C$ by $\minor{R, C}(g)$.
	
	For $w \in S_n$, we choose a distinguished lift $\dot{w}$ of $w$ to $G$. The lift satisfies
	\[\dot{w}_{ij}=\begin{cases}
		\pm 1 &\text{ if } i=w(j)\\
		0 & \text{ else }
	\end{cases}\]
and the signs of entries are determined by the condition that $\minor{w[j], [j]}(\dot{w})=1$ for all $j \in [n]$. If the particular lift of $w$ to $G$ does not matter, we also write $w$ for the lift (e.g. we write $B w B$ rather than $B \dot{w}B$).
	
	We identify the flag variety $\Fl_n$ with the quotient $G/B$. Concretely, a matrix $g \in G$ represents the flag $V_\bullet=(V_1 \subset V_2 \subset \cdots \subset V_n=\bbc^n)$ where $V_i$ is the span of $g_1, \dots, g_i$. 
The flag variety has two well-known decompositions into cells, the Schubert decomposition
	\[G/B = \bigsqcup_{w \in S_n} B w B/B = \bigsqcup_{w \in S_n} C_w\]
	and the opposite Schubert decomposition
		\[G/B = \bigsqcup_{w \in S_n} B_- w B/B= \bigsqcup_{w \in S_n} C^w.\]
	The stratum $C_w$ is a \emph{Schubert cell} and is isomorphic to $\bbc^{\ell(w)}$. The stratum $C^w$ is an \emph{opposite Schubert cell} and is isomorphic to $\bbc^{\ell(w_0)-\ell(w)}$. For a fixed lift $w$, it is well-known that the projection map $G \to G/B$ restricts to isomorphisms
	\begin{equation}\label{eq:schub-iso}
		N w\cap w N_- \xrightarrow{\sim} C_w \qquad \qquad N_-w \cap wN_- \xrightarrow{\sim} C^w.
	\end{equation}
	Or, more concretely, each coset in $C_w$ (resp. $C^w$) has a unique representative matrix which differs from $w$ only in entries that lie both above and to the left (resp. both below and to the left) of a nonzero entry of $w$ (see e.g. \cite{Fulton}).
	
	
We are concerned with the intersection of an opposite Schubert cell and a Schubert cell
\[\rich{v, w}:= C^v \cap C_w \]
which is called an \emph{(open) Richardson variety}. We usually drop the adjective ``open." The Richardson variety $\rich{v,w}$ is nonempty if and only if $v \leq w$, in which case it is a smooth irreducible affine variety of dimension $\ell(w)-\ell(v)$ \cite{Deodhar}. 

Open Richardson varieties were studied in the context of Kazdhan-Lusztig polynomials \cite{KL79}; the number of $\bbf_q$ points of $\rich{v,w}$ is exactly the $R$-polynomial indexed by $(v, w)$ \cite{Deodhar}, which can be used to recursively compute Kazhdan-Lusztig polynomials. The $\bbf_q$-point counts and more generally the cohomology of $\rich{v,w}$ are also related to knot homology \cite{GLCohom}. Real points of $\rich{v,w}$, and in particular \emph{positive points}, feature in work of Lusztig and Rietsch \cite{LusTPGr,Rietsch} on total positivity. Special cases of Richardson varieties include the \emph{(open) positroid varieties} of \cite{KLS}, which are Richardson varieties $\rich{v,w}$ where $w$ has a single descent. Richardson varieties themselves are special cases of \emph{braid varieties} (see e.g. \cite{CGGS}).

We identify $\rich{v, w}$ with two different subsets of $G$, one for Ingermanson's construction and one for Leclerc's. We will later use these identifications to define functions on $\rich{v, w}$. Below, we use the involutive automorphism $g \mapsto g^\theta$ of $G$ from \cite[(1.11)]{FZ99}; the $(i,j)$ entry of $g^\theta$ is the minor of $g$ obtained by deleting the $i$th row and $j$th column. It is not hard to check that $B^{\theta}=B_-$, $N^{\theta}=N_-$, and $\dot{v}^\theta$ is another lift of $v$ to $G$.

\begin{lem}\label{lem:richIso}
	For $v\leq w$, let 
		\[N_{v,w}:= N \cap \dot{w}N_- \dot{w}^{-1} \cap B_- vB \dot{w}^{-1} \quad \text{and} \quad N'_{v,w}:= N \cap \dot{v}^{-1}N \dot{v} \cap \dot{v}^{-1} B_- w B_-. \]
	
	Also, let $D: N_{v, w} \to G$ be the renormalization map sending $g \to g d_g$, where $d_g$  is the unique diagonal matrix so that $\minor{v[j], w[j]}(gd_g)=1$ for all $j$.
	
	We have isomorphisms
	\begin{align*}
		\alpha: D(N_{v, w}) &\to \rich{v,w}  &\quad \beta: N'_{v, w} & \to \rich{v, w}\\
		g d_g &\mapsto g d_g \dot{w} B& \quad g &\mapsto (\dot{v} g)^{\theta} B.
	\end{align*}
\end{lem} 
\begin{proof}
	If $g \in N_{v, w}$, then $g \dot{w}$ is in $B_- v B$. In particular, the minors $\minor{v[j], [j]}(g \dot{w})= \minor{v[j], w[j]}(g)$ are nonzero. This implies the map $D$ is well-defined. It is also an isomorphism onto its image. 
	
	The map $\alpha$ can be written as a composition of two maps
	\begin{alignat*}{3}
		D(N_{v, w}) & \xrightarrow{D^{-1}}& N_{v, w} &\xrightarrow{\alpha'} \rich{v, w}\\
		g d_g& \longmapsto  &g & \longmapsto g \dot{w} B
	\end{alignat*}
since $g d_g \dot{w}B$ is equal to $g \dot{w} B$.

	Now, it follows easily from \eqref{eq:schub-iso} that $\alpha'$ and $\beta$ are both isomorphisms, noting in the first case that $g \dot{w}$ is in $N \dot{w}\cap \dot{w}N_- \cap B_- vB$ and in the second that $(\dot{v} g)^\theta$ is in $N_- \tilde{v} \cap \tilde{v}N_- \cap B w B$ where $\tilde{v}=\dot{v}^\theta$.
\end{proof}

\begin{rmk}\label{rmk:differencesLeclerc}
	Leclerc identifies the flag variety with $B_- \backslash G$ rather than $G/B$, and so considers the variety
	\[^- \rich{v, w}:=B_-\backslash {(B_- v B \cap B_- w B_-)}\] 
	which is different from, though isomorphic to, $\rich{v, w}$. We fix an isomorphism so that we can pullback functions on $^-\rich{v, w}$ to functions on $\rich{v,w}$. The isomorphism we choose is
	\[
  \rich{v, w}\xrightarrow{\Theta} {(B v B_- \cap B_- w B_-)}/B_- \xrightarrow{\delta_1} \dot{v} N'_{v, w} \xrightarrow{\delta_2} 	{^-\rich{v, w}}.
	\]
	The map $\Theta: g B\mapsto g^{\theta} B_-$ is induced by the involution $g \mapsto g^\theta$ on $G$; from \cite[Section 2]{FZ99}, one can see that $B^{\theta}=B_-$ and that $\dot{w}^\theta$ is another lift of $w$ so it is indeed an isomorphism. The maps $(\delta_1)^{-1}$ and $\delta_2$ are the natural projections from $\dot{v} N'_{v,w}$ to $G/B_-$ and $B_-\backslash G$ respectively; these are isomorphisms using the appropriate analogue of \eqref{eq:schub-iso}. The composition $\delta= \delta_2 \circ \delta_1$ is called the \emph{left chiral map} in \cite[Definition 2.2]{GLRichardson}.
\end{rmk}

\subsection{Background on wiring diagrams and chamber minors}

Before describing Ingermanson's and Leclerc's seeds, we need some combinatorial background.

Given $w\in S_n$, a \emph{reduced expression} for $w$ is an expression $\bw=s_{h_1} \dots s_{h_\ell}$ where $\ell$ is as small as possible. The number $\ell$ is the \emph{length} of $w$, denoted $\ell(w)$. We use the notation
\[\wl{i}:=s_{h_1} \dots s_{h_{i-1}} \quad \text{ and } \quad \wu{i}:=s_{h_\ell} \dots s_{h_{i}} = w^{-1}\wl{i}\]
for prefixes of $\bw$ and prefixes of $\bw^{-1}$, setting $\wl{1}=e$.

As a shorthand, we write $v \leq \bw$ to indicate a pair of permutations $v \leq w$ and a choice of reduced expression $\bw$ for $w$.

\begin{defn}\label{defn:PDS}
	Let $v \leq \bw=s_{h_1} \dots s_{h_\ell}$. A \emph{subexpression} for $v$ in $\bw$ is an expression for $v$ of the form $v=s_{h_1}^v \dots s_{h_\ell}^v$ where $s_{h_i}^v \in \{e, s_{h_i}\}$. 
	As for $\bw$, we define
	\[\vl{i}:=s_{h_1}^v \dots s_{h_{i-1}}^v \quad \text{ and } \quad \vu{i}:=s_{h_{\ell}}^v \dots s_{h_{i}}^v.\]
	
	The indices $i$ where $s_{h_i}^v \neq e$ is the \emph{support} of the subexpression. The subexpression is \emph{reduced} if the support has size $\ell(v)$. The \emph{positive distinguished subexpression} (PDS) for $v$ in $\bw$ is the reduced subexpression whose support is lexicographically largest. If $\bw$ is fixed, we denote the PDS for $v$ by $\bv$.
	
	We denote the support of the PDS by $\jhol{v}$, and call these the \emph{hollow crossings} of $\bw$. The complement of the support is $\jsol{v}$; we call these the \emph{solid crossings} of $\bw$. Note that $|\jsol{v}|= \ell(w)-\ell(v)=\dim \rich{v, w}$.
\end{defn}

\begin{example}
Let $\bw=s_1 s_2 s_1 s_3 s_2 s_1$ and let $v=3214$. Reduced subexpressions for $v$ in $\bw$ include
\[ e e s_1 e s_2 s_1 \qquad e s_2 s_1 e s_2 e \qquad s_1 s_2 e e e s_1 \qquad s_1s_2s_1eee.\]
The first subexpression has support $\{3, 5, 6\}$ and is the PDS for $v$ in $\bw$. So the hollow crossings are $\jhol{v}=\{3, 5, 6\}$ and the solid crossings are $\jsol{v}=\{1, 2, 4\}$.
\end{example}

\begin{rmk}
	Alternatively, the PDS for $v$ can be defined using a greedy procedure, moving from right to left. Set $\vl{\ell+1}=v$. If $\vl{i+1}$ is already determined, then $\vl{i}$ is equal to either $\vl{i+1}$ or $\vl{i+1}s_{h_{i}}$, whichever is smaller. In the first case, $s_{h_i}^v=e$; in the second, $s_{h_i}^v=s_{h_i}$.
\end{rmk}

\begin{rmk}
	The notion of positive distinguished subexpressions (and more generally, distinguished subexpressions) is due to Deodhar \cite{Deodhar}. Our notation for the support and complement of the support is inspired by \cite{MR}, as is the terminology ``solid" and ``hollow" crossing. The $+$ in the superscript of $\jhol{v}$ is to indicate that $\jhol{v}$ records where the length of $\vl{i}$ increases. The $\bullet$ in the superscript of $\jsol{v}$ is to remind the reader that these are the solid crossings.
\end{rmk}

For $v\leq \bw$, we will draw both the reduced expression $\bw$ and the PDS $\bv$ in the plane as wiring diagrams. Since $\bw$ is itself a positive distinguished subexpression of $\bw$, we make all definitions for the PDS $\bv$.

\begin{defn}
The \emph{wiring diagram} $\wiring{\bv}$ is obtained by replacing each simple transposition $s_{i}$ in $\bv$ with the configuration of strands on the left, and each $e$ in $\bv$ with the configuration of strands on the right.
\begin{center}
	\includegraphics[width=0.4\textwidth]{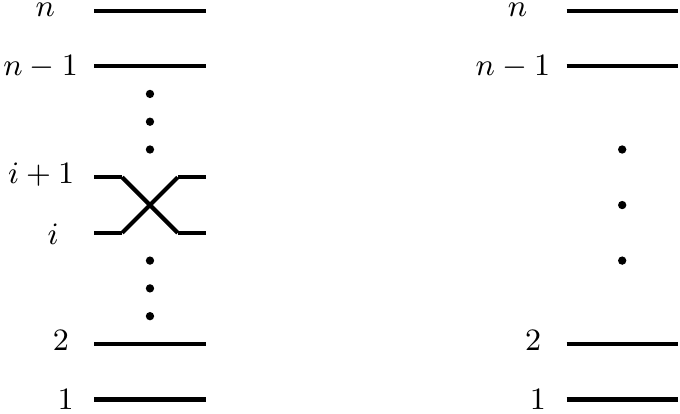}
\end{center}
	
	\noindent We label the crossings of $\wiring{\bv}$ with $\jsol{v}$ in the natural way. We also label the endpoints of the strands from $1$ to $n$, going from bottom to top. Each crossing $c$ has a \emph{rising strand}, whose height immediately to the right of $c$ is higher than immediately to the left of $c$, and a \emph{falling strand}.
\end{defn}
If a strand $\gamma$ in $\wiring{\bv}$ has right endpoint $h$, it has left endpoint $v(h)$. Since $\bv$ is reduced, then no two strands cross more than once.

A \emph{chamber} of a wiring diagram is a connected component of the complement of the strands. We denote chambers by $\chi$; the chamber to the left of crossing $c$ is $\chi_c$. We can label each chamber with a subset of $[n]$.

\begin{defn}[Right and left labeling of chambers]
	Let $\chi_c$ be a chamber of $\wiring{\bv}$. The \emph{left label} of $\chi_c$ is
	\[ \vl{c}[h_c] = \{i \in [n]: i \text{ is the left endpoint of a strand }\gamma \text{ below }\chi_c\}\]
	and the \emph{right label} is 
		\[ \vu{c}[h_c] = \{i \in [n]: i \text{ is the right endpoint of a strand }\gamma \text{ below }\chi_c\}.\]
\end{defn}

\begin{figure}
	\includegraphics[width=0.6\textwidth]{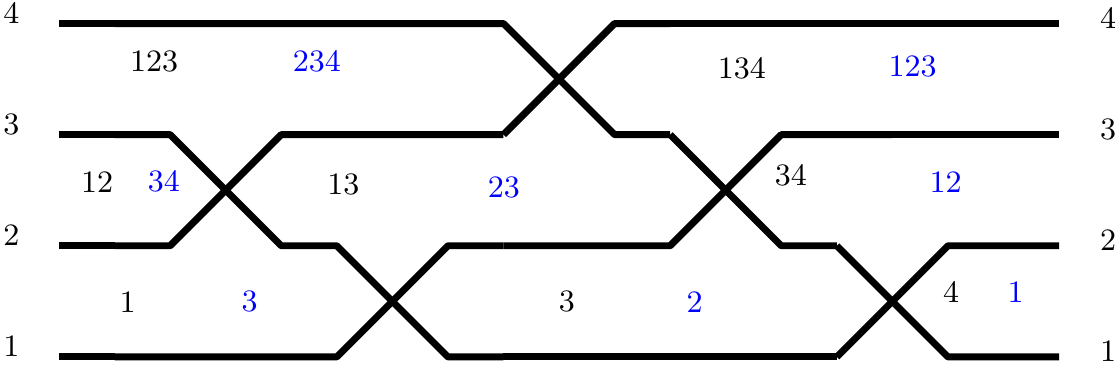}
	\caption{	\label{fig:wiringEx} A wiring diagram $\wiring{\bv}$ for $\bv=s_2 s_1 s_3 s_2 s_1$. The left and right labels of chambers are shown respectively in the left and right of each chamber.}
\end{figure}

The equalities above are easy to see by induction on $\ell$. Since $\bv$ is reduced, the left label can be obtained from the right label by applying $v$. See \cref{fig:wiringEx} for an example of a wiring diagram and its right and left labels.

The following combinatorial object will encode two seeds for $\rich{v, w}$, one in Ingermanson's cluster algebra and one in Leclerc's.

\begin{defn}
	Let $v \leq \bw$ and let $\bv$ be the PDS for $v$ in $\bw$. The \emph{stacked wiring diagram} $\wiring{v, \bw}$ is the union of the two wiring diagrams $\wiring{\bw}$ and $\wiring{\bv}$. We emphasize that the crossings of $\wiring{\bv}$ are drawn directly on top\footnote{This is in contrast to the ``double wiring diagrams" of \cite{FZ99}.} of the corresponding crossings of $\wiring{\bw}$. We call the strands of $\wiring{\bw}$ the \emph{$w$-strands} of $\wiring{v, \bw}$, and the strands of $\wiring{\bv}$ the \emph{$v$-strands}. We sometimes also call $w$-strands just ``strands". A \emph{chamber} of $\wiring{v, \bw}$ is a chamber of $\wiring{\bw}$. For $c \in [\ell]$, we denote by $\chi_c$ the chamber of $\wiring{v, \bw}$ which is to the left of crossing $c$. We call a chamber \emph{frozen} if it is open on the left.
\end{defn}

\begin{figure}
	\includegraphics[width=\textwidth]{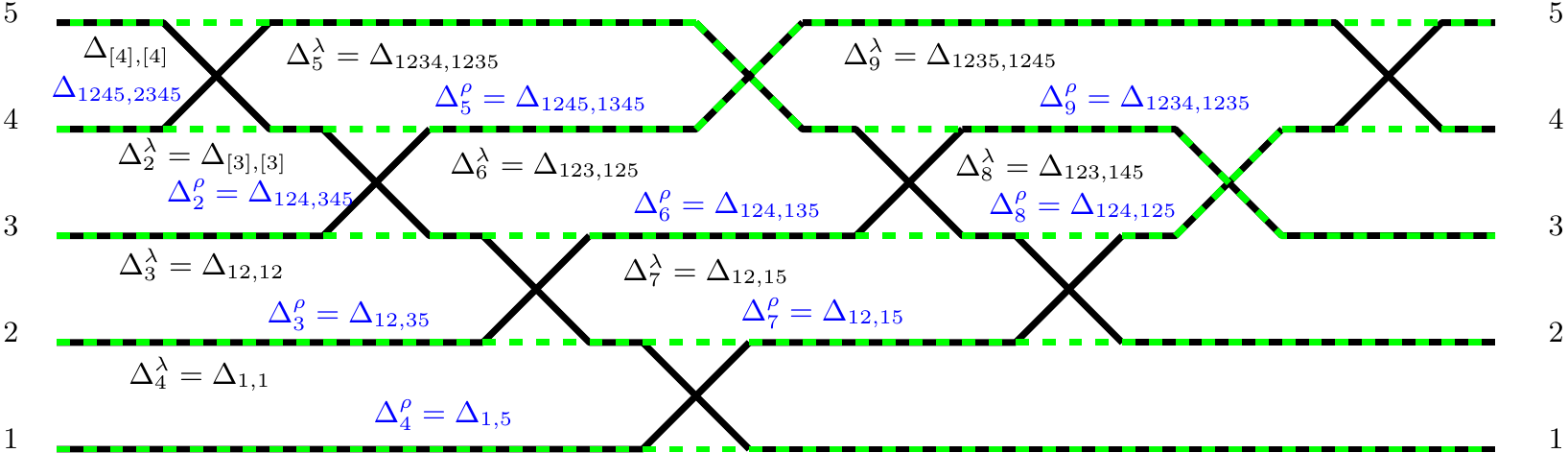}
	\caption{\label{fig:chamberMinorEx}The stacked wiring diagram $\wiring{v, \bw}$ for $\bw=s_4s_3s_2s_1s_4s_3s_2s_3s_4$ and $v=12534$. The $w$-strands are solid black; the $v$-strands are dashed green. Left and right chamber minors  (c.f. \cref{def:chamberMinor}) are black and blue, respectively.}
\end{figure}

See \cref{fig:chamberMinorEx} for an example of a stacked wiring diagram. Note that a crossing $c$ of $\wiring{v, \bw}$ is hollow (c.f. \cref{defn:PDS}) if it is a crossing of both $\wiring{\bv}$ and $\wiring{\bw}$; it is solid if it is only a crossing of $\wiring{\bw}$. 

Each chamber $\chi$ of $\wiring{v, \bw}$ has two left labels, one from the $w$-strands passing below $\chi$ and one from the $v$-strands. Similarly, $\chi$ has two right labels. We use these labels to define two regular functions for each chamber.

\begin{defn}\label{def:chamberMinor}
	Fix $v \leq \bw$ and stacked wiring diagram $\wiring{\vw}$. Let $c \in [\ell]$ and suppose $s_{h_c}=s_j$. For $gB\in \rich{v, w}$, let $x:= \alpha^{-1}(gB)$ and let $y := \beta^{-1}(gB)$, where $\alpha, \beta$ are the isomorphisms from \cref{lem:richIso}. The \emph{left chamber minor} of $\chi_c$ is the function
	\begin{align*}
		\lmin{c} : \rich{v, w}& \to \bbc\\
		gB &\mapsto \minor{\vl{c}[j], \wl{c}[j]} (x)
	\end{align*}
and the \emph{right chamber minor} is the function
	\begin{align*}
	\rmin{c} : \rich{v, w}& \to \bbc\\
	gB &\mapsto \minor{\vu{c}[j], \wu{c}[j]} (y).
\end{align*}
\end{defn}

See \cref{fig:chamberMinorEx} for an example of the left and right chamber minors. Note that chamber minors are defined only for chambers which are to the left of some crossing. This is because the analogous minors for the chambers which are open to the right are equal to $1$ on $\rich{v, w}$.

\begin{rmk}
	Left chamber minors were introduced by \cite{MR} in their study of the Deodhar stratification of $\rich{v, w}$ (in fact, their chamber minors were evaluated on elements of $N_{v,w}$ rather than elements of $D(N_{v,w})$, so they are monomially related to the left chamber minors defined here). They showed that the subset of $\rich{v, w}$ where the left chamber minors are nonzero is an algebraic torus, called the \emph{Deodhar torus}. This torus will be a cluster torus in Ingermanson's cluster structure. 
\end{rmk}

\begin{rmk}\label{rmk:our-minor-vs-leclerc}
	Leclerc uses functions $f_c:  {^-\rich{v, w}} \to \bbc$ defined by $f_c: B_- \dot{v} g \mapsto \minor{\vu{c}[h_c], \wu{c}[h_c]} (g)$, where $g \in N'_{v,w}$. This is related to the right chamber minor defined above by 
	\[\rmin{c}= f_c \circ \delta \circ \Theta\]
	where $\delta$ and $\Theta$ are as in \cref{rmk:differencesLeclerc}.
\end{rmk}

We will later use the left (resp. right) chamber minors to define cluster variables in Ingermanson's seed (resp. Leclerc's seed). The chamber minors are not algebraically independent; the chamber minors around a hollow crossing satisfy a binomial relation.

\begin{lem}\label{lem:hollowRel}
	Fix $v \leq \bw$ and a hollow crossing $c \in \jhol{v}$. Say the chambers surrounding $c$ are $\chi_{c^{\uparrow}}, \chi_{c^{\rightarrow}}, \chi_{c^{\downarrow}}, \chi_{c^{\leftarrow}}$. Then the chamber minors in those chambers satisfy
	\[\frac{\lmin{c^\uparrow} \lmin{c^\downarrow}}{\lmin{c^\leftarrow} \lmin{c^\rightarrow}} = 1 \quad \text{ and } \quad \frac{\rmin{c^\uparrow} \rmin{c^\downarrow}}{\rmin{c^\leftarrow} \rmin{c^\rightarrow}} =1.\]
\end{lem}
\begin{proof} 
	See the discussion in \cite{gracie} below Formula III.25 for the first relation; it follows from \cite{MR} and the Desnanot-Jacobi identity. 
	
	For the second relation, let $u:= \vl{c+1}$, $x:=\wl{c+1}$ and $i:=h_c$. Then the Desnanot-Jacobi relation implies 
	\[\minor{u[i+1], x[i+1]}\minor{u[i-1], x[i-1]}= \minor{u[i], x[i]}\minor{us_i[i], xs_i[i]}-\minor{us_i[i], x[i]}\minor{u[i], xs_i[i]}\]
	on $G$. Since $c \in \jhol{v}$ and $\bv$ is the rightmost subexpression for $v$ in $\bw$, we conclude that $us_i \nleq x$. Because $u \leq x$, this means that $us_i[i] \nleq x[i]$ in the Gale order (see \cref{defn:latticePath}). So $\minor{us_i[i], x[i]}$ vanishes identically on $B$. Noting that
	\[\minor{u[i+1], x[i+1]}= \rmin{c^\uparrow}, \quad \minor{u[i-1], x[i-1]}=\rmin{c^\downarrow}, \quad \minor{u[i], x[i]}=\rmin{c^\rightarrow}, \quad \minor{us_i[i], xs_i[i]}=\rmin{c^\leftarrow},\]
	this gives the desired relation.
\end{proof}

The left and right chamber minors are related by a \emph{twist automorphism} of $\rich{v, w}$, recently defined in \cite{GLRichardson}. The precise definition of the twist will not be needed, so we omit it.

\begin{prop}\cite[Theorem 11.6]{GLRichardson}\label{prop:twist}
	Fix $v \leq w$. There is a regular automorphism $\tau_{v,w}: \rich{v, w} \to \rich{v, w}$ such that for all $c \in [\ell]$, 
	\[\lmin{c}= \rmin{c} \circ \tau_{v,w}.\]
\end{prop}

\begin{proof}
	We translate \cite[Theorem 11.6]{GLRichardson} into our conventions. Recall the maps $\delta, \Theta$ from \cref{rmk:differencesLeclerc}. Galashin-Lam identify the flag variety with $G/B_-$  rather than $G/B$. Let $\rich{v,w}^-:=(B v B_- \cap B_- w B_-)/B_-$ denote the Richardson variety in $G/B_-$. They define an isomorphism $\vec{\tau}^{\text{pre}}_{v,w}: \rich{v, w}^- \xrightarrow{\sim} {^-\rich{v, w}}$, and set $\vec{\tau}_{v,w}:= \delta^{-1} \circ \vec{\tau}^{\text{pre}}_{v,w}.$
	
	Theorem 11.6 of \cite{GLRichardson} shows that 
	\[f_c \circ  \vec{\tau}^{\text{pre}}_{v,w} = \lmin{c} \circ \Theta^{-1} \]
	as maps on $\rich{v, w}^-$, where $f_c$ are the functions in \cref{rmk:our-minor-vs-leclerc}. So  we see that 
	\begin{align*}
		\lmin{c}&= f_c \circ \delta \circ \Theta \circ \Theta^{-1} \circ \delta^{-1} \circ \vec{\tau}^{\text{pre}}_{v,w} \circ \Theta\\
		&= \rmin{c} \circ \Theta^{-1} \circ \vec{\tau}_{v,w} \circ \Theta
	\end{align*}
where the second equality uses \cref{rmk:our-minor-vs-leclerc}. So the automorphism in the theorem statement is $\Theta^{-1} \circ \vec{\tau}_{v,w} \circ \Theta$.
\end{proof}

\subsection{Background on lattice paths}

\begin{figure}
	\includegraphics[width=0.3\textwidth]{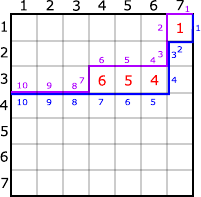}
	\caption{\label{fig:latticePathEx} Let $I:=\{1,3,4\}$ and $J:=\{2,3,7\}$. The path $\pth{I}$ is shown in blue and $\pth{J}$ is shown in purple in $\gridn$. The skew shape $\skew_I^{J}$ has two components, one with content $\{1\}$ and one with content $\{4,5,6\}.$}
\end{figure}

Let $\gridn$ denote the $n \times n$ grid with rows and columns indexed as in a matrix. A \emph{lattice path} in $\gridn$ is a path which begins at the upper right corner of the grid, takes unit length steps down or left, and ends on the left edge of the grid. We label the steps of the lattice path with $1, 2, \dots$ so the labels increase from the beginning of the path to the end (see \cref{fig:latticePathEx} for examples).

\begin{defn}\label{defn:latticePath}
	Let $I \in \binom{[n]}{h}$. We denote by $\pth{I}$ the lattice path of length $n+h$ in $\gridn$ whose vertical steps are labeled with $I$. We denote by $\young{I}$ the Young diagram (in English notation) whose lower boundary is $\pth{I}$.
	
 For $I, J \in \binom{[n]}{h}$, $I \leq J$ in the Gale order if $\pth{I}$ is weakly below $\pth{J}$; or equivalently $\young{J} \subset \young{I}$; or equivalently, writing $I=\{i_1 < \cdots < i_h\}$ and $J = \{ j_1 < \cdots < j_h\}$, we have $i_a \leq j_a$ for $a=1, \dots ,h$.
\end{defn}

Note that for $I \in \binom{[n]}{h}$, the steps $n+1, \dots, n+h$ of $\pth{I}$ are horizontal and the Young diagram $\young{I}$ has $h$ parts, which are all at least $h$.  

For $I \leq J$, we abuse notation and denote by $\skew_I^J$ the skew shape $\young{I}/\young{J}$. We will need to keep track of the ``connected components" of this skew shape.

\begin{defn}\label{defn:content}
	Consider a box $b$ in row $r$ and column $c$ of $\gridn$. The \emph{content} of $b$ is $r-c+n$. So the content of the box in the upper right corner is $1$, and content increases moving down or left. 
	
	For a skew shape $\lam / \mu$ in $\gridn$, the \emph{content} of $\lam / \mu$ is
	\[\{i:\lam / \mu \text{ has a box of content }i\}.\]
	A \emph{component} of $\lam / \mu$ is a maximal-by-inclusion element of 
	\[\{\nu/\rho \subset \lam/\mu: \text{the content of }\nu/\rho \text{ is an interval}\}.\]
\end{defn}

See \cref{fig:latticePathEx} for examples of these definitions.

Note that for $I \leq J \subset [n]$, the content of $\skew_I^{J}$ is contained in $[n-1]$. The components of $\skew_I^{J}$ are the connected components of $(\skew_I^{J}) \setminus (\pth{I} \cap \pth{J})$.

The following easy lemma will be useful to us.

\begin{lem}\label{lem:factors-minors-on-N}
	Consider $I \leq J$. In $\bbc[N]$,
	\begin{enumerate}
		\item if the components of $\skew_I^{J}$ are $\skew_{I_1}^{J_1}, \dots, \skew_{I_r}^{J_r}$, then
		\[\minor{I,J}= \prod_{k=1}^r \minor{I_k, J_k}\]
		and each minor on the right hand side is irreducible. 
		\item if $R \leq S$ is a pair of subsets such that $\skew_R^{S}$ is a translation of $\skew_I^{J}$ parallel to the line $y=-x$, then $\minor{R, S}=\minor{I,J}$.
	\end{enumerate}

\end{lem}
\begin{proof} For a pair of subsets $A \leq B$, say $\pth{A}$ and $\pth{B}$ intersect in steps $C \subset [n]$. Let $A':= A \setminus C$ and let $B':= B\setminus C$.
	
For (1): If $n \in N$, $\minor{I, J}(n)= \minor{I',J'}(n)$.  The submatrix $n'$ of $n$ on rows $I'$ and columns $J'$ is block-upper triangular. The blocks intersecting the diagonal are on rows $I_k'$ and columns $J_k'$. So we have 
\[\minor{I,J}= \minor{I', J'}= \prod_{k=1}^r \minor{I'_k, J'_k}= \prod_{k=1}^r \minor{I_k, J_k}.\]
Irreducibililty follows from \cite[Lemma 3.3]{GLPositroid}, since the skew shapes $\skew_{I_k}^{J_k}$ are connected by construction.

For (2): the content of a box in $\skew_I^{J}$ is the same as the content of the corresponding box of $\skew_R^{S}$. The right edge of a content $k$ box in $\gridn$ is $k$ steps from the upper right corner of $\gridn$. This implies $I'=R'$ and $J'=S'$, which gives the desired equality of minors.
\end{proof}

\begin{cor}\label{cor:irred-factors-skew-shapes}
The irreducible factors of $\rmin{c}$ (as an element of $\bbc[N]$), are the minors corresponding to the components of $\skew_{\vu{c}[h_c]}^{\wu{c}[h_c]}$.
\end{cor}

\begin{rmk}
	In light of \cref{lem:factors-minors-on-N}, we consider skew shapes in $\gridn$ only up to translation parallel to $y=-x$.
\end{rmk}

\section{Ingermanson's and Leclerc's cluster algebras}\label{sec:background2}

We will define Ingermanson's cluster algebra and Leclerc's in parallel, using similar symbols for objects in each construction. ``Ingermanson" is before ``Leclerc" in alphabetical order, so the symbols will follow the same rule (e.g., Ingermanson's cluster variables will be $A_d$, while Leclerc's will be $B_d$).

\subsection{Cluster algebras}
In this section, we set conventions and notation for cluster algebras and related concepts. We refer the reader to e.g. \cite{CAbook} for most definitions.

An \emph{ice quiver} is a directed graph with no loops or 2-cycles, each vertex of which is either \emph{mutable} or \emph{frozen}. A \emph{seed} $\Sigma=(\mathbf{A}, Q)$ in a field $\mcf$ consists of an ice quiver $Q$ together with a tuple $\mathbf{A}$ of elements of $\mcf$, called \emph{cluster variables}, which are indexed by the vertices of the quiver. Cluster variables indexed by mutable vertices are \emph{mutable}; the others are \emph{frozen}. The tuple $\mathbf{A}$ is the \emph{cluster} of $\Sigma$. For each mutable cluster variable $A_i$, we have the corresponding \emph{exchange ratio}
\[\hat{y}_i := \prod_{j \in Q} A_j ^{\# \text{ arrows } A_j \to A_i \text{ in }Q}.\]
By convention, if there are $b$ arrows from $A_i$ to $A_j$ in $Q$, then there are $-b$ arrows from $A_j$ to $A_i$.

There is an involutive operation called \emph{mutation} which can be performed at any mutable vertex of $Q$; this produces a new seed $\Sigma'=(\mathbf{A}', Q')$. The collection of all seeds which can be obtained from $\Sigma$ by a sequence of mutations is the \emph{seed pattern} of $\Sigma$. The \emph{cluster algebra} $\mca(\Sigma) \subset \mcf$ is the $\bbc$-algebra generated by all mutable variables in the seed pattern of $\Sigma$, the frozen variables, and the inverses of the frozen variables.

Let $V$ be an affine variety. We say $\mca(\Sigma)$ is a \emph{cluster structure} on $V$ if $\mca(\Sigma)=\bbc[V]$. If $\mca(\Sigma)$ and $\mca(\Sigma')$ are cluster structures on $V$, then of course $\mca(\Sigma)$ and $\mca(\Sigma')$ are equal as rings. However, their seeds may differ. Two cluster structures $\mca(\Sigma)$ and $\mca(\Sigma')$ are equal if $\Sigma$ and $\Sigma'$ are related by a sequence of mutations; they are \emph{quasi-equivalent} if $\Sigma$ and $\Sigma'$ are related by a sequence of mutations and rescalings by Laurent monomials in frozens which preserve all exchange ratios (see \cite{Fraser} for additional details).
We emphasize that a variety $V$ may have many different cluster structures, which may be quasi-equivalent or not; indeed, Richardsons which are open positroid varieties are known to have many cluster structures \cite{FSB}.

\subsection{Ingermanson's cluster structure}

Fix a Richardson variety $\rich{v, w}$. In \cite{gracie},  Ingermanson defined a seed $\iseed{v, \bw}$ for $\rich{v, w}$ using the \emph{unipeak} expression for $w$. We review her results in this section.

\begin{defn}
	Let $w \in S_n$. A reduced expression $\bw$ is \emph{unipeak} if in $\wiring{\bw}$, no strand travels down and then up.
\end{defn}

\begin{figure}
	\includegraphics[width=\textwidth]{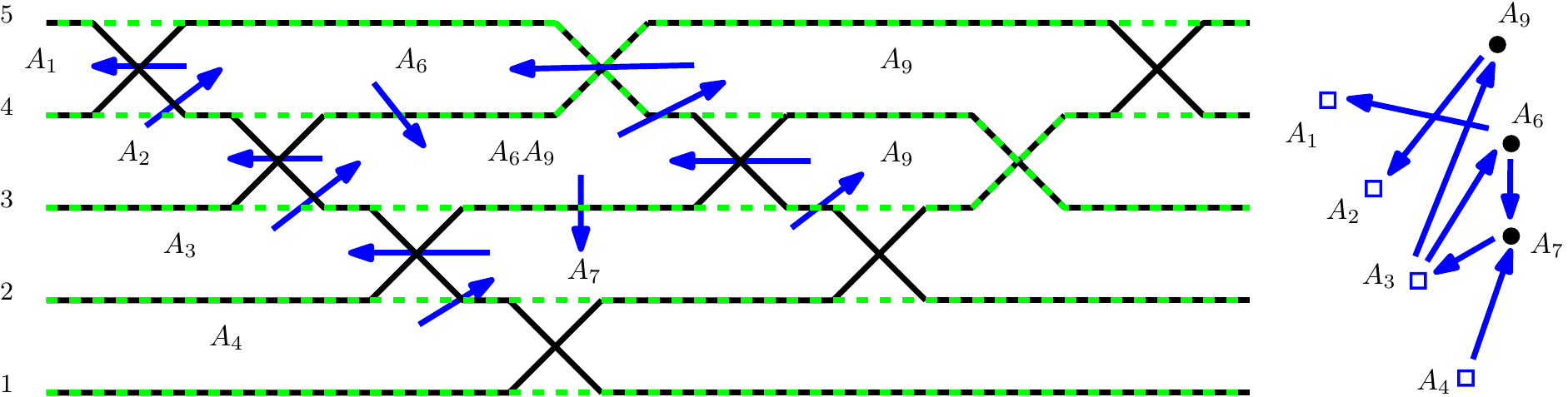}
	\caption{\label{fig:IngEx} Left: a stacked wiring diagram $\wiring{v, \bw}$ for $v=12534$ and $\bw=s_4s_3s_2s_1s_4s_3s_2s_3s_4$. Chambers are labeled by the cluster monomials from \cref{prop:chamber-mono-in-var-gracie}. The wiring diagram quiver is drawn on top. Right: the seed $\iseed{v, \bw}.$}
\end{figure}

The unipeak expressions for $w$ form a nonempty commutation class of reduced expressions for $w$; in particular, every permutation has a unipeak expression \cite{KLR}. \cref{fig:wiringEx} shows a non-unipeak expression; \cref{fig:IngEx} shows a unipeak expression. 

For the remainder of this section, let $\bw$ denote a unipeak expression for $w$.  We will define Ingermanson's seed $\iseed{v, \bw}=(\bA_{v, \bw}, \iquiv{v, \bw})$ in $\bbc[\rich{v, w}]$. The cluster variables $\bA_{v, \bw}$ are indexed by the solid crossings $\jsol{v}$ of $\wiring{\vw}$. We define the cluster variables by giving a monomial map from the left chamber minors to the set of cluster variables. The reader should look ahead to \cref{ex:Ing} for an example.

\begin{defn}\cite[Definition IV.6, Proposition IV.7]{gracie}
	Let $J \subset [n]$ and $u \in S_n$. We define $$\Piv_J(u) := \min_{I \leq J} u(I)$$ where minimum is taken in the Gale order\footnote{The collection $\{I : I \leq J\}$ is a matroid, so it follows from the maximality property for matroids that this minimum is unique (see \cite[Section 1.3]{CoxMatroid} and replace maximum with minimum everywhere.) Ingermanson used the notation Pivots$_J(u)$.}.
\end{defn}

For $1 \leq c<d \leq \ell$, let $L(c,d):=s_{h_{d-1}} s_{h_{d-2}} \cdots s_{h_c}[h_c]$; that is, $L(c,d)$ records the heights of the $w$-strands below $\chi_c$ immediately before crossing $d$. For example, $L(c,c)=[h_c]$. Let $M=(m_{c,d})$ be the matrix whose rows are indexed by $[\ell]$ and whose columns are indexed by $\jsol{v}$, with entries
\begin{equation}\label{eq:is-variable-in-chamber}
	m_{c,d}=\begin{cases}
	0 &\text{ if } c>d\\
	1 &\text{ if } \Piv_{L(c,d)}(\vl{d}s_{h_d})>\Piv_{L(c,d)}(\vl{d}) \\
	0 &\text{ if } \Piv_{L(c,d)}(\vl{d}s_{h_d})=\Piv_{L(c,d)}(\vl{d}).
\end{cases}\end{equation}

Deleting the rows of $M$ indexed by hollow crossings gives a upper unitriangular matrix with 0/1 entries; we denote its inverse by $P=(p_{d,c})$.

\begin{defn} \label{def:Ivar}
	For $d \in \jsol{v}$, we define the cluster variable
	\[A_d := \prod_{c \in \jsol{v}} (\lmin{c})^{p_{d,c}}.\]
\end{defn}

Using \cref{lem:hollowRel} we may express all left chamber minors in terms of cluster variables.

\begin{prop}\cite[Proposition V.1]{gracie}\label{prop:chamber-mono-in-var-gracie}
	For $c \in [\ell]$, we have
	\[\lmin{c}= \prod_{d \in \jsol{v}} (A_d)^{m_{c,d}} .\]
\end{prop}

\begin{defn}
	We say that a cluster variable $A_d$ \emph{appears} in chamber $\chi_c$ of $\wiring{\vw}$ if $m_{c,d}=1$. We denote the (closure of) the union of chambers in which $A_d$ appears by $\jc{d}$.``$\Spr$" stands for ``spread."\footnote{Ingermanson used the notation $JC(j)$ instead; the ``JC" stands for ``jump chambers"  since in these chambers, there is a ``jump" between the two pivot sets used to compute $m_{c,d}$.}
	The cluster variable $A_d$ is frozen if $A_d$ appears in a chamber which is open on the left.
\end{defn}

\begin{example}\label{ex:Ing}
	Let $v=12534$ and $\bw=s_4s_3s_2s_1\underline{s_4}s_3s_2\underline{s_3}s_4$; the hollow crossings are underlined. See \cref{fig:IngEx} for the stacked wiring diagram. We have 
	\[L(6, 9)=s_{h_8}s_{h_7}s_{h_6}[h_6]=s_3 s_2 s_3[3]=\{1,3,4\}  \qquad \text{ and } \qquad \vl{9}=s_4s_3.\]
	We compute $m_{6,9}$, which determines if $A_9$ appears in $\ch_6$.
	\[
		\Piv_{L(6,9)} \vl{9}= \Piv_{134} s_4s_3= \min_{I \leq 134} s_4s_3(I)=123
	\]
	and 
	\[\Piv_{L(6,9)} \vl{9} s_4= \Piv_{134} s_4s_3 s_4=\min_{I \leq 134} s_4s_3s_4(I)= 124.
	\]
	Since $123<124$, $m_{6,9}=1$ and by \cref{prop:chamber-mono-in-var-gracie} $A_9$ appears in $\ch_6$.
\end{example}

To summarize, so far we have a labeling of chambers of $\wiring{\vw}$ by monomials in cluster variables $A_d$ (see \cref{fig:IngEx} for an example). Moving from right to left, $A_d$ first appears in $\ch_d$, and then spreads to other chambers. The appearance of $A_d$ in $\chi_c$ is governed by the two pivot sets in \eqref{eq:is-variable-in-chamber}.

To define $\iquiv{\vw}$, we first draw a quiver on the wiring diagram, following \cite{CA3}.

\begin{figure}
	\includegraphics[width=0.1\textwidth]{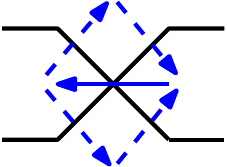}
	\caption{\label{fig:halfArrows} The half arrow configuration used to define the wiring diagram quiver. The horizontal arrow is two half-arrows.}
\end{figure}

\begin{defn}
The \emph{wiring diagram quiver }$\wquiv{\vw}$ has vertices labeled by left chamber minors of $\wiring{\vw}$. The chamber minors in frozen chambers of $\wiring{\vw}$ are frozen; all others are mutable. To determine the arrows, place the configuration of half arrows in \cref{fig:halfArrows}	around each crossing of $\wiring{\bw}$ and sum up the contributions\footnote{That is, to determine the number of arrows from $\lmin{c}$ to $\lmin{d}$: count the number of half-arrows from $\lmin{c}$ to $\lmin{d}$, subtract the number of half-arrows from $\lmin{d}$ to $\lmin{c}$ and divide by 2.}. Delete all arrows between frozen variables.
\end{defn}

\begin{defn} \label{def:IngQuiv}
	Let $B$ denote the square signed adjacency matrix of $\wquiv{v,\bw}$, with rows and columns indexed by $[\ell]$. For $c \in \jsol{v}$ and $c \neq d \in [\ell]$, 
	in $\iquiv{\vw}$ we have
	\begin{align*}\#(\text{arrows }A_c \to A_d)&=\sum_{\chi_a \ni A_c, \chi_b \ni A_d} B_{a, b}\\
		&=\sum_{a, b \in[\ell]} m_{a, c} B_{a,b} m_{b, d}\\
		&= (M^t B M)_{c,d}.
	\end{align*}

Equivalently, let $\hat{y}^W_\ch$ denote the $\hat{y}$-variable for a mutable vertex in $\wquiv{\vw}$ (this is a ratio of left chamber minors). Then
\begin{equation}\label{eq:yHats}
	\hat{y}_c= \prod_{\ch \ni A_c} \hat{y}^W_\ch.
\end{equation}

\end{defn}

In words, for each arrow in $\wquiv{\vw}$ between chambers containing $A_c$ and $A_d$, put an arrow between $A_c$ and $A_d$ in $\iquiv{v, w}$. Then delete 2-cycles.

Ingermanson showed that the upper cluster algebra $\mcu(\iseed{\vw})$ is equal to $\bbc[\rich{v,w}]$. Further, \cite[Corollary 5.8, Remark 7.18]{GLSBS1} shows that the quiver $\iquiv{\vw}$ is locally acyclic, which implies by work of Muller \cite{Muller} that $\mca(\iseed{\vw})= \mcu(\iseed{\vw})$. So we have the following theorem.

\begin{thm}\cite{gracie,GLSBS1} \label{thm:Gracie-clusterstruc}
	Fix $v \leq w$ and let $\bw$ be a unipeak expression for $w$. Then $\mca(\iseed{\vw})= \bbc[\rich{v, w}]$.
\end{thm}

\subsection{Leclerc's cluster structure}\label{}
We recall Leclerc's construction of conjectural cluster structure on $\bbc[\rich{v,w}]$ \cite{Leclerc}. One of the main results of this paper is that his construction does in fact yield a cluster structure. 

Leclerc defines a cluster category inside the module category of a preprojective algebra. Via the cluster character map, this gives rise to a cluster subalgebra of $\bbc[\rich{v,w}]$.  For a more detailed exposition of the representation theoretic construction we refer to \cite[Section 5]{SSBW}.

The preprojective algebra $\Lambda_{n-1}$ of type $A$ is the path algebra of the quiver 
$$P= \xymatrix@!{1 \ar@/^/[r]^{\alpha_1}& \ar@/^/[l]^{\alpha_1^*} 2 \ar@/^/[r]^{\alpha_2}& \ar@/^/[l]^{\alpha_2^*} 3 \ar@/^/[r]^{\alpha_3}  &\cdots   \ar@/^/[l]^{\alpha_3^*} \ar@/^/[r]^{\alpha_4} & \ar@/^/[l]^{\alpha_4^*} n-1}$$  
with relations 
$$\sum_{i} \alpha_i \alpha_i^*-\alpha_i^*\alpha_i=0.$$
A module $U$ over $\Lambda_{n-1}$ is obtained by placing a $\mathbb{C}$-vector space $U_i$ at each vertex $i$ of $P$ and linear maps between these vector spaces $\phi_{\alpha_i}:U_i\to U_{i+1}$ and $\phi_{\alpha_i^*}:U_{i+1}\to U_{i}$ for each arrow of $P$, such that the maps satisfy the relations given above.  Let $\dim\, U:= (\dim\, U_{i})_{i \in [n-1]}$ denote the dimension vector of $U$. The support of $U$ is the set of all vertices $i$ in the quiver such that $U_i\not=0$.  Let $\left|U\right|$ denote the number of pairwise non-isomorphic indecomposable direct summands of $U$, and let $\add U$ denote the full subcategory of the module category whose objects are direct sums of summands of $U$.

We will be interested in a special type of $\Lambda_{n-1}$-modules, which correspond to skew shapes in $\gridn$. Let $\lam/\mu \subset \gridn$ be a skew shape with content in $[n-1]$. The $\Lambda_{n-1}$-module $U_{\lam/\mu}$ is as follows. Recall the boxes of $\gridn$ are indexed as in a matrix. Each box $b=b_{i,j}$ of $\lam/\mu$ with content $c$ yields a basis vector $e_{i,j}$ of $(U_{\lam/\mu})_{c}$ and the maps are defined as follows. 
\[ \phi_{\alpha_{c}}(e_{i,j}) = \begin{cases}
e_{i+1,j} & \text{if } b_{i+1,j}\in  \lam/\mu \\
0 & \text{else}
\end{cases}   \hspace{1cm} 
\phi_{\alpha^*_{c-1}}(e_{i,j}) = \begin{cases}
e_{i,j+1} & \text{if } b_{i,j+1}\in  \lam/\mu \\
0 & \text{else}
\end{cases} \]

For example, if $\lam/\mu$ is the $(n-k)\times k$ rectangle whose lower right corner has content $k$, then $U_{\lam/\mu}$ is the indecomposable injective $\Lambda_{n-1}$-module at vertex $k$, and if $\lam/\mu$ is a single content $k$ box then $U_{\lam/\mu}$ is the simple $\Lambda_{n-1}$-module at vertex $k$, which we will denote by $S(k)$.

From this description it follows that the top of $U_{\lam/\mu}$ is a direct sum of simple modules $S(k)$, one for each box $b_{i,j} \in  \lam/\mu$ with content $k$ such that $b_{i-1,j}, b_{i,j-1}\not\in \lam/\mu$.  Similarly, the socle of $U_{\lam/\mu}$ is a direct sum of simple modules $S(k)$, one for each box $b_{i,j} \in  \lam/\mu$ with content $k$ such that $b_{i+1,j}, b_{i,j+1}\not\in \lam/\mu$. That is, $b_{i,j}$ is a content $k$ corner of the northwest or southeast boundary of $\lam/\mu$ respectively. 

A module $U_{\lam'/\mu'}$ is a submodule of $U_{\lam/\mu}$ if $\lam'/\mu' \subseteq \lam/\mu$ and for every $b_{i,j}\in \lam'/\mu'$ whenever $b_{i+1,j}, b_{i,j+1} \in \lam/\mu$ then $b_{i+1,j}, b_{i,j+1} \in \lam'/\mu'$ respectively.  On the other hand,  a module $U_{\lam'/\mu'}$ is a quotient of $U_{\lam/\mu}$ if $\lam'/\mu' \subseteq \lam/\mu$ and for every $b_{i,j}\in \lam'/\mu'$ whenever $b_{i-1,j}, b_{i,j-1} \in \lam/\mu$ then $b_{i-1,j}, b_{i,j-1} \in \lam'/\mu'$ respectively.  Any map of two modules $f: U_{\lam/\mu}\to U_{\nu/\rho}$ is determined by its image $\text{im}\,f$, where $\text{im}\,f=U_{\lam'/\mu'}$ is a quotient of $U_{\lam/\mu}$ and a submodule of $U_{\nu/\rho}$.  

\begin{rmk}\label{rmk:components-give-indec-summands}
	Consider a skew shape $\lam/\mu$ in $\gridn$ with content contained in $[n-1]$. The components of $\lam/\mu$ (c.f. \cref{defn:content}) give the indecomposable summands of $U_{\lam/\mu}$.
\end{rmk}

Recall that a pair of subsets $I \leq J \in \binom{[n]}{h}$ determines a pair of Young diagrams $\young{I} \supset \young{J}$ and a skew shape $\skew_I^{J}$ (c.f. \cref{defn:latticePath}). So the pair $I \leq J$ also determines a $\Lambda_{n-1}$-module, which we denote $U_{I,J}$.

Given $v\leq w$, Leclerc defines a certain subcategory $\mathcal{C}_{v,w}$ of the module category of $\Lambda_{n-1}$ which he showed admits a cluster structure in the sense of \cite{BIRS}. The category $\mathcal{C}_{v,w}$ is equipped with a cluster character map
\begin{align*}
	\varphi : \obj \mathcal{C}_{v,w}&\to \bbc[\rich{v,w}]\\
	U & \mapsto \varphi_U
\end{align*}
satisfying $U \oplus U' \mapsto \varphi_U \varphi_{U'}$. Each (reachable) \emph{cluster tilting module} $U$ of $\mcc_{v,w}$ corresponds to a seed in $\bbc(\rich{v,w})$; the cluster variables are the images of the indecomposable summands of $U$ under $\varphi$.

To obtain a cluster tilting module, we define a module for each chamber of the wiring diagram.

\begin{figure}
	\includegraphics[width=\textwidth]{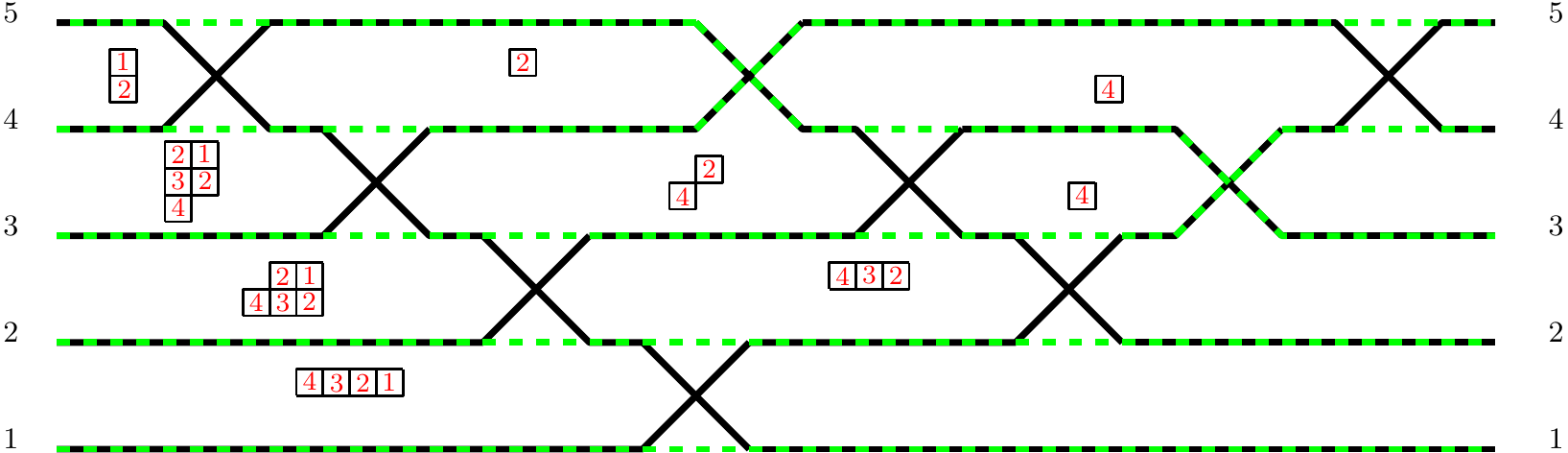}
	\caption{\label{fig:LecEx} The stacked wiring diagram $\wiring{v, \bw}$ for $\bw=s_4s_3s_2s_1s_4s_3s_2s_3s_4$ and $v=12534$, with chambers labeled by chamber modules (or by skew shapes).}
\end{figure}

\begin{defn}
	Fix $v \leq \bw$. For $c \in [\ell]$, the \emph{chamber module} is
		\begin{equation}\label{eq:chamberModule}
		U_c :=U_{\vu{c}[h_c],\wu{c}[h_c]}.
	\end{equation}
\end{defn}
See \cref{fig:LecEx} for an example of a stacked wiring diagram with chambers labeled by chamber modules. By results of Leclerc, $\varphi$ maps $U_c \mapsto \rmin{c}$.

The main result of \cite{Leclerc} can be formulated as follows. 

\begin{theorem}\cite[Theorem 4.5]{Leclerc}\label{Lec-seed} Fix $v \leq \bw$.
		$$U_{v,{\bf w}}:= 
		\bigoplus_{c\in \jsol{v}} U_c$$ is a cluster tilting object in $\mathcal{C}_{v,w}$. The corresponding seed $\lseed{\vw}=(\bB_{\vw}, \lquiv{\vw})$ in $\mathbb{C}[\rich{v, w}]$ can be described as follows.

\begin{itemize}

\item[(a)] The cluster variables are the $\ell(w)-\ell(v)$ irreducible factors\footnote{We mean the irreducible factors of $\rmin{c}$ as a function on $N$.} of 
\[\prod_{c \in \jsol{v}} \varphi_{U_c} = \prod_{c \in \jsol{v}} \rmin{c}.\]
The set of cluster variables is the $\varphi$-image of the set of indecomposable summands of the $U_c$.  

\item[(b)] A cluster variable is frozen if it is a factor of the right chamber minor of a frozen chamber in $\wiring{v, \bw}$. The frozen variables are $\varphi$-images of the  
	indecomposable summands of $\bigoplus_{i\in I} U_{v^{-1}([i]),w^{-1}([i])}$ (which 
		are the projective-injective objects).

	\item[(c)] The quiver $ \lquiv{v, \bw}$ is the endomorphism quiver of the cluster tilting module.  In particular, the vertices of the quiver are nonisomorphic indecomposable summands of $U_{v,{\bf w}}$ and the arrows are irreducible morphisms in $\add U_{v,{\bf w}}$ between these summands.
\end{itemize}
Finally, the cluster algebra
$\mca(\lseed{v, \bw})$
is a subalgebra of 
$\mathbb{C}[\rich{v, w}]$.
\end{theorem}

\begin{rmk}
	Analogously to Ingermanson's construction, in Leclerc's construction the right chamber minors of $\wiring{\vw}$ are cluster monomials. This is clear for chambers to the left of a solid crossing by construction; for chambers to the left of a hollow crossing, this follows from \cref{lem:hollowRel}. We say that a cluster variable $B \in \bB_{\vw}$ \emph{appears} in a chamber $\chi_c$ if $B$ is an irreducible factor of the right chamber minor $\rmin{c}$.
\end{rmk}

In certain special cases, Leclerc showed that $\mca(\lseed{v, \bw}) = \mathbb{C}[\rich{v, w}]$. He conjectured that this equality holds in general.

\begin{conj}\label{conj:leclerc} \cite{Leclerc}
	The cluster algebra $\mca(\lseed{\vw})$ is equal to $\bbc[\rich{v,w}]$.
\end{conj}

Our main result is that this conjecture is true. Moreover Leclerc's seeds for different reduced expressions $\bw, \bw'$ are related by mutation (c.f. \cref{prop:LecSeedMutationEquiv}), so Leclerc's construction gives a single cluster structure on $\bbc[\rich{v, w}]$.

	\section{Correspondence between cluster variables} \label{sec:clusterVar}

In this section, we show the relationship between the clusters $\bA_{v, \bw}$ and $\bB_{\vw}$ for $\bw$ unipeak. In Ingermanson's seed, the cluster variables are labeled by $\jsol{v}$; the variable $A_d$ ``first appears" in the chamber $\ch_d$. We start by pointing out Leclerc's cluster variables have an analogous labeling.

\begin{lem}\label{lem:Lec-var-label-by-solid}
	Let $v \leq \bw$ and let $d \in \jsol{v}$. There is a unique cluster variable $B_d \in \bB_{\vw}$ which appears in $\ch_d$ and does not appear in $\ch_c$ for $c>d$.
\end{lem}

\begin{proof}
	We proceed by induction on $\ell(w)$. The base case is $\ell(w)=0$, where the desired statement is vacuously true.
	
	The right chamber minor $\rmin{c}$ and its irreducible factors (as functions on $N$) depend only on the prefixes $\vu{c}, \wu{c}$ of $\bw^{-1}$ and $\bv^{-1}$. So we may assume $d=1$ is the first crossing of $\bw$. Let $\bw':= s_{i_2} \dots s_{i_\ell}$ and $v':= s_{i_2}^{v} \dots s_{i_\ell}^{v}$. Cutting off the first crossing of $\wiring{\vw}$ gives the wiring diagram $\wiring{v', \bw'}$. Applying \cref{Lec-seed} for $\rich{v',w'}$, we see that the right chamber minors appearing in $\wiring{v',w'}$ have $\ell(w')-\ell(v')$ irreducible factors. Since $\ell(w)=\ell(w')+1$ and $\ell(v)=\ell(v')$, \cref{Lec-seed} tells us that $\rmin{1}$ has exactly one irreducible factor which is not an irreducible factor of any chamber minor $\rmin{c}$ for $c>1$. 
\end{proof}

If $\bw$ is unipeak, we can describe $B_d$ in more detail.

\begin{lem}\label{lem:LecVar}
	Let $v \leq \bw$ with $\bw$ unipeak and let $d \in \jsol{v}$. Then $B_d$ is the minor corresponding to the northeast-most component of $\skew_{\vu{d}[h_d]}^{\wu{d}[h_d]}$. If $k$ is the final vertical step on the northwest boundary of this component, then $B_d$ is equal to the minor $\minor{[k] \cap \vu{d}[h_d], [k]\cap \wu{d}[h_d]}$ on $N$.
\end{lem}

\begin{proof}
	Let $\chi_{d'}$ denote the chamber directly to the right of crossing $d$.  Let $\rmin{d}=\Delta_{I,J}$ and $\rmin{d'}=\Delta_{I',J'}$ be the corresponding right chamber minors. If $\chi_{d'}$ is open on the right, then $I'=J'=[i_\ell]$ and $\Delta_{I',J'}=1$.  
	
	Since $d \in \jsol{v}$, we have $I'=I$.  Moreover, $J$ and $J'$ differ by a single element, say $J=J'\setminus\{r'\}\cup\{r\}$.  Note that $r'<r$ because the wiring diagram is reduced, and $r$ and $r'$ are the right endpoints of the rising and falling strands of $d$, respectively.
	
	 In addition, since $\bw$ is unipeak, the falling strand at $d$ must always go down to the right of $d$.  This implies that the strands with right endpoints $1, \dots, r'-1$ are below $\chi_d$ and $\chi_{d'}$.  Otherwise, these strands would need to cross above the falling strand of $d$ at some crossing to the right of $d$, which is impossible. Thus, $[r']\subset J'$.  Since $I\leq J'$ then all of $I, J, J'$ contain $[r'-1]$, and in particular all of the lattice paths $\pth{I}, \pth{J}, \pth{J'}$ agree on the first $r'-1$ steps, which are all vertical.  Then in the lattice path $\pth{J'}$, step $r'$ is vertical and step $r$ is horizontal; in the lattice path $\pth{J}$ the step $r'$ is horizontal and step $r$ is vertical.  Hence the skew shapes $\skew_I^{J'}$ and $\skew_I^{J}$ differ by a strip on the northwest boundary between steps $r'$ and $r$. This strip is contained in the northeast-most component of $\skew_I^{J}$, so all other components of the two skew shapes are the same.
\end{proof}

The goal of this section is to prove the following statements.

\begin{thm}\label{thm:factorizationsAgree}
	Let $v \leq \bw$ with $\bw$ unipeak. Choose $d \in \jsol{v}$ and chamber $\chi$ in the wiring diagram $W_{v, \bw}$. Then the cluster variable $A_d$ appears in $\chi$ if and only if the cluster variable $B_d$ appears in $\chi$.
\end{thm}

Using \cref{thm:factorizationsAgree}, we can prove the precise relationship between $A_d$ and $B_d$. Recall the maps $\alpha, \beta$ from \cref{lem:richIso}.

\begin{thm}\label{thm:varCorrespondence} Fix $v \leq \bw$ with $\bw$ unipeak and choose $d \in \jsol{v}$. Moreover, let $gB \in \rich{v,w}$ and set $x:= \alpha^{-1}(gB)$ and $y:=\beta^{-1}(gB)$.
	\begin{enumerate}
		\item $A_d = B_d \circ \tau_{v,w}$.
		\item If $B_d(gB)= \Delta_{I, J}(y)$, then $A_d(gB)=\Delta_{v(I), w(J)}(x)$. That is, if we express $B_d$ and $A_d$ as minors, the row sets are related by an application of $v$ and the column sets are related by an application of $w$.
	\end{enumerate}
	
\end{thm}

\begin{rmk}
	\cref{thm:factorizationsAgree} and \cref{thm:varCorrespondence} do not follow immediately from \cref{prop:twist}. For example, \cref{prop:twist} does not rule out the possibility that for some $d$, $A_d = A_c (B_d \circ \tau)$ where $A_c$ is a frozen variable. 
\end{rmk}

For the remainder of this section, we fix $v \leq \bw$ with $\bw=s_{i_1} \dots s_{i_\ell}$ unipeak. 

We will prove \cref{thm:factorizationsAgree} in the next three sections. We first show that \cref{thm:factorizationsAgree} holds in a special case, and then show that we can reduce to the special case. We prove \cref{thm:varCorrespondence} in \cref{sec:proofOfVarCorresp}

\subsection{Base case} In this section we show the following lemma.

\begin{lemma}\label{lem:base_case}
Suppose the final crossing $\ell$ is solid. Then $A_\ell$ appears in $\lmin{1}$ if and only if $B_\ell$ appears in $\rmin{1}$.
\end{lemma}

\begin{proof}
Since $\ell\in \jsol{v}$, $s_{i_\ell}$ is not in the PDS for $v$ in $\bw$. Set $j:=i_{\ell}$ and $k=i_1$. Then $B_\ell = \minor{j, j+1}$ is a one-by-one minor; the corresponding skew shape is a single box with content $j$.  It follows that $B_\ell$ appears in $\rmin{1}$ if and only if all of the following conditions hold:
\begin{itemize}
\item[(1)] $j\not\in w^{-1}[k]$ and $j+1\in w^{-1}[k]$;
\item[(2)] $j\in v^{-1}[k]$ and $j+1\not\in v^{-1}[k]$;
\item[(3)] $\left| w^{-1}[k]\cap [j-1] \right| = \left| v^{-1}[k]\cap [j-1] \right|$.
\end{itemize}

Now let $w'=ws_{j}$.  Then on the other hand, $A_\ell$ appears in $\Delta_1^\lambda$ if and only if $\Piv_L(vs_{j})>\Piv_L(v)$ where $L=(w')^{-1}[k]$.  First, $v^{-1}\leq (w')^{-1}$ implies that $v^{-1}[k] \leq L$, so we have $\Piv_L(v)=[k]$. Then the same argument implies that $\Piv_L(vs_{j})>[k]$ if and only if $s_{j}v^{-1}[k]\not\leq (w')^{-1}[k]$.  

To prove the lemma it suffices to show that conditions (1)-(3) above are equivalent to the condition that $s_{j}v^{-1}[k]\not\leq (w')^{-1}[k]$.  Since $v\leq w'$, it is easy to see that $s_{j}v^{-1}[k]\not\leq (w')^{-1}[k]$ if and only if all of the following conditions hold: 
\begin{itemize}
\item[(1')] $j\in (w')^{-1}[k]$ and $j+1\not\in (w')^{-1}[k]$;
\item[(2)] $j\in v^{-1}[k]$ and $j+1\not\in v^{-1}[k]$;
\item[(3')] $\left| (w')^{-1}[k]\cap [j-1] \right| = \left| s_{j}v^{-1}[k]\cap [j-1] \right|$.
\end{itemize}

From $w=w's_{j}$, (1) and (1') are equivalent, as are (3) and (3').
\end{proof}

\subsection{Leclerc's factorizations are stable under left and right multiplication}
In this section, we show that the apearance of $B_d$ in a chamber $\ch$ does not change under removing prefixes and suffixes from $\bw$.

We need some definitions involving wiring digrams which differ by a single crossing at the left or right. Notice that if $\bw < \bw \cdot s_i$, each chamber $\chi$ of $\wiring{\bw}$ corresponds naturally to a chamber $\chi$ of $\wiring{\bw \cdot s_i}$, and similarly if $\bw < s_i \cdot \bw$.

\begin{defn}\label{def:right}
	Let $\bw'=\bw_{(\ell-1)}$ and $v'=v_{(\ell-1)}$. Note that $\jsol{v}\setminus \{\ell\}=\jsol{v'}$. We denote cluster variables in $\iseed{v', \bw'}$ and $\lseed{v', \bw'}$ by $A'_j$ and $B_j'$, respectively.	For $d \in \jsol{v'}$, the appearance of $A_d'$ (resp. $B_d'$) is \emph{stable under right multiplication} if for all chambers $\chi$ of $\wiring{v', \bw'}$, $A'_d$ (resp. $B_d'$) appears in $\chi$ in $\wiring{v', \bw'}$ if and only if $A_d$ (resp. $B_d$) appears in $\ch$ in $\wiring{v, \bw}$. 
\end{defn}

We make exactly analogous definitions for left multiplication.     

\begin{defn}\label{def:left}
	
	Let $\bw = s_{i_1} \dots s_{i_\ell}$ and $v\leq w$. Let $\bw'=s_{i_2} \dots s_{i_\ell}$ and $v'=s^v_{i_2} \dots s^v_{i_\ell}$. Note that the crossings of $\bw'$ are indexed by $2, \dots, \ell$, and that $\jsol{v}\setminus \{1\}=\jsol{v'}$. We denote cluster variables in $\iseed{v', \bw'}$ and $\lseed{v', \bw'}$ by $A'_j$ and $B_j'$, respectively. For $d \in \jsol{v'}$, the appearance of $A_d'$ (resp. $B_d'$) is \emph{stable under left multiplication} if for all chambers $\chi$ of $\wiring{v', \bw'}$, $A'_d$ (resp. $B_d'$) appears in $\chi$ in $\wiring{v', \bw'}$ if and only if $A_d$ (resp. $B_d$) appears in $\ch$ in $\wiring{v, \bw}$. 
\end{defn}

\begin{lem}\label{lem:Lec_left}
	In the setup of \cref{def:left}, let $j \in \jsol{v'}$. The appearance of $B'_j$ is stable under left multiplication.
\end{lem}
\begin{proof}
	Notice that all right chamber minors ${(\rmin{c})}'$ of $\wiring{v', \bw'}$ are equal to the corresponding right chamber minor $\rmin{c}$ of $\wiring{v, \bw}$. Also, $B'_j=B_j$. The claim follows.
\end{proof}

\begin{prop}\label{lem:Lec_right}
In the setup of \cref{def:right}, let $j \in \jsol{v'}$. The appearance of $B'_j$ is stable under right multiplication. 
\end{prop}

\begin{proof}
	Fix a chamber $\ch_c$. Note that the appearance of a cluster variable $B'_j$ in a right chamber minor $\rmin{c}$ does not depend on the prefix $s_{i_1} \dots s_{i_{c-1}}$. So we may assume without loss of generality that $c=1$ is the first crossing of $\bw'$. Also, we set $k:=i_\ell$.

First, suppose that $B_j'$ appears in $(\rmin{1})'=\Delta_{v'^{-1}[i_1],{w'}^{-1}[i_1]}$.  Then $B_j' =\Delta'_{R,S}$ for some $R \subset [r, s]$ and $S \subset [ r+1, s+1]$ with $s\geq r$ where $r\in R$ and $s+1\in S$.  The corresponding lattice paths for $B_j'$ and $(\Delta_1^\rho)'$ are shown in Figure~\ref{LecFact} on the left.  Next, it will be more convenient to work with modules and lattice paths instead of minors.   Let $M_j', M_j$ be indecomposable modules that correspond to $B_j', B_j$ respectively.  Also, let $X,Y$ be summand of $U_1'$ that are adjacent to $M_j'$ as in the figure.  To prove the proposition, we will compare the indecomposable module $M_j'$ with $M_j$ and also the chamber module $U'_1$ with $U_1$ (c.f. \eqref{eq:chamberModule}).

\begin{figure}
 \scalebox{.6}{\Large\input{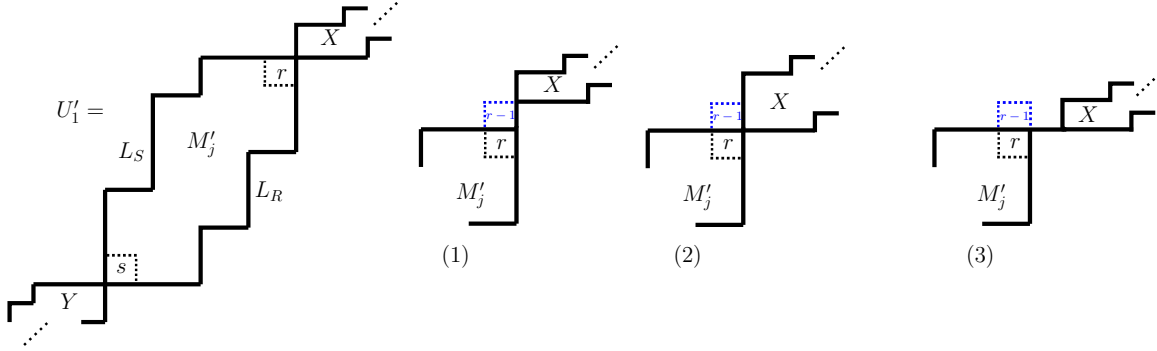}}
  \vspace*{-4cm}
\caption{The relation between $U_1$ and $U_1'$ in the proof of Proposition~\ref{lem:Lec_right}.}
\label{LecFact}
\end{figure}

If $k<r-1$ or $k>s+1$, then $M_j'=M_j$ and $U_1$ is obtained from $U'_{1}$ by possibly adding and/or removing a content $k$ box. In particular, we see that $M_j$ is a summand of $U_1$ as desired. Similarly, if $k\in[r,s]$ then $M_j$ is obtained from $M_j'$ by possibly adding and/or removing a content $k$ box, while $U_1$ is obtained from $U_{1}'$ in the exact same way by adding and/or removing the corresponding box to the summand $M_j'$ of $U'_{1}$. Hence, we conclude again that $M_j$ is a summand of $U_1$.  

Next, suppose that $k=r-1$. There are several cases to consider (see Figure~\ref{LecFact}, right).  Note that by assumption we have $r\not\in {w'}^{-1}[i_1]$ and $r\in (v')^{-1}[i_1]$.

(1) Suppose that $k=r-1 \in {w'}^{-1}[i_1]$, and $r-1\in (v')^{-1}[i_1]$.   Then $U_1$ is obtained from $U'_1$ by adding a content $r-1$ box, and the summands $X,Y$ do not change.  Similarly, since $M'_j, M_j$ are the topmost summands in $U'_j, U_j$ respectively, then $M_j$ is obtained from $M'_j$ by adding a content $r-1$ box.  Hence, we conclude that $M_j$ is a summand of $U_1$ as desired.

(2) Suppose that $k=r-1 \in {w'}^{-1}[i_1]$ and $r-1\not \in (v')^{-1}[i_1]$.  Since $r\in (v')^{-1}[i_1]$ and $r-1\not \in (v')^{-1}[i_1]$, then we conclude that $s_{r-1}(v')^{-1}[i_1] < (v')^{-1}[i_1]$, so then $v's_{r-1}<v'$.  Thus, $v=v'$, but in the expression for $w=w's_{r-1}$ the reduced expression for $v'=v$ contains $s_{r-1}$. This contradicts the assumptions in Definition~\ref{def:right} that $v'=v_{(\ell-1)}$.  

(3) Suppose $k=r-1\not\in (w')^{-1}[i_1]$.  Then $r-1\not\in(v')^{-1}[i_1]$.  In this case, it may happen that $M_j$ is obtained from $M_j'$ by adding a content $r-1$ box.  However, $U'_1=U_1$, so we observe that $M_j'$ and not $M_j$ is a summand of $U_1$.  We will show that this leads to a contradiction.  Let $U_p$ contain $M_j'$ as a summand with $p$ being maximal, that is $M_j'=M_p$.   By \cref{lem:LecVar} the module $M_j'$ is the topmost summand in $U_p$. Note that $p\not=j$ since $M_j\not=M_p$.  Then $p\in \jsol{v'}$ and $M_j'$ is then also the topmost summand of $U'_p$.  Again \cref{lem:LecVar} implies that $p$ is the maximal index such that $M_j'$ is a summand of $M'_p$.  Therefore, $p=j$, which is a contradiction. 

This completes the proof in the case $k=r-1$ and the other remaining situation is when $k=s+1$, which can be shown in a similar way.  This shows that if $M_j'$ is a summand of  $U'_1$ then $M_j$ is a summand of  $U_1$.   It remains to show the converse. 

Now, suppose that $M_j$ is a summand of $U_1$, and we want to show that $M_j'$ is a summand of $U'_1$.   Again let $i_{\ell}=k$ and we consider several cases. 

If $M_j = M_j'$ and $M_j$ is a summand of $U_1$ then one cannot add a content $k$ box to the top of $M_j'$ in $U_1$ or remove a content $k$ box from the bottom of $M_j'$ in $U_1$.  Then $M_j'$ is also a summand of $U'_1$ because $U'_1$ is obtained from $U_1$ by possibly removing a content $k$ box from the top and/or adding it to the bottom.  Hence, it remains to show the case when $M_j\not=M_j'$.

If $M_j = S(k)$ is a simple module represented by a single box and $M_j'=0$, then $s_{k}=s_{i_{\ell}}$ at the end of $\bw$ is not part of the reduced expression for $v$ inside $\bw$, as otherwise $M_j$ would be zero.  But then  $M_{\ell} = S(k)$ so $\ell=j$ which is a contradiction, since we assume $j<\ell$. 

Now suppose that $M_j'$ is obtained from $M_j$ just by removing a content $k$ box from the top.  Then if in addition $U_1\not= U'_1$, then $U'_1$ is similarly obtained from  $U_1$ by removing a content $k$ box from the top.  Hence, $M_j'$ is then a summand of $U'_1$ as desired.  Otherwise, if  $U_1= U'_1$, then $M_j$ and not $M'_j$ is a summand of  $U'_1$.  But then $M_j=M'_p$ for some $p\not=j$, and so $M_j=M'_p=M_p$, where the last equality holds since $M'_p$ already has content $k$ box in the top.  Hence, $j=p$, which is a contradiction.
The same argument applies if $M_j'$ is obtained from $M_j$ just by adding a box with content $k$ to the bottom.  


Finally, suppose that $M_j'$ is obtained from $M_j$ by both removing a content $k$ box from the top and by adding a content $k$ box to the bottom.  Then $k$ would be a vertical step while $k+1$ would be a horizontal step for both the top and the bottom contour of $M_j'$. This implies that $M_j'$ and $M_j$ both contain a box with content $k-1$ and a box with content $k+1$.  Hence, $k$ is not a minimal or maximal content of a box in $M_j'$ and $M_j$.  This implies that if $M_j$ is a summand of $U_1$ then $M'_j$ is then a summand of $U'_1$ as desired. 

This completes all the cases and proves the proposition. 
\end{proof}

%
%

\subsection{Ingermanson's factorizations are stable under left and right multiplication} Here, we show that the apearance of $A_d$ in a chamber $\ch$ does not change under removing prefixes and suffixes from $\bw$.

Let $c<d$ index two crossings in a diagram for $(v, \bw)$. Recall from \eqref{eq:is-variable-in-chamber} and \cref{prop:chamber-mono-in-var-gracie} that $A_d$ appears in the chamber $\chi_c$ if and only if $$\Piv_{L(c,d)}(\vl{d-1} s_{i_d})> \Piv_{L(c,d)}(\vl{d-1}),$$ where $L(c,d):= s_{i_{d-1}}s_{i_{d-2}} \cdots s_{i_c}[i_c]$.

\begin{lem}\label{lem:IRightStable}
	In the setup of \cref{def:right}, let $d \in \jsol{v'}$. Then the appearance of $A'_d$ is stable under right multiplication.
\end{lem}

\begin{proof} Choose a chamber $\ch_c$.
	The condition above for the appearance of $A'_d$ in $\ch_c$ depends only on $\vl{d-1}$, $s_{i_d}$ and the subword $s_{i_c} \cdots s_{i_{d-1}}$ of $w$. None of this data changes under the right multiplication taking $\bw'$ to $\bw$ and $v'$ to $v$.
\end{proof}

To show that the appearance of $A_d$ is stable under left multiplication, we first need a lemma about how pivot sets interact with left multiplication.

\begin{lem}\label{lem:leftMultPivots}
	Fix $I \subset [n]$ and $u \in S_n$. Choose $s_i \in S_n$ such that $u < s_iu$. Set 
	\[u(P):= \min_{J \leq I} u(J) \quad \text{and} \quad s_iu(Q):=\min_{J \leq I} s_iu(J) . \]
	Then either 
	\begin{enumerate}
		\item $s_iu(Q)=u(P)$ and $Q=P$.
		\item $s_iu(Q)=u(P)$ and $Q=(P \setminus \{u^{-1}(i)\} )\cup \{u^{-1}(i+1)\}$.
		\item $s_iu(Q)=(u(P) \setminus\{i\}) \cup \{i+1\}$ and $Q=P$.
	\end{enumerate}
In particular, $u(P) \leq s_iu(Q)$, that is, $\Piv_I(u) \leq \Piv_I(s_iu)$.
\end{lem}

\begin{proof}
	Suppose first that $Q=P$ so that $s_iu(Q)=s_iu(P)$. If neither or both of $i,i+1$ is in $u(P)$, then we are in situation $(1)$. So we may assume only one of $i, i+1$ is in $u(P)$. Suppose for contradiction that $i+1 \in u(P)$ and $i \notin u(P)$. Because $u<s_i u$, we have $u= \dots i \dots i+1 \dots$. But this means that $P'=P \setminus \{u^{-1}(i+1)\} \cup \{u^{-1}(i)\}$ is smaller than $P$ and satisfies $u(P')<u(P)$, a contradiction.
	
	Now, suppose that $Q \neq P$. Then $u(P)<u(Q)$ and $s_iu(Q)<s_iu(P)$, by the definition of $Q$ and $P$. Notice that $s_iu(Q)$ and $u(Q)$ differ by at most one element, and similarly for $s_iu(P)$ and $u(P)$. It is not hard to see (for example, by considering the lattice paths $\pth{u(P)}, \pth{u(Q)}$, etc.) that this implies $u(P)=s_iu(Q)$.
\end{proof}

\begin{lem} \label{lem:IleftStable}
	In the setup of \cref{def:left}, let $d \in \jsol{v'}$. Then the appearance of $A'_d$ is stable under left multiplication.
\end{lem}

\begin{proof}
Say $i_1=k$ and fix $c<d$. The variable $A_d$ appears in the chamber $\chi_c$ if and only if $\Piv_I(\vl{d-1} s_{i_d})> \Piv_I(\vl{d-1})$, where $I:= s_{i_{d-1}}s_{i_{d-2}} \cdots s_{i_c}[i_c]$. Clearly the subset $I$ is the same for both $\bw$ and $\bw'$. If $v=v'$, then $\vl{d-1}=\vl{d-1}'$, so the condition above also determines the appearance of $A_d$ in $\chi_c$ in the diagram for $(v', \bw')$.
	
	 So it suffices to consider the situation where $v=s_k v'$, and in fact, $v>v'$. Let $u:=\vl{d-1}'=s^v_{i_2} \cdots s^v_{i_{d-1}}$ and say $i_d=h$. We would like to show the following: 
	 \begin{align*}
	 	\min_{G \leq I} u(G) \neq \min_{G \leq I} u s_h(G)  \iff \min_{G \leq I} s_k u(G) \neq \min_{G \leq I} s_k u s_h(G).
	 \end{align*}
 Say that the minima above are acheived at $M, N, Q, P$ respectively, so we would like to show that 
\begin{align}\label{eq:pivot}
	u(M) \neq us_h(N) \iff s_ku(Q) \neq s_k u s_h(P).
\end{align}

Note that, since $d \in \jsol{v'}$ and $d\in \jsol{v}$ we have $u<u s_h$ and $s_k u < s_k u s_h$.

By \cref{lem:leftMultPivots}, of the following possiblities
\begin{enumerate}[label={(\Alph*)}]
	\item $u(M)=s_k u(Q)$
	\item $M=Q$ and $s_k u(Q)=u(M)\setminus \{k\} \cup \{k+1\}$
	\item $u s_h(N) =s_k u s_h (P)$
	\item $N=P$ and $s_k u s_h(P)=u s_h(N)\setminus \{k\} \cup \{k+1\}$
\end{enumerate}
exactly one of (A), (B) hold and exactly one of (C), (D) hold.

If (A) and (C) hold, we clearly have \eqref{eq:pivot}.

If (A) and (D) hold, then tracing inequalities gives $s_k u(Q) < s_k u s_h(P)$. Suppose for the sake of contradiction that $u(M)=us_h(N)$. Various assumed equalities imply $s_k u s_h(P)=s_k u (Q)\setminus \{k\} \cup \{k+1\}$.  Then \cite[Proposition IV.45]{gracie} implies that in the diagram for $(v, \bw)$, the $v$-strand $\gamma_{k+1}$ with left endpoint $k+1$ lies above $\chi_c$ and the $v$-strand $\gamma_k$ lies below $\chi_c$. But these strands cross at the first crossing, so this is impossible.

If (B) and (C) hold, tracing inequalities gives $u(M)<u s_h(N)$. If $s_k u(Q) = s_k u s_h(P)$, then our assumptions imply $u s_h (N)=u(M) \setminus \{k\} \cup \{k+1\}$. This is impossible: because $s_k u s_h >u s_h$, in $u s_h$ we see $\dots k \dots k+1 \dots$, so by the minimality of $u s_h (N)$, if $u s_h (N)$ contains $k+1$ it also contains $k$.

If (B) and (D) hold, then $s_k u(Q)=u(M) \setminus \{k\} \cup \{k+1\}$ and $s_k u s_h(P)=u s_h(N) \setminus \{k\} \cup \{k+1\}$. So \eqref{eq:pivot} clearly holds.
\end{proof}

We can now prove \cref{thm:factorizationsAgree}.

\begin{proof}[Proof of \cref{thm:factorizationsAgree}]

Let $v \leq w$ and let $\bw=s_{i_1} \cdots s_{i_\ell}$ be a unipeak expression for $w$. Consider $d \in \jsol{v}$ and a chamber $\chi_c$ in $\wiring{v, \bw}$ which is to the left of crossing $c$. We would like to show that $A_d$ appears in $\chi_c$ if and only if $B_d$ appears in $\chi_c$. We may assume that $c<d$, since in both constructions, the cluster variable indexed by $d$ only appears in chambers to the left of crossing $d$.

Consider the pair $(v', w')$ where $\bw'=s_{i_c} \cdots s_{i_d}$ and $v'=s^v_{i_c} \cdots s^v_{i_d}$ (where crossings are labelled $c, c+1, \dots, d$ to avoid confusion). Repeated application of \cref{lem:IRightStable} and \cref{lem:IleftStable} implies that $A_d$ appears in $\chi_c$ in the diagram $\wiring{v, \bw}$ if and only if $A'_d$ appears in $\chi_c$ in the diagram $\wiring{v', \bw'}$. Similarly, repeated application of \cref{lem:Lec_left,lem:Lec_right} implies an identical statement for $B_d$. In the diagram $\wiring{v', \bw'}$, we are in exactly the situation of \cref{lem:base_case}, so $A_d$ appears in $\chi_c$ if and only if $B_d$ appears in $\chi_c$.

\end{proof}

\subsection{Relating $A_i$ and $B_i$}\label{sec:proofOfVarCorresp}
In this section, we prove \cref{thm:varCorrespondence}. 

%

We first recall some of \cite[Chapter V]{gracie}. Choose a crossing $d \in \jsol{v}$. The chamber minor $\lmin{d}$ is $\Delta_{R,S}$ for some $R, S \subset [n]$. By \cite[Proposition V.8]{gracie}, the cluster variable $A_d$ is equal to $\Delta_{I, J}$ where $I=R\cap v[h]$ and $J=S \cap w[h]$ for some $h \in [n]$. The number $h$ is the height of the right endpoint of the path $\pi_d$ (defined below).

\begin{figure}
	\includegraphics[width=\textwidth]{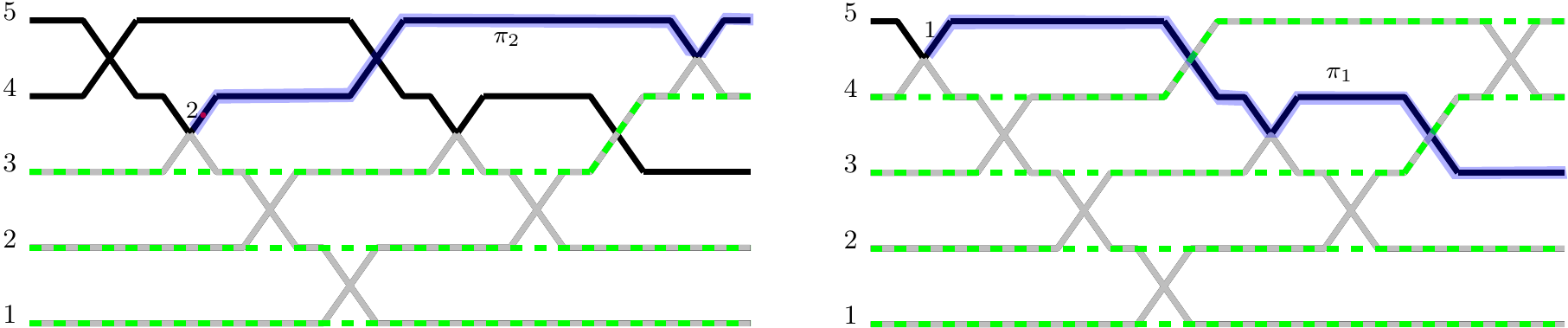}
	\caption{\label{fig:piEx} Two examples of paths $\pi_d$, which are highlighted in blue. The edges of $W'$ are shown in black; the deleted $w$-strand segments are grey. For clarity, only the $v$ strands below $\ch_d$ are drawn.}
\end{figure}

\begin{defn}\label{defn:wiringPath}
	Fix $(v, \bw)$ and choose a crossing $d \in \jsol{v}$. A \emph{segment} of a strand is a connected component of $\wiring{v, \bw} \setminus \{\text{intersections of $w$-strands}\}$. We modify $\wiring{v, \bw}$ by deleting $w$-strand segments which touch a $v$-strand that passes below $\chi_d$; call the resulting graph $W'$.
	
	Then $\pi_d$ is a path on $W'$ which is (the closure of) a union of strand segments. It begins at crossing $d$ and travels to the right, ending at the right edge of $W'$. The first strand segment $\pi_d$ travels along is the rising strand of $d$. If $\pi_d$ arrives at a crossing of $W'$, immediately after the crossing $\pi_d$ always follows the rising strand. See \cref{fig:piEx} for examples.
\end{defn}

\begin{rmk}\label{rmk:pathFacts}\ 
	\begin{enumerate}
		\item Drawing $\pi_d$ on the full wiring diagram $\wiring{v, \bw}$, note that $\pi_d$ can only go from height $h$ to height $h-1$ at a hollow crossing. \label{itm:pathDownHollowOnly}
		\item It follows from the proof of \cite[Proposition 5.8]{gracie} that if $\pi_d$ approaches a crossing $c$ along the rising strand, it will also exit $c$ on the rising strand. In particular, the rising strand of $c$ will be present in $W'$. \label{itm:pathRisingAlwaysOk}
		\item It follows from \cite[Proposition 5.8]{gracie} that if a strand $\alpha$ goes from above $\pi_d$ to below, then $\alpha$ does not pass below $\chi_d$. \label{itm:pathStrandNotBelow}
	\item  It follows from \cite[Proposition 3.16]{gracie} that if $\pi_d$ follows a segment of strand $\alpha$, then $\alpha$ passes below $\chi_d$. \label{itm:pathStrandBelow}
	\end{enumerate}

\end{rmk}

\begin{proof}[Proof of \cref{thm:varCorrespondence}]
	For 1): 	\cref{thm:factorizationsAgree} implies that the same upper unitriangular matrix $P$ gives the monomial map from left chamber minors to $\bA_{v, \bw}$ and from right chamber minors to $\bB_{\vw}$. That is, we have
	\[A_d = \prod_{c \in \jsol{v}} (\lmin{c})^{p_{d,c}}  \qquad \text{ and } \qquad B_d = \prod_{c \in \jsol{v}} (\rmin{c})^{p_{d,c}}.\]
	To see that $B_d \circ \tau_{v,w}=A_d$, precompose both sides with $\tau_{v,w}$ and use \cref{prop:twist}.
	
	For 2): Fix $\bw=s_{i_1} \cdots s_{i_\ell}$ and $v \leq w$. Choose $d \in \jsol{v}$.
	
Let $\Delta_{R,S}$ be the chamber minor $\lmin{d}$, so $\rmin{d}=\Delta_{v^{-1}(R),w^{-1}(S)}$. Let $\mu$ be the skew shape bounded by lattice paths $\pth{v^{-1}(R)}$ and $\pth{w^{-1}(S)}$ and let $\mu_{[q]}$ denote the region between lattice paths $\pth{v^{-1}(R) \cap [q]}$ and $\pth{w^{-1}(S) \cap [q]}$. 

From \cref{lem:LecVar}, Leclerc's cluster variable $B_d$ is the minor $\Delta_{I, J}$ where $I= v^{-1}(R) \cap [q]$ and $J=w^{-1}(S) \cap [q]$ and $q \in [n]$ is the smallest number such that the region $\mu_{[q]}$ is a skew shape with a single connected component. That is, the content of $\mu_{[q]}$ is a nonempty interval ending at $q-1$.

On the other hand, by the discussion above, Ingermanson's cluster variable $A_d$ is the minor $\Delta_{R\cap v[h], S \cap w[h]}$ where $h$ is the height of the right endpoint of $\pi_d$. We would like to show that $h=q$. 

First, note that $v^{-1}(R) \cap [h]$ and $w^{-1}(S)\cap [h]$ are the same cardinality, so $\mu_{[h]}$ is a skew shape. Further, \cref{rmk:pathFacts} (4) and the definition of $\pi_d$ implies that $h \in w^{-1}(S)$ but $h \notin v^{-1}(R)$. So $\mu_{[h]}$ has a box with content $h-1$. Thus, we just need to show that the content of $\mu_{[h]}$ is an interval.

We will proceed by induction.

If $d$ is the rightmost crossing, this is true by inspection, as the cluster variables in the two constructions are equal to chamber minors. Otherwise, consider $\bw':=s_{i_1} \cdots s_{i_{\ell-1}}$ and $v':=\vl{\ell -1}$. Let $\pi_d'$ be the path in $\wiring{v', \bw'}$ and $h'$, $\mu'$ and $q'$ defined analogously as above. Say ${i_\ell}=a$. Then $\mu$ is obtained from $\mu'$ by first adding a content $a$ box along the top boundary if possible and then deleting a content $a$ box along the bottom boundary if possible (where ``if possible" means ``if the resulting collection of boxes is a skew shape"). 

If $h' \notin \{a, a+1\}$, then $\pi_d'$ and $\pi_d$ end at the same height so $h=h'$. If $a>h$, then the inductive hypothesis easily implies $\mu_{[h]}$ is equal to $\mu'_{[h]}$. If $a+1<h$, the fact that the number of content $a$ boxes in $\mu$ and $\mu'$ differ by at most one and straightforward casework confirms that the content of $\mu_{[h]}$ is an interval. In either case, $q=h$ as desired.

If $h'=a$, then $\pi_d$ ends at $a+1$ by \cref{rmk:pathFacts} (2). The top connected component of $\mu'$ ends at content line $a-1$ by the inductive hypothesis. By \cref{rmk:pathFacts} (3) and (4), it is possible to add a box of content $a$ along the top boundary of $\mu'$, but it is easy to check one cannot delete a box of content $a$ along the bottom. Thus, $\mu$ has top connected component ending at content line $a=h-1$.

If $h'=a+1$, there are two cases depending on whether $\pi_d$ goes down or straight at the final crossing. If $\pi_d$ goes down, then $h=a$. Since $\mu$ and $\mu'$ are identical at and above content line $a-1$, by induction the content of $\mu_{[h]}$ is an interval. If $\pi_d$ goes straight at the final crossing, then $h=h'$. By induction the content of $\mu'_{[h]}$ is an interval ending at $a$. Since the content of $\mu_{[h]}$ and $\mu_{[h]}'$ can differ only by $a$ and $\mu_{[h]}$ contains a box of content $a$, we are done.

%

\end{proof}

\section{Leclerc's quiver in terms of wiring diagrams}\label{sec:LQuivFromWiring}

In this section we analyze the morphisms between the indecomposable summands of the cluster tilting object $U_{v, \bw}$.  Our goal is to characterize irreducible morphisms in terms of the wiring diagram, so that we can ultimately compare $\lquiv{v, \bw}$ and $\iquiv{v, \bw}$.

Throughout this section we fix $v\leq \bw$ where $\bw=s_{i_1}\dots s_{i_\ell}$ is unipeak.  For $j\in\jsol{v}$ let $M_j$ denote the indecomposable summand of $U_j$ corresponding to $B_j$. Recall that by \cref{lem:LecVar} the module $M_j$ is the topmost summand of $U_j$. 

\subsection{Morphisms coming from neighboring chambers}\label{sec:map-descripiton}

Let $\chi$ and $\chi'$ be chambers adjacent to a solid crossing $i \in \jsol{v}$, and let $U_\ch$ and $U_{\ch'}$ be the corresponding chamber modules. Let $\alpha$ be the falling $w$-strand at the crossing $i$ with right endpoint $a$ while $\alpha'$ be the rising $w$-strand at $i$ with right endpoint $a'$. Let $\Delta_{I,J}, \Delta_{I',J'}$ denote the right chamber minors for the chambers $\chi,\chi'$ respectively. We will define explicit morphisms between the modules $U_{\chi}$ and $U_{\chi'}$.  There are three cases depending on the relative positions of the chambers $\chi$ and $\chi'$. Recall the notation $\chi_{i^{\uparrow}}, \chi_{i^{\rightarrow}}, \chi_{i^{\downarrow}}, \chi_{i^{\leftarrow}}$ for the chambers above, to the right, below, and to the left of $i$ respectively.  

\begin{figure}
 \scalebox{.6}{\Large\input{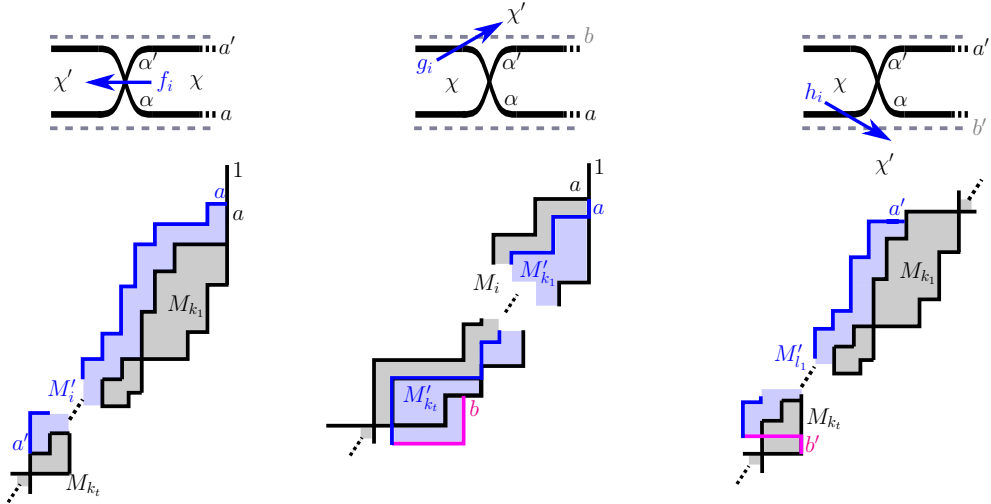}}
\caption{Construction of the maps $f_i, g_i, h_i$. On the left the module $U_{\chi'}$ is obtained from $U_\chi$ by adding a strip between $a$ and $a'$ to the top, in the middle the module $U_{\chi'}$ is obtained from $U_\chi$ by removing a strip starting at $a$ from the top and adding a strip starting at $b$ to the bottom, and on the right the module $U_{\chi'}$ is obtained from $U_\chi$ by adding a strip to the top starting at $a'$ and removing a strip from the bottom starting at $b'$.}
\label{Fig:fgh}
\end{figure}

First, suppose that $\chi = \chi_{i^{\rightarrow}}$ and $\chi'=\chi_{i^{\leftarrow}}$ (Figure~\ref{Fig:fgh}, left).   Since $i\in \jsol{v}$, we see that $I=I'$ and $J'=J\setminus\{a\}\cup\{a'\}$.  
By the same reasoning as in the proof of \cref{lem:LecVar}, we see that $a'>a$ and $1, \dots, a \in I\cap J$.  
In particular, we see that $M'_i$, the topmost summand in $U_{\chi'}$, is obtained from the first few topmost summands $M_{k_1}, \dots, M_{k_t}$ of $U_{\chi}$ by adding a strip to the top between indices $a$ and $a'$. So there is an inclusion $f_i:\oplus_{j=1}^t M_{k_j} \to M'_i$.


Second, suppose that $\chi=\chi_{i^{\leftarrow}}$ and $\chi'=\chi_{i^{\uparrow}}$ (Figure~\ref{Fig:fgh}, middle).   Then $J'=J\cup \{a\}$ and $I'=I\cup \{b\}$, where $b$ is the right endpoint of the $v$-strand above $i$.  Since the strand $\alpha$ moves downward to the right of $i$, we see that $b>a$.  Moreover, similarly to the case above, we see that $1, \dots, a-1 \in J$.   Hence, $a$ is the smallest integer that is not in $J$. From this, the module $U_{\chi'}$ is obtained from $U_{\chi}$ by removing a strip from the top of $M_i$ and adding a smaller strip to the bottom of $M_i$.  Hence, $M_i$ results from a direct sum of the first few topmost summands $\oplus_{j=1}^t M'_{k_j}$ of $U_{\chi'}$ and there is a surjection $g_i: M_i \to \oplus_{j=1}^t M'_{k_j}$. 

Third, suppose that $\chi=\chi_{i^{\leftarrow}}$ and $\chi'=\chi_{i^{\downarrow}}$ (Figure~\ref{Fig:fgh} right).   Then $I'=I\setminus\{b'\}$ and $J'=J\setminus\{a'\}$ where $a', b'$ are the right endpoints of the rising $w$-strand at $i$ and the $v$-strand traveling below $i$ respectively.  Then $U_{\chi'}$ is obtained from $U_{\chi}$ by both adding a strip to the top and also removing a strip from the bottom.  This gives a map $h_i: \oplus_{j=1}^t M_{k_j}\to \oplus_{j=1}^{t'} M'_{l_j}$  from a certain collection of summands of $U_{\chi}$ to a certain collection of summands of $U_{\chi'}$.  Observe, that if $a'\geq b'$ then $t=1$, and if $a'<b'$ then $t'=1$.  Note that in general, the map $h_i$ is neither injective nor surjective.  Moreover, in special cases it may even be zero. 

Note that $f_i, g_i, h_i$ are defined above on certain summands of $U_{\chi}, U_{\chi'}$. We extend them to the modules $U_{\chi}, U_{\chi'}$ by mapping the remaining summands either via the identity map when possible or the zero map.  

We summarize the results of this subsection below.

\begin{proposition}\label{prop:maps}
Let $\chi,\chi'$ be two chambers adjacent to the crossing $i\in \jsol{v}$.  Let $a,a'$ denote the right endpoints of the falling and rising $w$-strands at the crossing $i$ respectively.  Similarly, let  $b, b'$ be the right endpoints of the $v$-strands that pass just above and just below the crossing $i$ respectively.

\begin{itemize}
\item[(a)] If $\chi =  \chi_{i^{\rightarrow}}$ and $\chi'=\chi_{i^{\leftarrow}}$ then $U_{\chi'}$ is obtained from $U_{\chi}$ by adding a strip to the top between steps $a$ and $a'$.  In particular, there exists an injective morphism $f_i: U_\chi\to U_{\chi'}$. 
\item[(b)] If $\chi =  \chi_{i^{\leftarrow}}$ and $\chi'=\chi_{i^{\uparrow}}$ then $U_{\chi'}$ is obtained from $U_{\chi}$ by removing a strip from the top starting at $a$ and adding a strip to the bottom starting at $b$.  In particular, there exists a surjective morphism $g_i: U_\chi\to U_{\chi'}$. 
\item[(c)] If $\chi =  \chi_{i^{\leftarrow}}$ and $\chi'=\chi_{i^{\downarrow}}$ then $U_{\chi'}$ is obtained from $U_{\chi}$ by adding a strip to the top starting at $a'$ and removing a strop from the bottom starting at $b'$.  In particular, there exists a morphism $h_i: U_\chi\to U_{\chi'}$ whose  image is $U_{\chi}$ without the strip starting at $b'$. 
\end{itemize}
\end{proposition}

By construction it is easy to see that the maps defined above have the largest possible images.  If $q: M\to M'$ is one of $f_i, g_i, h_i$, then there may be other morphisms $r: M \to M'$. However by construction $r$ would factor through $q$ as $r(M)$ would lie inside $q(M)$.  Since we are interested in describing irreducible maps between indecomposable summands in Leclerc's seed, it suffices to consider only the maps $f_i, g_i, h_i$ between the appropriate summands.  

\subsection{Arrows in Leclerc's quiver}

The goal of this section is to show that arrows in Leclerc's quiver come from morphisms between the modules appearing in the neighboring chambers.   We begin with some preliminary lemmas. 

\begin{lemma}\label{top:k}
Let $\bw = s_{i_1} \dots s_{i_\ell}$ and $v \leq \bw$. Let $w' = w_{(\ell-1)}$ and $v' = v_{(\ell-1)}$.  The simple module $S(i_\ell)$ is not a summand of the top of $U'_j$ for any $j \in [\ell-1]$.
\end{lemma}

\begin{proof}
Let $i_\ell = k$. If the simple $S(k)$ is a summand of the top of $U'_j$ then $k+1\in (w')^{(j)}[i_j]=s_{i_{l-1}}\dots s_{i_{j+1}}s_{i_j}[i_j]$ and $k\not\in (w')^{(j)}[i_j]$.  Then $k\in w^{(j)}[i_j]$ and $k+1\not\in w^{(j)}[i_j]$, since $w$ is obtained by multiplying $w'$ on the right by $s_k$.  However, this contradicts the assumption that $(w')^{(j)}<w^{(j)}$.   
\end{proof}

\begin{lemma}\label{lem:maps}
Let $\bw = s_{i_1} \dots s_{i_\ell}$ and $v \leq w$. Let $w' = w_{(\ell-1)}$ and $v' = v_{(\ell-1)}$.  For every $i,j<\ell$ there is a bijection between irreducible morphisms $f: M_i\to M_j$ in  $\textup{add}\,U_{v,\bw}$ and irreducible morphisms $f': M_i'\to M_j'$ in $\textup{add}\,U_{v',\bw'}$ such that $\textup{im}\,f$ and $\textup{im}\,f'$ differ possibly at vertex $k:=i_\ell$.  In particular, there is an exact sequence 
\[S(k)\to \textup{im}\, f' \to \textup{im}\, f \to S(k).\]
\end{lemma}

\begin{proof}
First, suppose that $\ell\in \jhol{v}$.  Suppose that there exists a nonzero morphism $f':M'_i\to M'_j$ such that $\textup{im}\,f'$ is indecomposable.  Now consider $\textup{im}\,f'$ inside $M_i'$.  If the top contour of the image, and hence also the top contour of $M'_i$, has vertical step $k$ and horizontal step $k+1$ then let $\textup{im}\,f$ be obtained from $\textup{im}\,f'$ by adding $S(k)$ to the top.  Moreover, if $\textup{im}\,f'$ contains $S(k)$ in the socle then in addition let $\textup{im}\,f$ be obtained from $\textup{im}\,f'$ by removing $S(k)$.   Otherwise, let $\textup{im}\,f=\textup{im}\,f'$.    Recall that $M_i, M_j$ are obtained from $M_i', M_j'$ respectively by adding $S(k)$ to the top when possible and then also removing $S(k)$ from the socle when possible.  It is easy to see that $f'$ induces the corresponding nonzero map $f:M_i\to M_j$ with image $\textup{im}\,f$ as described above.  Note that $\textup{im}\,f'\not\cong S(k)$ by Lemma~\ref{top:k}, which means that the induced map $f$ is indeed nonzero. 

Conversely, we can also see that every nonzero $f: M_i\to M_j$ induces a corresponding nonzero $f':M_i'\to M_j'$ by the same reasoning.  Here note that $\textup{im}\,f\not\cong S(k)$ as $M_j$ cannot have $S(k)$ in the socle since $\ell\in \jhol{v}$, so indeed $f'$ is nonzero.    This shows that there is a bijection between nonzero morphisms $f: M_i\to M_j$ in the seed $\textup{add}\,U_{v,\bw}$ and nonzero morphisms $f': M_i'\to M_j'$ in $\textup{add}\,U_{v',\bw'}$ such that the images of $f$ and $f'$ differ at most by $S(k)$ in the top and socle.  Since every such morphism is a composition of irreducible ones, we obtain the desired bijection between irreducible morphisms as in the statement of the lemma.   This completes the proof in the case $\ell\in \jhol{v}$.

Now, suppose that $\ell\in \jsol{v}$.  Then a similar argument as above implies that there is a bijection between nonzero morphisms $f':M_i'\to M_j'$ and nonzero morphisms $f:M_i\to M_j$ such that $\textup{im}\,f\not\cong S(k)$.  Note that if $\textup{im}\,f\cong S(k)$ then by construction the corresponding morphism $M_i'\to M_j'$ would be zero.  However, in $\textup{add}\,U_{v,\bw}$ the module $M_{\ell}=S(k)$ since $\ell\in \jsol{v}$. Hence, for $i,j<\ell$ no morphism $M_i\to M_j$ with image $S(k)$ can be irreducible in $\textup{add}\,U_{v,\bw}$, as it would factor through $M_{\ell}$.  
Therefore, we obtain the desired bijection between irreducible morphisms $f':M_i'\to M_j'$ in $\textup{add}\,U_{v',\bw'}$ and irreducible morphisms $f:M_i\to M_j$ in $\textup{add}\,U_{v,\bw}$ as claimed. 
\end{proof} 

\begin{proposition}\label{prop:neigh_ch}
Fix an arrow $\alpha$ in $\lquiv{v, \bw}$. Then there is a solid crossing $j \in \jsol{v}$ such that $\alpha$ comes from an arrow between chambers $\chi, \chi'$ of $\wiring{v, \bw}$ that are adjacent to $j$.
\end{proposition}

\begin{proof}
Let $\bw = s_{i_1} \dots s_{i_\ell}$ and $v \leq \bw$.  Set $k:=i_\ell$.  We proceed by induction on the length of $w$.  If $\ell=1$ then there are no arrows in Leclerc's quiver so the statement holds.  Now suppose that $\ell>1$, and let $w' = w_{(\ell-1)}$ and $v' = v_{(\ell-1)}$.  Observe that the wiring diagrams $\wiring{v, \bw}$ and $\wiring{v', \bw'}$ differ by a single chamber and arrows between the common chambers remain the same for the two wiring diagrams.  Similarly, by \cref{lem:maps} we know that the irreducible morphisms between the corresponding indecomposable summands of $U_{v,\bw}$ and $U_{v',\bw'}$ also remain the same.  Therefore, it suffices to consider the irreducible morphisms of $U_{v,\bw}$ that start and end at $S(k)=M_{\ell}$ in the case that $\ell \in \jsol{v}$, and show that they can be realized as arrows in the wiring diagram between the neighboring chambers. 

First, suppose that there is a nonzero map $S(k)\to M_j$ for some $j<\ell$ with $j\in  \jsol{v}$.  Then $M_j$ has $S(k)$ in the socle and $M_j\not\cong S(k)$. Hence, the bottom contour of $M_j$ has a vertical step $k$ and a horizontal step $k+1$.   This means that $M_j$ belongs to a chamber $\chi_j$ that lies below the $v$-strand with right endpoint labeled $k+1$ and above the $v$-strand with right endpoint labeled $k$. 
Since $j\in   \jsol{v}$ then the $v$-strands do not cross at $j$.  This means that the chamber $\chi_{j'}$ directly to the right of $\chi_j$ corresponds to a module $U_{j'}$ that also has $S(k)$ in the socle, as the bottom contours of $U_j$ and  $U_{j'}$ remain the same.  Moreover, by \cref{prop:maps}(a) there exists an injective morphism $f_j:M'\to M_j$ where $M'$ is an indecomposable direct summand of $U_{j'}$.   Since $U_{j'}$ has $S(k)$ in the socle then so does $M'$, and hence we obtain the following injective morphism $S(k)\to M'\to M_j$.   This implies that our starting map $S(k)\to M_j$ is irreducible if and only if $S(k)=M'$.  In particular, this occurs if and only if $S(k)$ is a summand of $U_{j'}$, and hence such a morphism arrises from an arrow between the neighboring chambers $\chi_j$ and $\chi_{j'}$. 

Now, suppose that there is a nonzero map $M_j\to S(k)$ for some $j<\ell$ with $j\in  \jsol{v}$. 
Then $S(k)$ belongs to the top of $M_j$ and $M_j\not\cong S(k)$.  Hence, the top contour of $M_j$ has a horizontal step $k$  and a vertical step $k+1$.   In particular, this means that $M_j$ belongs to a chamber $\chi_j$ that lies below the $w$-strand with right endpoint $k$ and below the $w$-strand with right endpoint $k+1$. Let $\chi_{j'}$ be the chamber directly above the crossing $j$.  Next, we consider two cases.  

If the corresponding module $U_{j'}$ also has $S(k)$ at the top, then by  \cref{prop:maps}(b) there exists a surjective morphism $g_j: M_j\to M'$ where $M'$ is an indecomposable direct summand of $U_{j'}$.   Hence, composing the two maps we obtain a surjective morphism $M_j\to M' \to S(k)$.   We obtain that the map $M_j\to S(k)$ is irreducible if and only if $S(k)=M'$, and hence such a morphism arrises from an arrow between the neighboring chambers $\chi_j$ and $\chi_{j'}$. This completes the proof in the first case. 

Now, suppose that $U_{j'}$ does not have $S(k)$ at the top.  By unipeakness the $w$-strand with the right endpoint $k+1$ is still below $\chi'$, but then the $w$-strand with the right endpoint $k$ must also be below $\chi'$. 
Then in the wiring diagram the falling $w$-strand at crossing $j$ must have right endpoint $k$.  Moreover, by unipeakness this strand must continue to move weekly downward as it moves to the right after the crossing $j$.   This implies that if $b'$ is the right endpoint of the $v$-strand that passes just below the crossing $j$ then $b'> k$. Note that $b'\not=k$ since the last crossing $\ell$ of $w$ is solid and between height $k$ and $k+1$.  Now, let $a'$ denote the right endpoint of the rising $w$-strand at $j$.  We see that $a'>k+1$.  By  the description of modules and morphisms in \cref{prop:maps}(c), it follows that $S(k)$ is still at the top of $U_{j''}$ where $\chi_{j''}$ is the chamber in $\wiring{v, \bw}$ directly below the crossing $j$ and the composition $M_j \to U_{j''} \to S(k)$ is nonzero.  This shows that $M_j \to S(k)$ is irreducible if and only if $S(k)$ is a summand of $U_{j''}$, and hence such a morphism arrises from an arrow between the neighboring chambers $\chi_j$ and $\chi_{j''}$.
\end{proof}

\subsection{Irreducible morphisms from $\jc{j}$}\label{sec:map-JC}

We showed that each arrow in Leclerc's quiver, which corresponds to an irreducible morphism between the indecomposable modules in the seed, comes from one of the maps $f, g, h$ between two chambers adjacent to a solid crossing. Next we fix $j \in \jsol{v}$ and determine which maps between chambers give an arrow between $M_j$ and $M_i$ for $i<j$. It is easy to see that if $M_j$ is a direct summand of both $U_{\chi}$ and $U_{\chi'}$, that is $\chi, \chi'\in \jc{j}$, then for any map $f: U_{\chi}\to U_{\chi'}$, the image $f(M_j)$ is contained inside $M_j$.  Since Leclerc's quiver does not have loops, no such map yields an arrow in Leclerc's quiver. Therefore, when analyzing irreducible morphisms involving $M_j$ it suffices to consider arrows between two neighboring chambers where exactly one of them belongs to $\jc{j}$.

We start with some lemmas about the boundary of $\jc{j}$.

\begin{defn}\label{def:ends}
	Let $c$ be a crossing in $\wiring{\vw}$. A crossing $d$ is a \emph{left end} of $\jc{c}$ if the chamber $\ch^\leftarrow$ to the left of $d$ is not in $\jc{c}$ and the chamber $\ch^\rightarrow$ to the right of $d$ is in $\jc{c}$. A \emph{right end} of $\jc{c}$ is defined similarly. A crossing $d$ on the boundary of $\jc{c}$ is a \emph{cusp} if an odd number of the surrounding chambers are in $\jc{c}$. Note that a cusp may not be an end, and vice versa (see \cref{fig:irredMaps} for examples).
\end{defn}

\begin{lem}\label{lem:propogateRightAndUp}
	Let $c \in \jsol{v}$ and let $d<c$ be a solid crossing to its left. If the chamber $\ch_d$ is in $\jc{c}$, then so is the chamber $\ch^\rightarrow$ to the right of $d$ and the chamber $\ch^\uparrow$ above $d$.
\end{lem}

\begin{proof}
Recall from the proof of \cref{lem:LecVar} that the skew-shapes $\skew_I^{J}$ and $\skew_I'^{J'}$ respectively labeling $\ch_d$ and $\ch^\rightarrow$ differ only by a strip contained in the topmost component of $\skew_I^{J}$. The component $\mu$ of $\skew_I^{J}$ corresponding to cluster variable $B_c$ is not the topmost component of $\skew_I^{J}$ since $d<c$ and $d$ is solid. So $\mu$ is also a component of $\skew_I'^{J'}$ and thus $B_c$ appears in $\ch^\rightarrow$. That is, $\ch^\rightarrow$ is in $\jc{c}$. A similar argument shows the same statement for $\ch^\uparrow.$
\end{proof}

\begin{cor}\label{cor:rightEndHollow}
	Let $c$ be a solid crossing in $\wiring{\vw}$. If $d \neq c$ is a right end of $\jc{c}$, then $d$ is a hollow crossing.
\end{cor}

\begin{cor}\label{cor:forbidden-cusp}
	No cusp of $\jc{c}$ is of the form 
	\begin{center}
		\includegraphics[width=0.3\textwidth]{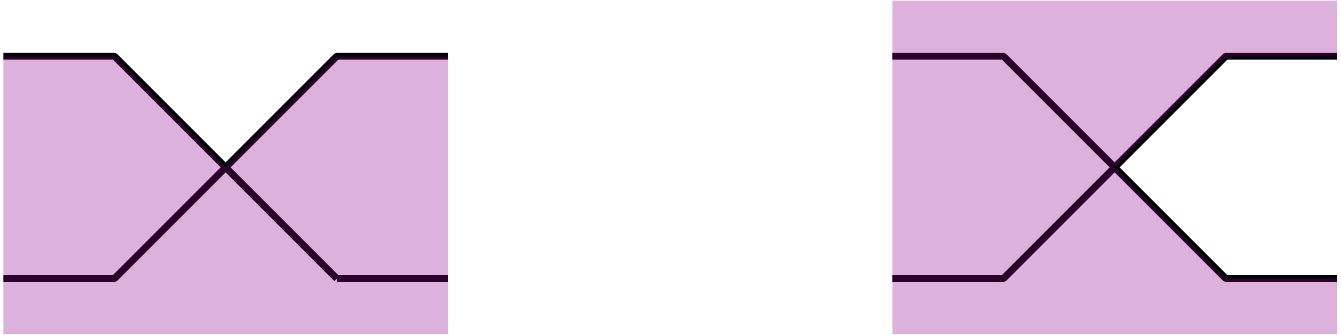}
	\end{center}
	where the shaded chambers are in $\jc{c}$ and white chambers are not.
\end{cor}

\begin{rmk}
	Using \cref{cor:rightEndHollow} and the fact that all arrows in $\lquiv{v, \bw}$ come from two chambers around a solid crossing, we may consider only arrows between a chamber $\chi \in \jc{j}$ and a chamber $\chi'\notin \jc{j}$ which touches a left end or a cusp of $\jc{j}$.
\end{rmk}

In the next few lemmas we consider several cases. Recall from \cref{prop:maps} the definition of the maps $f, g, h$ that we associate to arrows in the wiring diagram. 

\begin{figure}
	\includegraphics{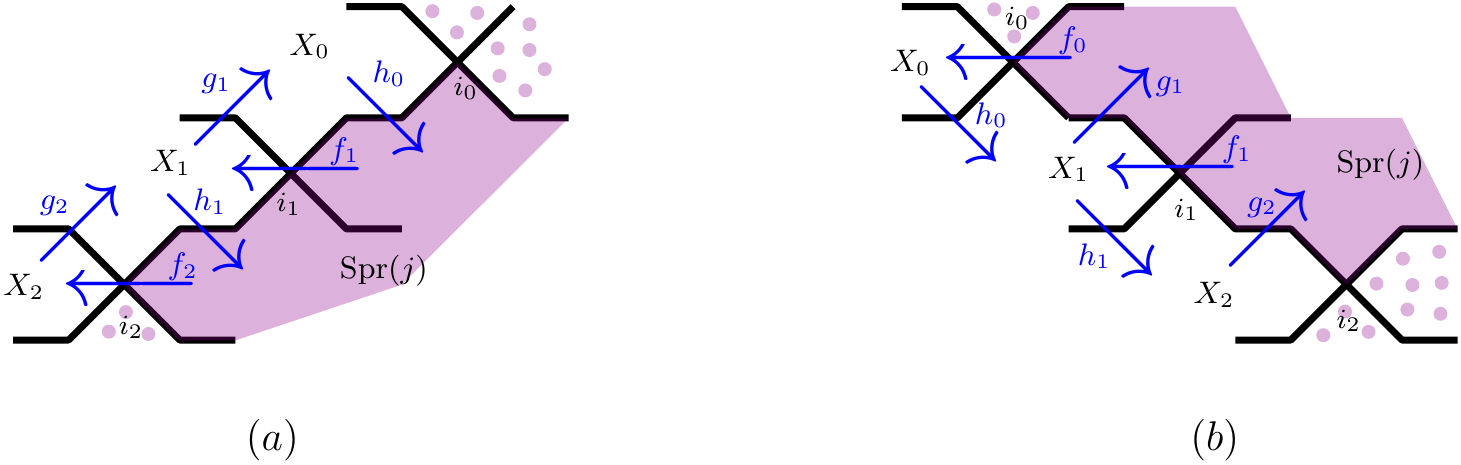}
	\caption{The boundary of $\jc{j}$ considered in \cref{lem:new1} and \cref{lem:new2}.}
	\label{Fig_new1}
\end{figure}

\begin{lemma}\label{lem:new1}
Suppose that the boundary of $\jc{j}$ is as in \cref{Fig_new1} (a), then the following statements regarding the morphisms in the figure hold.
\begin{itemize}
\item[(i)] $f_1=g_2f_2$ or $f_1=0$;
\item[(ii)] $h_1=h_0g_1$ or $h_1=0$.
\end{itemize}
In particular, the morphisms $f_1, h_1$ are reducible. 
\end{lemma}

\begin{proof}
Since irreducible morphisms are preserved under right multiplication, see \cref{lem:maps}, we may assume that $j=1$ and $M_j=S(j)$ is a simple module.   
We also remark that the maps in \cref{prop:maps} were defined for solid crossings. Nevertheless, in the situation of \cref{Fig_new1}(a) we will extend the definition of these maps to hollow crossings.  

Let $\alpha$ denote the rising strand bounding $\jc{j}$ on the left. Let $X_0, X_1, X_2$ denote certain modules appearing in the chambers to the left of crossings $i_0, i_1, i_2$ along $\alpha$.   More precisely,  $X_j$ is the unique indecomposable module in its chamber supported at vertex $j$ and if no such summand exists then we set $X_j=0$.  


(i) Suppose that $f_1\not=0$, that is $X_1$ has $S(j)$ in the socle.  If $i_2\in \jhol{v}$, then by \cref{lem:hollowRel} we have $X_1=X_2$, so $f_1=f_2$ and $g_2=1_{X_1}$ and the statement holds.  

If $i_2\in  \jsol{v}$ then it suffices to show that $S(j)$ that is in the socle of $X_2$ is mapped to the $S(j)$ in the socle of $X_1$ via $g_2$.  For this recall the construction of the morphism $g_{2}$ in \cref{prop:maps}(b).  The module $X_{1}$ is obtained from $X_{2}$ by removing a strip from the top of $X_{2}$ and adding a strip to the bottom starting at the vertical step labeled $b$, where $b$ is the right endpoint of the $v$-strand running just above the crossing $i_2$.  Therefore, the map $g_{2}$ has the desired property if and only if the modules $X_{2}$ and $X_1$ have the same dimension at vertex $j$ which occurs if and only if $b<j$. However, by \cite[Corollary IV.47]{gracie} we obtain that $b<j$ or $b=j+1$.  Note that $b=j+1$ implies that $X_1$ cannot have $S(j)$ in the socle, contradicting our assumption that $f_1\not=0$.  Thus, $b<j$, which completes the proof that $f_1$ factors through $f_2$ in the case $i_2\in  \jsol{v}$.

(ii) Suppose that $h_1\not=0$, that is $X_1$ has $S(j)$ in the top. If $i_1\in \jhol{v}$, then by \cref{lem:hollowRel} we have $X_1=X_0$, so $h_1=h_0$ and $g_1=1_{X_1}$ and the statement holds.  

If $i_1\in  \jsol{v}$ then it suffices to show that $S(j)$ that is at the top of $X_1$ is mapped to $S(j)$ in the top of $X_0$ via $g_1$.  It is easy to see that since $X_1$ has $S(j)$ at the top then so does $X_0$, so we obtain the desired conclusion since $g_1$ is induced by a surjective map on the chamber modules.
\end{proof}

\begin{lemma}\label{lem:new2}
Suppose that the boundary of $\jc{j}$ is as in \cref{Fig_new1} (b), then the following statements regarding the morphisms in the figure hold.
\begin{itemize}
\item[(i)] $f_1=h_0f_0$ or $f_1=0$;
\item[(ii)] $g_1=g_2h_1$ or $g_1=0$.
\end{itemize}
In particular, the morphisms $f_1, g_1$ are reducible. 
\end{lemma}

\begin{proof}
We proceed in the same way as in the proof of \cref{lem:new1}.  Hence, we may assume that $j=1$ and $M_j=S(j)$ is a simple module.   

Let $\alpha$ denote the falling strand bounding $\jc{j}$ on the left. Let $X_0, X_1, X_2$ denote certain modules appearing in the chambers to the left of crossings $i_0, i_1, i_2$ along $\alpha$.   More precisely,  $X_j$ is the unique indecomposable module in its chamber supported at vertex $j$ and if no such summand exists then we set $X_j=0$.  


(i) Suppose that $f_1\not=0$, that is $X_1$ has $S(j)$ in the socle.  If $i_0\in \jhol{v}$, then the chamber above $i_0$ is in $\jc{j}$ and by \cref{lem:hollowRel} we have $X_0=X_1$.  Hence, $f_1=f_0$ and $h_0=1_{X_1}$ and the statement holds.  

If $i_0\in  \jsol{v}$ then it suffices to show that $S(j)$ that is in the socle of $X_0$ is mapped to the $S(j)$ in the socle of $X_1$ via $h_0$.  Note that $X_0$ actually contains $S(j)$ in the socle, because $i_0\in  \jsol{v}$ so the map $f_0$ is injective.  Let $b'$ denote the right endpoint of the $v$-strand passing just below the crossing $i_0$. Then from the description of the map $h_0$ in \cref{prop:maps}(c), we can see that $h_0$ maps $S(j)$ in the socle of $X_0$ to the socle of $X_1$ if and only if $\textup{dim} (X_0)_j\leq \textup{dim} (X_{1})_j$.  This in turn occurs if $b'>j$.  Since this $v$-strand with endpoint $b'$ is part of a lower boundary of $\jc{j}$, then \cite[Corollary IV. 47]{gracie} implies that $b'\geq j$.  Note that if $b'=j$ then $X_{1}$ does not have $S(j)$ in the socle, which contradicts the assumption that $f_{1}\not=0$.  This completes the proof that $f_0$ factors through $f_1$ in the case $i_0\in  \jsol{v}$.

(ii) Suppose that $g_1\not=0$, that is $X_1$ has $S(j)$ in the top. If $i_1\in \jhol{v}$, then by \cref{lem:hollowRel} we have $X_1=X_2$, so $g_1=g_2$ and $h_1=1_{X_1}$ and the statement holds.  

If $i_1\in  \jsol{v}$ then it suffices to show that $S(j)$ that is at the top of $X_1$ is mapped to $S(j)$ in the top of $X_2$ via $h_1$.  It is easy to see that since $X_1$ has $S(j)$ at the top then so does $X_2$.
Let $a'$ denote the right endpoint of the rising $w$-strand at the crossing $i_1$.
By the description of the map $h_1$, we see that $h_1$ maps the $S(j)$ at the top of $X_1$ to the top of $X_{2}$ if and only if $\textup{dim} (X_1)_j\geq \textup{dim} (X_{2})_j$.  Moreover, the latter holds if $a'>j$, which in turn follows from \cite[Proposition IV.51]{gracie}.  
\end{proof}

\begin{figure}
\includegraphics{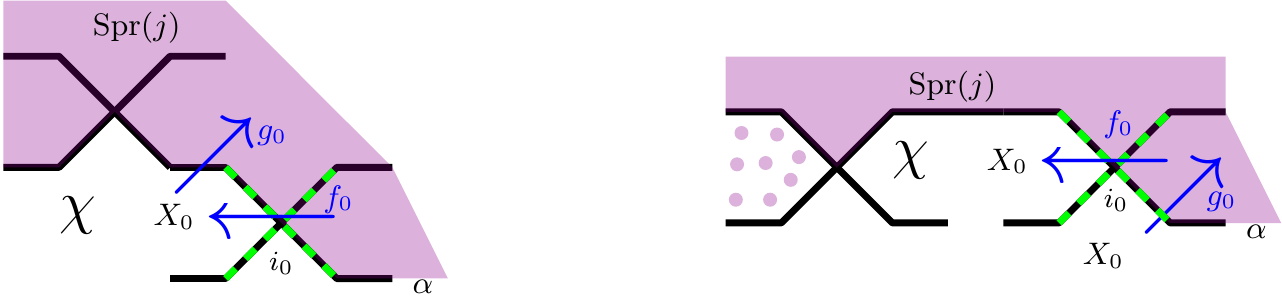}
\caption{The chamber $\ch$ considered in \cref{lem:new3}. The crossing $i_0$ is either hollow (as shown) or a cusp. If $i_0$ is a cusp, then we only consider the maps which begin or end in $\jc{j}$.}
\label{Fig_new2}
\end{figure}

\begin{lemma}\label{lem:new3}
Suppose $\ch$ is as in \cref{Fig_new2} and either $i_0$ is a cusp, or $i_0$ is hollow and traveling down the falling strand of $i_0$, one sees only hollow crossings before reaching a cusp of $\jc{j}$.
Then the morphisms $g_0$ and $f_0$ are reducible.
\end{lemma} 

\begin{figure}
\includegraphics[width=0.7\textwidth]{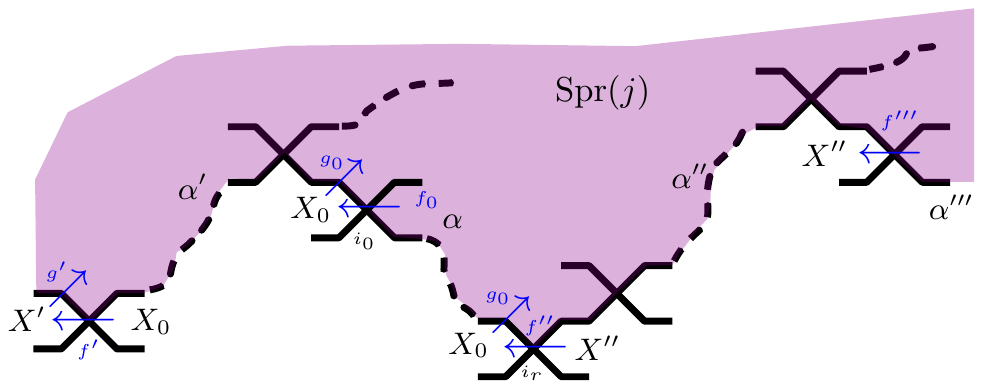}
\caption{The proof of \cref{lem:new3}. The strands $\alpha'$ and $\alpha''$ have right endpoints $a'$ and $ a''$.}
\label{Lem_5.6}
\end{figure}

\begin{proof}
Suppose that we are in the case of \cref{Fig_new2} shown on the left. The arguments in the other case follow similarly.  We proceed in the similar way as in the proof of \cref{lem:new1}.  Hence we may assume that $j=1$ and $M_j=S(j)$ is a simple module.  
Now, consider a more detailed picture of the situation given in \cref{Lem_5.6}.  Each $X_i$ in the figure denotes an indecomposable module supported at vertex $j$ appearing in a given chamber, and if no such module exists we let $X_i=0$.
Here we depict the more complicated case when the cusp at $i_r$ is as shown.  The shape of the cusp will not be important in the argument that $g_0$ is reducible.  Moreover, the other possibility is when the cusp is such that the only chamber around $i_r$ not contained in $\jc{j}$ is to the left of the crossing, and in this case we see that $f_0$ is reducible by \cref{lem:new1}(i).  Hence, we may consider \cref{Lem_5.6}. By assumption, either $i_0$ is a cusp or $i_r$ is a cusp and all crossings $i_0, \dots, i_{r-1}$ along the strand $\alpha$ and the boundary of $\jc{j}$ preceding $i_{r}$ are hollow.  In particular, this implies that all the maps $f_i, g_i$ and the indecomposable summands $X_i$ coming from these crossings are the same, so we denote them by $f_0, g_0$ and $X_0$ respectively.  By \cref{cor:rightEndHollow} it follows that all the crossings along the strands $\alpha', \alpha''$, as in the figure, along the right boundary of $\jc{j}$ except for the first and the last one are are also hollow.

Now, consider the map $g_0: X_0\to S(j)$. This map is nonzero as the crossing $i_{r}$ is solid.   Similarly, the map $g': X'\to S(j)$ is nonzero as well.  We claim that $g_0$ factors through $g'$.  Indeed, consider the injective morphism $f': X_0 \to X'$, and let $a'$ be the right endpoint of $\alpha'$.    By \cref{prop:maps}(a) we see that it maps $S(j)$ at the top of $X_0$ to the top of $X'$ if $a'<j$.  However, by \cite[Proposition IV.50]{gracie} $a'<j$ which shows that the composition $g'f': X_0\to X'\to S(j)$ is nonzero and so it equals $g_0$.  This completes the proof of the claim that $g_0$ factors through $g'$, and hence it is reducible. 

Now, consider the map $f_0: S(j)\to X_0$.  If $f_0=0$ then it is reducible and the lemma holds, hence suppose that $f_0\not=0$.  Consider the maps $f'': X''\to X_0$ and $f''': S(j) \to X''$ as shown in the figure.  Let $a''$ be the right endpoint of the strand $\alpha''$.  Observe that by the description of the map $f''$ in \cref{prop:maps}(a) the module $X''$ has $S(j)$ in the socle, because $f''$ is injective, $X_0$ has $S(j)$ in the socle, and $a''<j$.  This implies that the composition $f''f''': S(j) \to X''\to X_0$ is nonzero and equals $f_0$.  Hence, $f_0$ factors through $f'''$, so it is reducible.  This completes the proof of the lemma. 
\end{proof}


Next, we will show that the maps coming from the left boundary of $\jc{j}$ not covered by \cref{lem:new1}-\cref{lem:new3} are indeed irreducible.

\begin{proposition} \label{prop:Lec_arrows} Let $j \in \jsol{v}$ be a mutable solid crossing.
The complete list of irreducible morphisms involving $M_j$ and $M_d$ for $d>j$ are illustrated in \cref{fig:irredMaps}, with the following exception: in the rightmost two cases, if $d$ is hollow, and traveling down the falling strand of $d$, you pass through only hollow left ends of $\jc{j}$ and then reach a cusp, do not include the arrow for $d$ or the cusp. 
	\begin{figure}[h]
	\includegraphics[width=\textwidth]{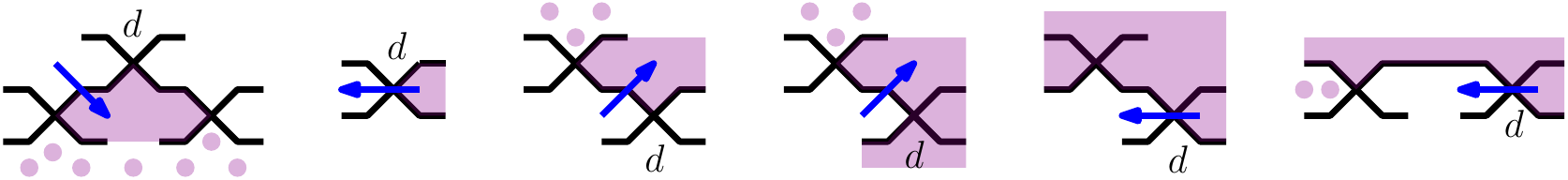}
	\caption{\label{fig:irredMaps} Arrows in the wiring diagram quiver corresponding to irreducible morphisms. The shaded chambers are in $\jc{j}$, white regions are not in $\jc{j}$, and dotted regions can be either. If $d$ is a hollow crossing in the rightmost 2 cases, $M_d$ should be interpreted as $M_{d'}$ where $d'$  is the first solid crossing you reach traveling down the falling strand of $d$.}
\end{figure}
\end{proposition}


\begin{figure}
\includegraphics[width=0.4\textwidth]{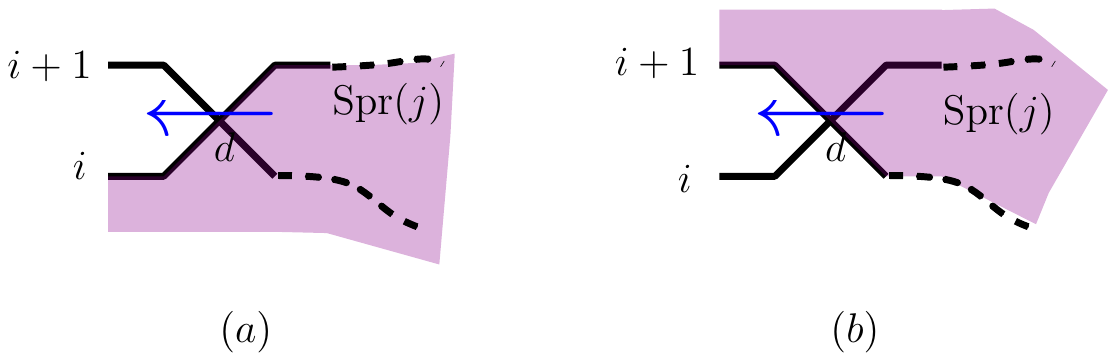}
\caption{ Arrows corresponding to irreducible morphisms if $\jc{j}$ is open on the left.}
\label{Fig:5.7new}
\end{figure}

\begin{proof}
	In fact, we will show something stronger than the proposition statement, as we also consider irreducible maps coming from the left boundary of $\jc{j}$ when $\jc{j}$ is not bounded. 
	In particular, we will show the proposition statement together with the following when $\jc{j}$ is open on the left:
	\begin{itemize}
		\item[(i)] In \cref{Fig:5.7new}(a) the arrow corresponds to an irreducible morphism except if traveling along the rising strand of $d$ you reach a peak as in \cref{fig:jcj} appearing in the first row second from the right.  
		\item[(ii)] In Figure~\ref{Fig:5.7new}(b) the arrow corresponds to an irreducible morphism except if traveling down the falling strand of $d$ you pass through only hollow left ends of $\jc{j}$ before reaching a cusp.
		\end{itemize}

As before we may assume $M_j=S(j)$. 
	
First, observe that the exceptional cases as in the statement of the proposition already follow from \cref{lem:new3}.  Moreover, in the exceptional situation of case (i) we see that the rising strand at $d$ actually has the right endpoint $j$, so in particular the corresponding morphism is actually zero since $M_j=S(j)$.  Similarly, the exceptional situation of case (ii) follows form \cref{lem:new3}.
Hence we may omit these special cases from the discussion.

Now, to prove the proposition we proceed by induction on the length of $w$.  If $\ell(w)=1$ and $v=e$ then there is only one nontrivial chamber and no irreducible morphisms. Otherwise, if $\ell(w)=1$ and $v=w$ then there are no nontrivial chambers and $M_j=0$, so again there is nothing to show. 

Let $\bw = s_{i_1} \dots s_{i_\ell}$ and $v \leq \bw$, and let $\bw' = s_{i_2} \dots s_{i_\ell}$ with $v' = v$ if $1\in \jsol{v}$ and $v'=s_{i_1}v$ otherwise. Let $j\in \jsol{v}$ and consider the corresponding region $\jc{j}$.  Note that if $j=1$ then this reduces to the base case above and the result follows.  Hence, we may suppose that $j>1$.  Let $\Spr'(j)$ denote the corresponding region in $\wiring{v', \bw'}$, and by induction we know that the irreducible maps coming from the left boundary of $\Spr'(j)$ are as in the statement of the proposition and parts (i) and (ii) given above.    We can think of the wiring diagram for $\wiring{v', \bw'}$ as being contained inside $\wiring{v, \bw}$, hence we can think of $\Spr'(j)$ as being contained inside $\jc{j}$ as the two only differ by at most one chamber $\chi_1$.  We consider several cases depending on whether the new chamber $\chi_1$ is in $\jc{j}$, adjacent to $\jc{j}$, or neither.   Recall that if $1\in \jsol{v}$ then $M_1$ denotes the unique indecomposable summand of $U_1$ that does not appear in the cluster tilting module $U_{v',\bw'}$.  

First, suppose that $\chi_1$ is not in $\jc{j}$ and not adjacent to $\jc{j}$.  Then $\jc{j}$ is the same as $\Spr'(j)$ and they have the same neighboring chambers.  If in addition $1\in \jhol{v}$ then the two modules $U_{v,\bw}$ and $U_{v',\bw'}$ are exactly the same in Leclerc's construction, so the irreducible morphisms remain the same.  Since we are in the case when $\jc{j}=\Spr'(j)$ then we are done by induction.  

Now, suppose that we are still in the situation when $\chi_1$ is not in $\jc{j}$ and not adjacent to $\jc{j}$ but $1\in \jsol{v}$. Then  $U_{v,\bw}$ contains a new indecomposable module $M_1$ that is not present in $U_{v',\bw'}$.  By \cref{prop:neigh_ch}, irreducible morphisms come from neighboring chambers, hence there are no irreducible morphisms between $M_j$ and $M_1$, since $\chi_1$ is not adjacent to $\jc{j}$ by assumption.  This implies that there are no new irreducible morphisms starting or ending at $M_j$ in the larger category $\text{add}\,U_{v,\bw}$.  Now suppose that a morphism $f: M_j\to M$ which was irreducible in $\text{add}\,U_{v',\bw'}$ now becomes reducible with the addition of the new module $M_1$.  Since every morphism is a composition of irreducible ones we have $f: M_j  \to X \to M$ for some irreducible map $f': M_j\to X$  in $\text{add}\,U_{v,\bw}$.  Then $X\not\cong M_1$ and moreover $f'$ must also be irreducible in $\text{add}\,U_{v',\bw'}$.  However, this means that the map $f$ is reducible in $\text{add}\,U_{v',\bw'}$, which is a contradiction to our assumption.   Note that analogous argument shows that any morphism $g: X\to M_j$ which was irreducible in $\text{add}\,U_{v',\bw'}$ also remains irreducible in $\text{add}\,U_{v,\bw}$.  This implies that the irreducible maps starting and ending in $M_j$ are the same in both $\text{add}\,U_{v,\bw}$ and $\text{add}\,U_{v',\bw'}$, and since $\jc{j}=\Spr'(j)$ the same conclusion holds on the level of wiring diagrams.  This completes the proof in the case when $\chi_1$ is not in $\jc{j}$ and not adjacent to $\jc{j}$.

Now, suppose that $\chi_1\in \jc{j}$.  From the representation theoretic point of view, we obtain no new morphisms starting or ending in $M_j$ in the larger category $\text{add}\,U_{v,\bw}$, and morphisms that were irreducible in $\text{add}\,U_{v',\bw'}$ still remain irreducible  in $\text{add}\,U_{v,\bw}$.  Indeed, if $1\in \jhol{v}$ then the two seeds contain the same set of modules, and if $1\in \jsol{v}$ then $M_1$ and $M_j$ have disjoint support so there are no nonzero morphisms between these modules.   Now looking at the combinatorics of the wiring diagrams, we see that if $\chi_1$ is in the interior of $\jc{j}$ then no new arrows appear.  If $\chi_1$ lies on the boundary of $\jc{j}$ then we are in one of the situation of \cref{Prop_5.7new}(a)-(c).  Note that in these cases the left boundary of $\jc{j}$ is the same as the left boundary of $\Spr'(j)$, hence we also get no new maps on the level of wiring diagrams.

\begin{figure}
\includegraphics[width=0.6\textwidth]{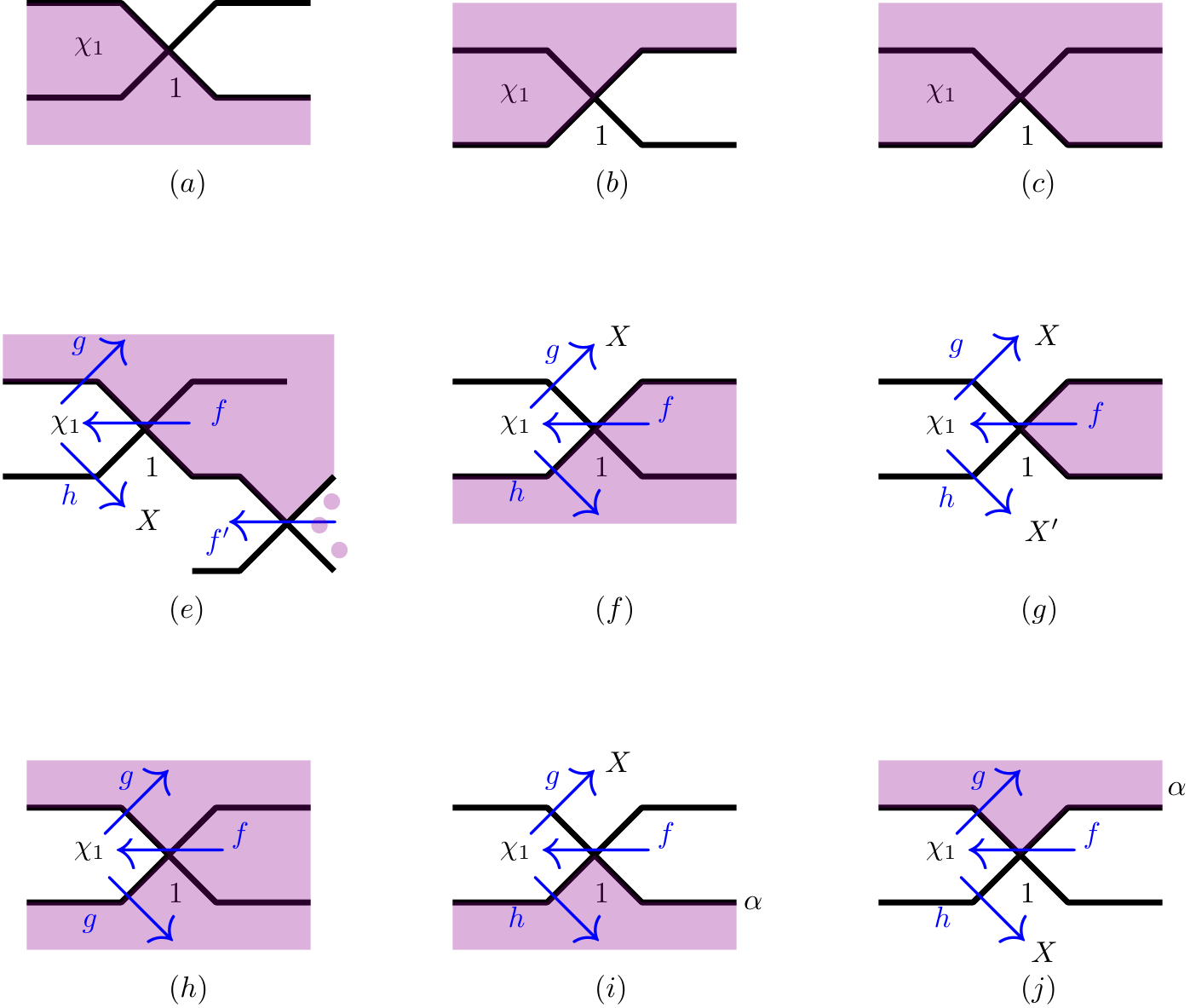}
\caption{Cases in the proof of \cref{prop:Lec_arrows}. Shaded chambers are in $\jc{j}$.}
\label{Prop_5.7new}
\end{figure}

The remaining cases is when $\chi_1$ is not contained in $\jc{j}$ but it is adjacent to $\jc{j}$, which are shown in \cref{Prop_5.7new}(e)-(j).  If $1\in\jhol{v}$  then  we are in the case of \cref{Prop_5.7new}(e),(f) and it is easy to see that the proposition holds by induction, as these correspond to cases (i) and (ii) when $\jc{j}$ is open to the left.  Therefore, we may assume that $1\in \jsol{v}$.   



First, consider the case of \cref{Prop_5.7new}(e), where $\chi_1$ is adjacent to the lower boundary of $\jc{j}$.
We see that the new maps in $\text{add}\,U_{v,\bw}$ starting and ending in $M_j$ due to the addition of the new chamber $\chi_1$ are $f, g$.  By \cref{lem:new2} the map $g$ is reducible and the map $f'$ factors through $f$.  Note that by induction $f'$ was irreducible in $\text{add}\,U_{v',\bw'}$ and now we need to show that $f$ is irreducible in $\text{add}\,U_{v,\bw}$ and that no other maps starting or ending in $M_j$ other than $f'$ become reducible.  If $f: M_j\to M_1$ is reducible then there is a nontrivial factorization $f=f''f''': M_j\xrightarrow{f'''} M\xrightarrow{f''} M_1$ where $f''$ is irreducible in $\text{add}\,U_{v,\bw}$.  Then by \cref{prop:neigh_ch} the map $f''$ comes from an arrow between the neighboring chambers.  The only arrow ending in $M_1$ is $f$, so $f''=f$, and we obtain a contradiction to the existence of the nontrivial factorization of $f$.  This shows that $f$ is indeed irreducible in $\text{add}\,U_{v,\bw}$.  Similarly, irreducible maps in $\text{add}\,U_{v',\bw'}$ cannot factor through $M_1$ so they remain irreducible in   $\text{add}\,U_{v,\bw}$. 
This completes the proof in the case of \cref{Prop_5.7new}(e).  Analogous argument applies in the cases of \cref{Prop_5.7new} (f) and (h) where $\chi_1$ is adjacent to the upper boundary of $\jc{j}$ or both.

Now, consider the case of \cref{Prop_5.7new}(h), where all the neighboring chambers of $\chi_1$ are in $\jc{j}$. Here we get new maps $g,f$ in $\text{add}\,U_{v,\bw}$.  Since these are the only maps starting or ending at $M_1$, then they must be irreducible by \cref{prop:neigh_ch}. Therefore, it suffices to show that an irreducible map in $\text{add}\,U_{v',\bw'}$ remains irreducible in $\text{add}\,U_{v,\bw}$.  Indeed, if not then it would have to factor through $gf: M_j\to M_j$ which is a map in $\text{add}\,U_{v',\bw'}$, a contradiction.  This completes the proof in the case of \cref{Prop_5.7new}(h). 


Next, consider the case of \cref{Prop_5.7new}(i), and we claim that the map $h$ is irreducible in $\text{add}\,U_{v,\bw}$.  If $h$ is reducible then it would have to factor though $g$.  By assumption $M_j=S(j)$, so the strand $\alpha$ has right endpoint $j$ and $X$ is not supported in $j$ since it lies in a chamber above $\alpha$.  This means that $g=0$.   Hence, we see that $h$ does not factor through $g$. Finally, by similar arguments as before, we can see that irreducible maps in $\text{add}\,U_{v',\bw'}$ remain irreducible in $\text{add}\,U_{v,\bw}$ since any such map would have to factor though $f$ which starts in a module in $\text{add}\,U_{v',\bw'}$.   This completes the proof in this case. 


Now, consider the remaining case of \cref{Prop_5.7new}(j), which is similar to the previous situation. Note that by the same reasoning as before irreducible maps in $\text{add}\,U_{v',\bw'}$ remain irreducible in $\text{add}\,U_{v,\bw}$.   We claim that $g$ is irreducible.  If $g$ is not irreducible than it factors though $h$, and recall that $M_j=S(j)$.  By \cite[Proposition IV.50]{gracie}, the right endpoint $a'$ of the strand $\alpha$ satisfies $a'<j$ or $a'=j+1$.  If $a'=j+1$ then $X$ is not supported in $j$, so $g$ cannot factor though $h$.  If $a'<j$ then we consider the structure of $M_1$ and $X$ and the morphism $h$ as described in \cref{prop:maps}(c).  We see that in the case $a'<j$ the image of $S(j)$ in the top of $M_1$ under $h$ does not map into the top of $X$.  This means that $g$ cannot factor though $h$, which shows the claim that $g$ is irreducible.    This completes the proof of the proposition.  
\end{proof}

\section{Correspondence between quivers}\label{sec:quiversEqual}

	In this section, we first analyze $\iquiv{v, \bw}$. Then in \cref{sec:quiversEqualProof}, we show that $\iquiv{v, \bw}=\lquiv{v, \bw}$.
	
	\subsection{Ingermanson's quiver}
Throughout this section, we fix $v \leq \bw$ with $\bw$ unipeak.

		The definition of Ingermanson's quiver (c.f. \cref{def:IngQuiv}) involves a lot of cancellation. In this section, we give a ``cancellation-free" description of Ingermanson's quiver in terms of the wiring diagram quiver, so that we can compare with \cref{prop:Lec_arrows}.
	
	\begin{defn} Let $S$ and $T$ be disjoint subsets of the cluster of $\Sigma_{v, \bw}$. A collection $\mcc$ of arrows in the wiring diagram is a \emph{witnessing collection} for $(S,T)$ if for every $A_s \in S$ and $A_t \in T$, the number of arrows $A_s \to A_t$ in $\iquiv{\vw}$ is the same as the number of arrows in $\mcc$ which point from a chamber in $\jc{s}$ to a chamber in $\jc{t}$ (this number may be negative).
	
	That is to say, one can compute the arrows in $\iquiv{v, \bw}$ between $S$ and $T$ just by considering the contributions of $\mcc$ and ignoring all other arrows in the wiring diagram.
	\end{defn}

Let $A_d$ be a mutable cluster variable and let 
\[\mcs_d=\{A_i: i \in [d-1], A_i \text{ does not appear in }\jc{d}.\}\] In this section, we will find a witnessing collection $\mcc_d$ for $(A_d, \mcs_d)$. The union of the witnessing collections will clearly be equal to Leclerc's quiver $\lquiv{\vw}$; we will later show that the union of witnessing collections is equal to Ingermanson's quiver $\iquiv{\vw}$ as well.

Before finding a witnessing collection for $(\ivar{d}, \mcs_d)$, we present an alternate definition of $\iquiv{\vw}$.

\begin{defn}[Crossing monomial]
	Let $c$ be a crossing in $\wiring{\vw}$. Say $\ch_{c^{\uparrow}}, \ch_{c^{\downarrow}}, \ch_{c^{\leftarrow}}, \ch_{c^{\rightarrow}}$ are the chambers above, below, to the left, and to the right of $c$, respectively. The \emph{crossing monomial} of $c$ is defined as
	\[t_c:= \frac{\lmin{c^\uparrow} \lmin{c^\downarrow}}{\lmin{c^\leftarrow} \lmin{c^\rightarrow}}.\]
\end{defn}

\begin{rmk} \label{rmk:crossingMonoFacts}\ 
	\begin{enumerate}
		\item \label{itm:crossingMonoOne} If $c \in \jhol{v}$, then $t_c =1$ (c.f. \cref{lem:hollowRel}).
		\item \label{itm:crossingMonoYHat} Say a chamber $\ch$ is bounded on the left by crossing $a$ and on the right by crossing $b$. It is not hard to check that 
		\begin{equation} \label{eq:yHatinCrossing}
			\hat{y}^W_\ch=\frac{t_a}{t_b}.
		\end{equation}
	\end{enumerate}
\end{rmk}

Recall from \cref{def:ends} the notion of left ends, right ends, and cusps of $\jc{j}$.

\begin{lem}\label{lem:onlyLeftEnds}
	Let $c$ be a solid crossing. Then 
	\[\hat{y}_c = \frac{1}{t_c} \cdot \prod_{d} t_d\]
	where the product is over all $d$ that are left ends of $\jc{c}$.
\end{lem}

\begin{proof}
	Combine \cref{eq:yHats,eq:yHatinCrossing} and note that all crossing monomials for crossings in the interior of $\jc{c}$ cancel. Then apply \cref{rmk:crossingMonoFacts,cor:rightEndHollow}.
\end{proof}

\begin{thm}\label{thm:witnessingCollections}
	Consider $c \in \jsol{v}$ a mutable solid crossing. To obtain a witnessing collection $\mcc_c$ for $(A_c, \mcs_c)$, take the arrows indicated in \cref{fig:specialArrows}, with the following exception: in the final two cases, if $d$ is hollow and traveling down the falling strand of $d$, you pass through only hollow left ends of $\jc{c}$ and then reach a cusp, do not include the arrow for $d$ or the cusp.
	\begin{figure}[h]
		\includegraphics[width=\textwidth]{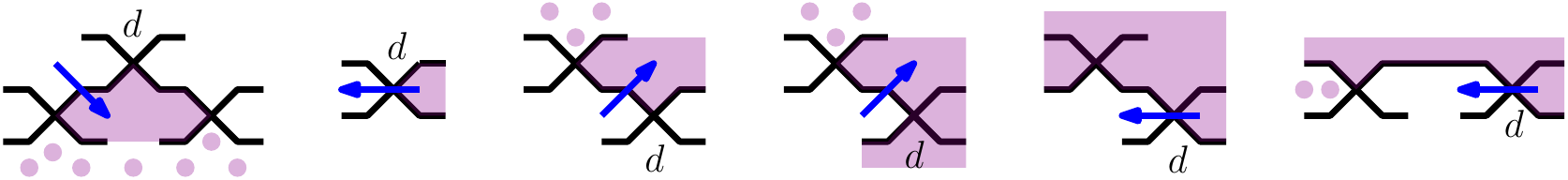}
			\caption{\label{fig:specialArrows} A witnessing collection for all arrows between $A_c$ and $A_d$ ($d<c$) in $\iquiv{v, \bw}$. If $d$ is hollow $A_d$ should be interpreted as $A_{d'}$ where $d'$ is the first solid crossing along the falling strand of $d$. The shaded chambers are in $\jc{c}$, white regions are not in $\jc{c}$, and dotted regions can be either.}
	\end{figure}

\end{thm}

Comparing \cref{thm:witnessingCollections,prop:Lec_arrows}, we have an immediate corollary.

\begin{cor}\label{cor:unionOfWitnessingEqualsLeclerc}
	The union of witnessing collections in Ingermanson's quiver
	\[\bigcup_{c \in \jsol{v}} \mcc_c\]
	is equal to Leclerc's quiver $\lquiv{\vw}$ (where vertices of both quivers are labeled by solid crossings). In particular, $\lquiv{\vw}$ is a subquiver of $\iquiv{\vw}$.
\end{cor}

\begin{figure}
	\includegraphics[width=0.7\textwidth]{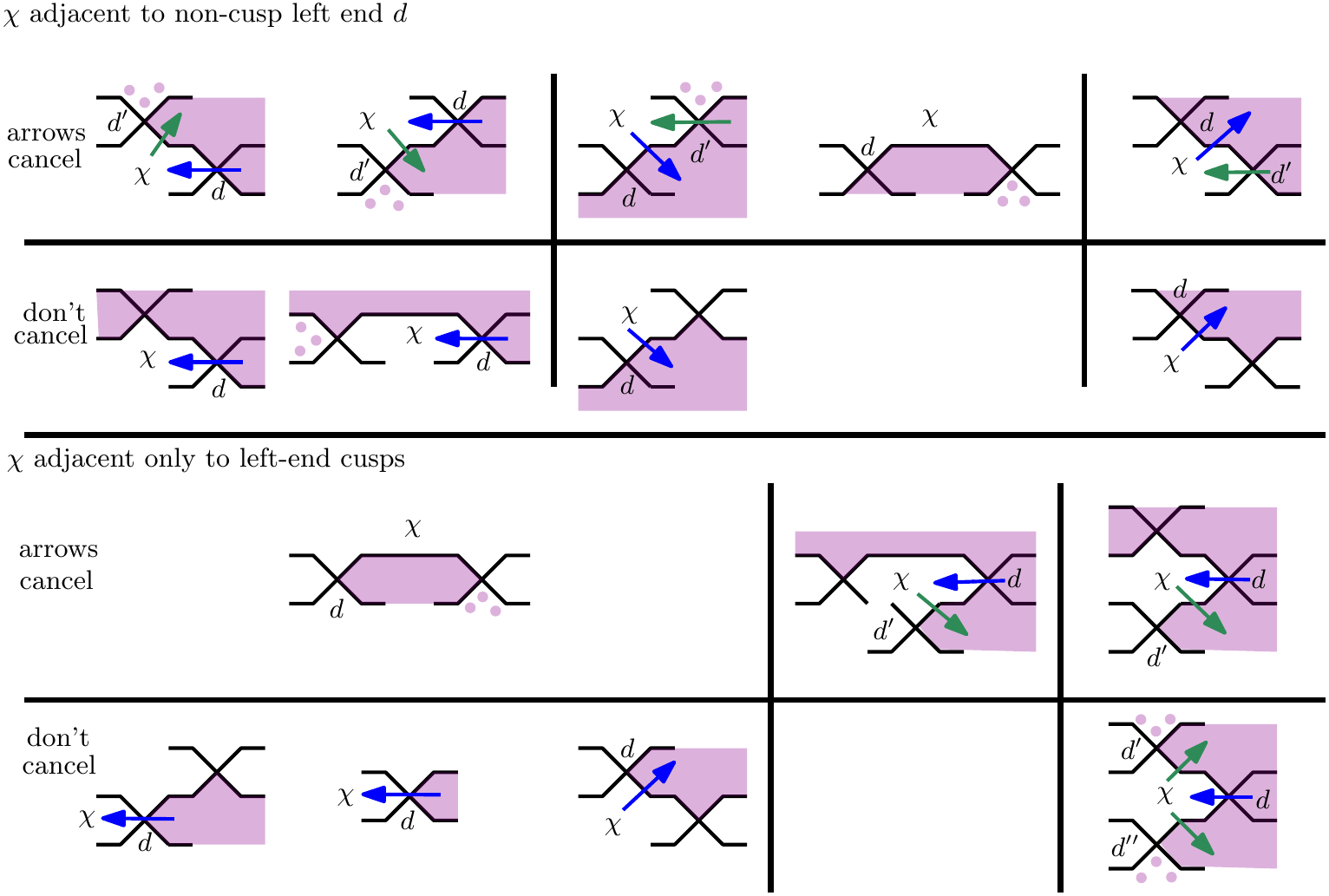}
	\caption{\label{fig:arrows-cancel-cases} The cases showing when left-end-arrows to chamber $\ch$ cancel or don't cancel in the proof of \cref{thm:witnessingCollections}. }
	\label{fig:jcj}
\end{figure}

\begin{proof}[Proof of \cref{thm:witnessingCollections}] Using \cref{lem:IRightStable}, without loss of generality we may assume $c=\ell$ is the final crossing. 
	
	By \cref{lem:onlyLeftEnds}, the arrows from $A_c$ to $A_a$ for $a<c$ are determined by product of crossing ratios for left ends of $\jc{c}$. Recall that we only concern ourselves with $A_a$ not appearing in $\jc{c}$, so in fact we can consider the product of modified crossing ratios
		\[t'_d:= \frac{\lmin{c^\uparrow} \lmin{c^\downarrow}}{\lmin{c^\leftarrow}}.\]
		The modified crossing ratio encodes the arrows 
		\begin{center}
			\includegraphics[width=0.1\textwidth]{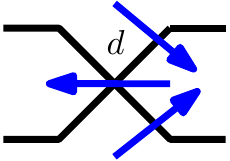}
		\end{center}
	between chambers around the left end $d$, which we call $d$-arrows.
		
	First, we analyze which chamber minors cancel in the product of modified crossing ratios, or, equivalently, which arrows around left ends cancel. Again, we only care about chambers which are not in $\jc{c}$. Suppose a chamber $\chi$ not in $\jc{c}$ contributes to $t'_d$. Whether or not the $d$-arrow involving $\chi$ cancels with another $d'$-arrow depends on if there is a nearby left end along one of the strands of $d$. The cases are summarized in \cref{fig:arrows-cancel-cases}. 
	
	The situations in which the $d$-arrow involving $\chi$ does not cancel is exactly when $\chi$ and the arrow are as pictured in \cref{fig:specialArrows}. We call such chambers $\chi$ and arrows \emph{special}. 
	
	Now, we show that each special chamber $\chi$ contains a unique cluster variable which is not in $\jc{c}$, which we denote $A_\chi$. Let $a$ be the crossing to the right of $\chi$. In the four leftmost cases of \cref{fig:specialArrows}, $a$ is solid, so $A_a$ appears in $\chi$. On the other hand, in all cases but the one in the far left of \cref{fig:specialArrows}, a chamber of $\jc{c}$ lies to the right or above $a$. Using \cref{lem:propogateRightAndUp}, all other cluster variables $A_r$ appearing in $\chi$ appear in an adjacent chamber in $\jc{c}$. Thus $A_a$ is the only candidate for $A_\chi$. 
	
	If we are in the far left case of \cref{fig:specialArrows}, the falling strand $\alpha$ of $a$ is also the falling strand of crossing $c$. Indeed, following $\alpha$ to the right of $a$, by unipeakness $\alpha$ continues to travel down and eventually leaves the boundary of $\jc{c}$, say immediately after crossing $b$. The crossing $b$ must be a cusp, and is thus solid. \cref{cor:forbidden-cusp} implies that $b$ is in fact a right end of $\jc{c}$. Then by \cref{cor:rightEndHollow}, $b=c$. Also, all crossings $a_1, \dots,a_k$ along $\alpha$ between $a$ and $b$ are right ends of $\jc{c}$ and so are hollow by \cref{cor:rightEndHollow}. So the product of crossing ratios $t_{a_1} \dots t_{a_k}$ is equal to $1$. We also have 
	\[t_{a_1} \dots t_{a_k}= \frac{\lmin{a_1^\uparrow} \lmin{a_k^\downarrow}}{\lmin{a_1^\leftarrow}\lmin{a_k^\rightarrow}}\]
	where $\chi_{a_k^\downarrow}$ is the chamber below $a_k$, etc.
	The chamber $\chi_{a_k^\rightarrow}$ is also the chamber above crossing $c$. Because $c=\ell$ is the last crossing, the chamber minor of this chamber is equal to $1$. The chamber $\chi_{a_1^\leftarrow}$ is in $\jc{c}$, so the above equality implies in particular that all cluster variables appearing in $\chi_{a_1^\uparrow}$ also appear in $\jc{c}$. By \cref{lem:propogateRightAndUp}, all cluster variables $A_r \neq A_a$ appearing in $\chi$ also appear in $\chi_{a_1^\uparrow}$, which is the chamber to the right of $a$. So in this case also, the only candidate for $A_\chi$ is $A_a$.
	
	For the two rightmost cases of \cref{fig:specialArrows}, the crossing $a$ may or may not be solid. If it is solid, the same argument as above shows that $A_a$ is the only candidate for $\chi$. If it is hollow, follow the falling strand $\alpha$ of $a$ to the right of $a$ until it hits a solid crossing $b$. This solid crossing is either a left end of $\jc{c}$ or a cusp, and is guaranteed to exist because $\alpha$ must leave the boundary of $\jc{c}$ eventually. By \cref{rmk:crossingMonoFacts}(1) and the fact that all crossings along $\alpha$ between $a$ and $b$ are hollow, $A_b$ appears in $\chi$ and all other cluster variables appearing in $\chi$ also appear in $\jc{c}$. So the only candidate for $A_\chi$ is $A_b$. Note that if $b$ is a cusp, then $A_b$ is also the candidate for $A_{\chi'}$ where $\chi'$ is the special chamber to the left of the cusp. The special arrows involving $\chi$ and $\chi'$ contribute a two-cycle between $A_b$ and $A_c$ in $\iquiv{\bv, w}$. So we exclude this case from consideration.
	
	Now, we need to verify that the candidate $A_a$ or $A_b$ for $A_\chi$ does not appear in $\jc{c}$ somewhere else. So long as we are not in the exception described in the theorem, this follows from \cref{thm:factorizationsAgree} and \cref{prop:Lec_arrows}, as there is an irreducible map involving the modules $M_c$ and $M_a$ (or $M_b$). This implies the contents of $M_c$ and $M_a$ (or $M_b$) overlap, so they cannot be summands of the same chamber module.
	
	The special arrows are the only arrows of the wiring diagram quiver which can contribute an arrow between $\mcs_c$ and $A_c$ in $\iquiv{\bv, w}$. We have already identified pairs of special arrows which give a 2-cycle in $\iquiv{\bv, w}$, and do not include these arrows in our collection $\mcc_c$. For all other pairs of special chambers $\chi, \chi'$, it is easy to see that $A_{\chi} \neq A_{\chi'}$, and so no arrows in $\iquiv{\bv, w}$ coming from the special arrows for $\chi$, $\chi'$ can cancel. This shows $\mcc_c$ is indeed a witnessing collection.
	
 \end{proof}
	
	\subsection{Equality of quivers} \label{sec:quiversEqualProof}
	We have shown that each arrow of Leclerc's quiver $\lquiv{\vw}$ is an arrow of Ingermanson's quiver $\iquiv{\vw}$. Moreover, we know that Leclerc's cluster algebra $\mca(\lseed{\vw})$ is a subalgebra of $\bbc[\rich{v,w}]$ and Ingermanson's cluster algebra $\mca(\iseed{\vw})$ is equal to $\bbc[\rich{v,w}]$. Finally, we have an automorphism $\tau_{v, w}^*$ of $\bbc[\rich{v,w}]$ which takes Ingermanson's cluster $\bA$ to Leclerc's cluster $\bB$. We will now use these facts to show that in fact, Ingermanson's quiver cannot have any additional arrows.

\begin{lem}\label{lem:quiverEqual}
	Consider two seeds $(\bx, Q)$ and $(\bx, Q')$ in the field of rational functions $\bbc(\bx)$ with the same cluster but different quivers. Let $\mca:=\mca(\bx, Q)$ be the cluster algebra of the first seed and $\mca':=\mca(\bx, Q')$ be the cluster algebra of the second.
	
	Suppose that $\mca' \subset \mca$ and that $Q'$ is a subquiver of $Q$ (identifying a vertex of $Q'$ with the vertex of $Q$ labeled by the same cluster variable). Then in fact $Q=Q'$.
\end{lem}	

\begin{proof}
	Consider a mutable vertex $k$ of $Q$ and $Q'$. We will argue that for all $j$, $\#(\text{arrows }k \to j)$ is the same in $Q$ and $Q'$.
	
	In $(\bx, Q')$, mutating at $k$ gives the cluster variable 
	\[\tilde{x}_k'= \frac{M_+ + M_-}{x_k}\] 
	where $M_+, M_-$ are monomials in $\bx$. Since $Q'$ is a subquiver of $Q$, mutating at $k$ in $(\bx, Q)$ gives the cluster variable 
	 	\[\tilde{x}_k= \frac{M_+N_+ + M_-N_-}{x_k}.\] 
	 	
	 By assumption, $\mca' \subset \mca$, so in particular $\tilde{x}'_k$ is a Laurent polynomial in $\mu_k(\bx, Q)$, say
	 \[\frac{M_+ + M_-}{x_k}= \sum_{\textbf{a}\in \bbz^n} c_{\textbf{a}} x_1^{a_1} \dots \tilde{x}_k^{a_k} \dots x_n^{a_n}.\]
	 Clearing the denominator on the left and writing $\tilde{x}_k$ in terms of the cluster $\bx$, we have
	 \[M_+ + M_-=x_k \cdot \sum_{\textbf{a}\in \bbz^n} c_{\textbf{a}} x_1^{a_1} \dots ( M_+N_+ + M_-N_- )^{a_k} x_k^{-a_k} \dots x_n^{a_n}.\]
	 	Note that the left hand side is a polynomial with degree $0$ in $x_k$, so the same must be true of the right hand side. Because the binomial $M_+N_+ + M_-N_-$ is also degree 0 in $x_k$, this implies that $a_k=1$ for all nonzero $c_\textbf{a}$. So we have
	 	\[M_+ + M_-= ( M_+N_+ + M_-N_- ) \cdot \sum_{\substack{\textbf{a}\in (\bbz_{\geq 0})^n \\ a_k=1} }c_{\textbf{a}} x_1^{a_1} \dots x_{k-1}^{a_{k-1}} x_{k+1}^{a_{k+1}} \dots x_n^{a_n}.\]
	 	Comparing degrees, we see that $a_1=\dots=a_n=0$ and $N_+=N_-=1$. That is, $\tilde{x}'_k=\tilde{x}_k$, which implies $\#(\text{arrows }k \to j)$ is the same in $Q$ and $Q'$ for all $j$.
\end{proof}

\begin{cor}\label{cor:quivEqual}
	Choose $v \leq \bw$ with $\bw$ unipeak. Label the vertices of both Ingermanson's quiver $\iquiv{\vw}$ and Leclerc's quiver $\lquiv{\vw}$ by the set of solid crossings $\jsol{v}$. Then $\iquiv{\vw}=\lquiv{\vw}$.
\end{cor}

\begin{proof}
	This follows directly from \cref{lem:quiverEqual}, with $(\bx, Q)=(\iclus{\vw}, \iquiv{\vw})=:\iseed{\vw}$ and $(\bx, Q')=(\lclus{\vw} \circ \tau_{v,w}, \lquiv{\vw})=:\tau_{v,w}^*(\lseed{\vw})$. Indeed, \cref{thm:Gracie-clusterstruc} and \cref{Lec-seed} together imply that $\mca(\lseed{\vw})$ is a subset of $\mca(\iseed{\vw})= \bbc[\rich{v, w}]$. Because $\tau_{v,w}$ is an automorphism of $\mca(\iseed{\vw})$, we also have that $\mca(\tau_{v,w}^*(\lseed{\vw})) \subset \mca(\iseed{\vw})$. By \cref{thm:varCorrespondence}, $\tau_{v,w}^*(\lvar{c})=\ivar{c}$, so the clusters of $\iseed{\vw}$ and $\tau_{v,w}^*(\lseed{\vw})$ are equal. \cref{cor:unionOfWitnessingEqualsLeclerc} shows that $\lquiv{\vw}$ is a subquiver of $\iquiv{\vw}$. So the assumptions of \cref{lem:quiverEqual} are satisfied.
\end{proof}

\section{Finishing up proofs}\label{sec:LSeedsMutationEquiv}

	In this section, we complete the proofs of \cref{thm:LecConjIsTrue,thm:compareStructures,thm:posPartLeclerc}.
	
	\begin{proof}[Proof of \cref{thm:compareStructures}]
		The equality of quivers is \cref{cor:quivEqual}. The fact that for all $c \in \jsol{v}$, $A_c = B_c \circ \tau_{v,w}$ is part $1)$ of  \cref{thm:varCorrespondence}.
	\end{proof}
	
	To prove \cref{thm:LecConjIsTrue}, we need the following proposition. 
	
	\begin{prop}\label{prop:LecSeedMutationEquiv}
		Let $v\leq w$ and fix two reduced words $\bw$ and $\bw'$ for $w$. The seeds $\lseed{v, \bw}$ and $\lseed{v, \bw'}$ are related by a sequence of mutations.
	\end{prop}
\begin{proof}
		Recall that all reduced expressions for $w$ are related by a sequence of commutation moves $s_i s_j \leftrightarrow s_j s_i$ where $|i-j|>1$ and braid moves $s_i s_{i+1} s_i \leftrightarrow s_{i+1} s_i s_{i+1}$. So it suffices to consider the case where $\bw, \bw'$ are related by a single commutation move or a single braid move.
		
		If $\bw$ and $\bw'$ are related by a commutation move, it's easy to see that the right chamber minors of $\wiring{v, \bw}$ are the same as the right chamber minors of $\wiring{v, \bw'}$, up to an indexing change. This means that the indecomposable summands of the cluster tilting objects $U_{\vw}$ and $U_{v, \bw'}$ are the same (though the summands may be indexed differently), and thus the seeds $\lseed{v, \bw}$ and $\lseed{v, \bw'}$ are the same.
		
		If $\bw$ and $\bw'$ are related by a braid move involving at least one hollow crossing, then the seeds $\lseed{v, \bw}$ and $\lseed{v, \bw'}$ are also related by reindexing. Let $(\rmin{d})'$ denote the right chamber minors for $\wiring{v, \bw'}$, and also denote the cluster variables of $\lseed{v, \bw'}$ with primes.  Say the braid move occurs at indices $c-1, c, c+1$. By inspection, $\rmin{d}= (\rmin{d})'$ for all crossings $d \notin \{c-1, c, c+1\}$, so in particular $B_d= B_d'$ for $d \in \jsol{v} \setminus \{c-1, c, c+1\}$. It is not hard to check that $\rmin{c-1}=(\rmin{c})'$ and $\rmin{c}=(\rmin{c-1})'$ by comparing the $v$ and $w$-strands of $\wiring{\vw}$ and $\wiring{v, \bw'}$; the key insight is that the same number of crossings in $\{c-1, c, c+1\}$ will be hollow in each wiring diagram. Now, by considering the various cases of which of $\{c-1, c, c+1\}$ can be hollow in $\wiring{\vw}$, one can see that if $d \in \{c-1, c, c+1\}$ is solid, then $B_d$ is an irreducible factor of either $\rmin{c}$ or $\rmin{c-1}$.  Since each of those chamber minors is also a chamber minor for $\wiring{v, \bw'}$, this means that $B_{d}$ is also a cluster variable in $\lseed{v, \bw'}$. Thus the indecomposable summands of $U_{\vw}$ and $U_{v, \bw'}$ are the same (though summands are indexed differently), and the corresponding seeds are the same.

		If $\bw$ and $\bw'$ are related by a braid move involving three solid crossings, say at indices $c-2, c-1, c$, then $\lseed{v, \bw}$ and $\lseed{v, \bw'}$ are the same or are related by mutation in direction $c$. As above, use primes to indicate the right chamber minors and cluster variables for $v \leq \bw'$. Again, by inspection, $\rmin{d}= (\rmin{d})'$ for $d \notin \{c-2, c-1, c\}$, $\rmin{c-1}=(\rmin{c-2})'$ and $\rmin{c-2}=(\rmin{c-1})'$. So the irreducible factors of 
		\[P=\prod_{c \neq d} \rmin{d}\]
		 coincide with the irreducible factors of 
		 \[Q=\prod_{c \neq d} (\rmin{d})'.\]
		  Now, if $\rmin{c}P$ and $(\rmin{c})'Q$ have the same set of irreducible factors, then $U_{\vw}$ and $U_{v, \bw'}$ have the same set of indecomposable summands and the corresponding seeds are equal. If they do not have the same set of factors, then in fact their factors differ only by $B_c \neq B_c'$. This implies that 
		  \[(T_{\vw}/M_c) \oplus M_c'=T_{v, \bw'}\]
		  where
		  \[T_{\vw}:= \bigoplus_{d\in \jsol{v}} M_d \qquad T_{v, \bw'}:= \bigoplus_{d\in \jsol{v}} M'_d\]
		  are the basic cluster tilting modules obtained from $U_{\vw}$ and $U_{v, \bw'}$ and $M_d, M_d'$ are the indecomposable modules which map to $B_d, B_d'$ under the cluster character map $\varphi$. The definition of mutation of basic cluster tilting objects (see e.g. \cite[Definition 3.9 c)]{Leclerc}) implies that $T_{\vw}$ and $T_{v, \bw'}$ are related by categorical mutation. This in turn implies $\lseed{\vw}$ and $\lseed{v, \bw'}$ are related by mutation in direction $c$.
	
\end{proof}

\begin{proof}[Proof of \cref{thm:LecConjIsTrue}]
	That $\mca(\lseed{v, \bw})$ does not depend on the choice of reduced expression $\bw$ is exactly \cref{prop:LecSeedMutationEquiv}.
	
	When $\bw$ is unipeak, \cref{thm:compareStructures} implies that $\mca(\iseed{v, \bw})= \tau_{v,w}^* (\mca(\lseed{v, \bw}))$. From \cref{thm:Gracie-clusterstruc}, we have $\mca(\iseed{v, \bw})= \bbc[\rich{v, w}]$. By \cref{prop:twist}, $\tau_{v, w}^*$ is an automorphism of $\bbc[\rich{v, w}]$, so
	\[\bbc[\rich{v, w}]= (\tau_{v, w}^*)^{-1}(\bbc[\rich{v, w}])= (\tau_{v, w}^*)^{-1}(\mca(\iseed{v, \bw}))= \mca(\lseed{v, \bw}).\]
\end{proof}

\begin{proof}[Proof of \cref{thm:posPartLeclerc}]
	The \emph{totally positive part} $\rich{v, w}^{>0, \Lus}$ of $\rich{v, w}$ defined by Lusztig is the subset of $\rich{v, w}$ where all the Marsh-Rietsch chamber minors are positive \cite[Proposition 12.1]{MR}. The left chamber minors $\lmin{c}$ are monomially related to the Marsh-Riestch chamber minors, and the cluster variables $\iclus{v, \bw}$ are in turn monomially related to the left chamber minors. Thus, $\rich{v, w}^{>0, \Lus}$ coincides with
	\[\rich{v, w}^{>0, \Ing}:=\{F \in \rich{v, w}: \text{all cluster variables in } \mca(\iseed{v, \bw}) \text{ are positive on }F \}.\]
	
	On the other hand, using \cref{thm:compareStructures}, $\rich{v, w}^{>0, \Lec}$ is equal to $\tau_{v,w}(\rich{v, w}^{>0, \Ing})= \tau_{v, w}(\rich{v, w}^{>0, \Lus})$. By \cite[Theorem 2.6, Remark 3.1]{GLRichardson} , the map $\tau_{v,w}$ preserves $\rich{v, w}^{>0, \Lus}$, so we are done.
\end{proof}

\bibliographystyle{alpha}
\bibliography{rich}

\newcommand{\etalchar}[1]{$^{#1}$}
\begin{thebibliography}{GLSBS22}

\bibitem[BFZ05]{CA3}
Arkady Berenstein, Sergey Fomin, and Andrei Zelevinsky.
\newblock Cluster algebras. {III}. {U}pper bounds and double {B}ruhat cells.
\newblock {\em Duke Math. J.}, 126:1--52, 2005.

\bibitem[BGW03]{CoxMatroid}
Alexandre~V. Borovik, I.~M. Gelfand, and Neil White.
\newblock {\em Coxeter matroids}, volume 216 of {\em Progress in Mathematics}.
\newblock Birkh\"{a}user Boston, Inc., Boston, MA, 2003.

\bibitem[BIRS09]{BIRS}
A.~B. Buan, O.~Iyama, I.~Reiten, and J.~Scott.
\newblock Cluster structures for 2-{C}alabi-{Y}au categories and unipotent
  groups.
\newblock {\em Compos. Math.}, 145(4):1035--1079, 2009.

\bibitem[CGG{\etalchar{+}}22]{CGGLSS}
Roger Casals, Eugene Gorsky, Mikhail Gorsky, Ian Le, Linhui Shen, and Jos\'e
  Simental.
\newblock Cluster structures on braid varieties.
\newblock {\em \tt{arXiv:2207.11607}}, 2022.

\bibitem[CGGS20]{CGGS}
Roger Casals, Eugene Gorsky, Mikhail Gorsky, and Jos\'e Simental.
\newblock Algebraic weaves and braid varieties.
\newblock {\em \tt{arXiv:2012.06931}}, 2020.

\bibitem[CK22]{CK}
Peigen Cao and Bernhard Keller.
\newblock On {L}eclerc's conjectural cluster structures for open {R}ichardson
  varieties.
\newblock {\em \tt{arXiv:2207.10184}}, 2022.

\bibitem[Deo85]{Deodhar}
Vinay~V. Deodhar.
\newblock On some geometric aspects of {B}ruhat orderings. {I}. {A} finer
  decomposition of {B}ruhat cells.
\newblock {\em Invent. Math.}, 79(3):499--511, 1985.

\bibitem[FG06]{FGDualityConj}
Vladimir Fock and Alexander Goncharov.
\newblock Cluster {$\mathcal X$}-varieties, amalgamation, and {P}oisson-{L}ie
  groups.
\newblock In {\em Algebraic geometry and number theory}, volume 253 of {\em
  Progr. Math.}, pages 27--68. Birkh\"{a}user Boston, Boston, MA, 2006.

\bibitem[FPST22]{FPST}
Sergey Fomin, Pavlo Pylyavskyy, Eugenii Shustin, and Dylan Thurston.
\newblock Morsifications and mutations.
\newblock {\em J. Lond. Math. Soc. (2)}, 105(4):2478--2554, 2022.

\bibitem[Fra16]{Fraser}
Chris Fraser.
\newblock Quasi-homomorphisms of cluster algebras.
\newblock {\em Adv. in Appl. Math.}, 81:40--77, 2016.

\bibitem[FSB22]{FSB}
Chris Fraser and Melissa Sherman-Bennett.
\newblock Positroid cluster structures from relabeled plabic graphs.
\newblock {\em Algebr. Comb.}, 5(3):469--513, 2022.

\bibitem[Ful97]{Fulton}
William Fulton.
\newblock {\em Young tableaux}, volume~35 of {\em London Mathematical Society
  Student Texts}.
\newblock Cambridge University Press, Cambridge, 1997.
\newblock With applications to representation theory and geometry.

\bibitem[FWZ16]{CAbook}
Sergey Fomin, Lauren Williams, and Andrei Zelevinsky.
\newblock Introduction to {C}luster {A}lgebras. {C}hapters 1-3.
\newblock {\em \tt{arXiv:1608.05735}}, 2016.

\bibitem[FZ99]{FZ99}
Sergey Fomin and Andrei Zelevinsky.
\newblock Double {B}ruhat cells and total positivity.
\newblock {\em J. Amer. Math. Soc.}, 12(2):335--380, 1999.

\bibitem[FZ02]{CA1}
Sergey Fomin and Andrei Zelevinsky.
\newblock Cluster algebras. {I}. {F}oundations.
\newblock {\em J. Amer. Math. Soc.}, 15:497--529, 2002.

\bibitem[GGS{\etalchar{+}}14]{GGSVV}
J.~Golden, A.~Goncharov, M.~Spradlin, C.~Vergu, and A.~Volovich.
\newblock Motivic amplitudes and cluster coordinates.
\newblock {\em Journal of High Energy Physics}, 2014(91), 2014.

\bibitem[GHKK18]{GHKK}
Mark Gross, Paul Hacking, Sean Keel, and Maxim Kontsevich.
\newblock Canonical bases for cluster algebras.
\newblock {\em J. Amer. Math. Soc.}, 31:497--608, 2018.

\bibitem[GL19]{GLPositroid}
Pavel Galashin and Thomas Lam.
\newblock Positroid varieties and cluster algebras.
\newblock {\em {\tt arXiv:1906.03501}}, 2019.

\bibitem[GL20]{GLCohom}
Pavel Galashin and Thomas Lam.
\newblock Positroids, knots, and $q,t$-{C}atalan numbers.
\newblock {\em \tt{arXiv:2012.09745}}, 2020.

\bibitem[GL22]{GLRichardson}
Pavel Galashin and Thomas Lam.
\newblock The twist for richardson varieties.
\newblock 2022.
\newblock arXiv:2204.05935.

\bibitem[GLS11]{GLS}
Christof Geiss, Bernard Leclerc, and Jan Schr\"{o}er.
\newblock Kac-{M}oody groups and cluster algebras.
\newblock {\em Adv. Math.}, 228:329--433, 2011.

\bibitem[GLSBS]{GLSBS2}
Pavel Galashin, Thomas Lam, Melissa Sherman-Bennett, and David~E Speyer.
\newblock Braid variety cluster structures, {II}: general type.
\newblock In preparation.

\bibitem[GLSBS22]{GLSBS1}
Pavel Galashin, Thomas Lam, Melissa Sherman-Bennett, and David~E Speyer.
\newblock Braid variety cluster structures, {I}: 3{D} plabic graphs.
\newblock {\em \tt{arXiv:2210.04778}}, 2022.

\bibitem[GSV10]{GSV}
Michael Gekhtman, Michael Shapiro, and Alek Vainshtein.
\newblock {\em Cluster algebras and {P}oisson geometry}, volume 167 of {\em
  Mathematical Surveys and Monographs}.
\newblock American Mathematical Society, Providence, RI, 2010.

\bibitem[Ing19]{gracie}
Grace Ingermanson.
\newblock Cluster algebras of open {R}ichardson varieties, 2019.
\newblock PhD thesis.

\bibitem[KL79]{KL79}
David Kazhdan and George Lusztig.
\newblock Representations of {C}oxeter groups and {H}ecke algebras.
\newblock {\em Invent. Math.}, 53(2):165--184, 1979.

\bibitem[KL80]{KL80}
David Kazhdan and George Lusztig.
\newblock Schubert varieties and {P}oincar\'{e} duality.
\newblock In {\em Geometry of the {L}aplace operator ({P}roc. {S}ympos. {P}ure
  {M}ath., {U}niv. {H}awaii, {H}onolulu, {H}awaii, 1979)}, Proc. Sympos. Pure
  Math., XXXVI, pages 185--203. Amer. Math. Soc., Providence, R.I., 1980.

\bibitem[KLR00]{KLR}
Christian Kassel, Alain Lascoux, and Christophe Reutenauer.
\newblock Factorizations in {S}chubert cells.
\newblock {\em Adv. Math.}, 150(1):1--35, 2000.

\bibitem[KLS13]{KLS}
Allen Knutson, Thomas Lam, and David~E. Speyer.
\newblock Positroid varieties: juggling and geometry.
\newblock {\em Compos. Math.}, 149:1710--1752, 2013.

\bibitem[Lec16]{Leclerc}
B.~Leclerc.
\newblock Cluster structures on strata of flag varieties.
\newblock {\em Adv. Math.}, 300:190--228, 2016.

\bibitem[Lus90]{LusCanonicalBases}
G.~Lusztig.
\newblock Canonical bases arising from quantized enveloping algebras.
\newblock {\em J. Amer. Math. Soc.}, 2:447--498, 1990.

\bibitem[Lus94]{LusTPGr}
G.~Lusztig.
\newblock Total positivity in reductive groups.
\newblock In {\em Lie theory and geometry}, volume 123 of {\em Progr. Math.},
  pages 531--568. Birkh\"{a}user Boston, Boston, MA, 1994.

\bibitem[M\'22]{Menard}
Etienne M\'enard.
\newblock Cluster algebras associated with open {R}ichardson varieties: an
  algorithm to compute initial seeds.
\newblock {\em \tt{arXiv:2201.10292}}, 2022.

\bibitem[MR04]{MR}
B.~J. Marsh and K.~Rietsch.
\newblock Parametrizations of flag varieties.
\newblock {\em Represent. Theory}, 8:212--242, 2004.

\bibitem[MS17]{MSTwist}
Greg Muller and David~E. Speyer.
\newblock The twist for positroid varieties.
\newblock {\em Proc. Lond. Math. Soc.}, 115:1014--1071, 2017.

\bibitem[Mul13]{Muller}
Greg Muller.
\newblock Locally acyclic cluster algebras.
\newblock {\em Adv. Math.}, 233:207--247, 2013.

\bibitem[Pos06]{Postnikov}
Alexander Postnikov.
\newblock Total positivity, {G}rassmannians, and networks.
\newblock {\em {\tt arXiv:math/0609764}}, 2006.

\bibitem[Rie99]{Rietsch}
K.~Rietsch.
\newblock An algebraic cell decomposition of the nonnegative part of a flag
  variety.
\newblock {\em J. Algebra}, 213:144--154, 1999.

\bibitem[Sco06]{Scott}
J.~S. Scott.
\newblock Grassmannians and cluster algebras.
\newblock {\em Proc. London Math. Soc.}, 92:345--380, 2006.

\bibitem[SSBW19]{SSBW}
Khrystyna Serhiyenko, Melissa Sherman-Bennett, and Lauren Williams.
\newblock Cluster structures in {S}chubert varieties in the {G}rassmannian.
\newblock {\em Proc. Lond. Math. Soc. (3)}, 119(6):1694--1744, 2019.

\bibitem[STWZ19]{STWZ}
Vivek Shende, David Treumann, Harold Williams, and Eric Zaslow.
\newblock Cluster varieties from {L}egendrian knots.
\newblock {\em Duke Math. J.}, 168(15):2801--2871, 2019.

\end{thebibliography}
	
\end{document}